\documentclass[a4paper,10pt]{scrartcl}

\setkomafont{title}{\sffamily\bfseries\huge}
\setkomafont{author}{\large}
\setkomafont{date}{\small}

\usepackage{lmodern}
\usepackage[T1]{fontenc}
\usepackage[left=2cm,right=2cm,top=2cm,bottom=3cm]{geometry}
\usepackage{framed,multirow}
\usepackage{microtype}
\usepackage[dvipsnames]{xcolor}
\usepackage{tikz}
\usepackage[hidelinks]{hyperref}
\microtypesetup{activate={true,nocompatibility},
                final,
                tracking=true,
                kerning=true,
                spacing=true,
                factor=1000,
                stretch=20,
                shrink=20}
\microtypecontext{spacing=nonfrench}

\usepackage{mathtools}
\usepackage{bm}
\usepackage{tabularx}
\usepackage{booktabs}

\graphicspath{{figures/}}

\usepackage{amsmath}
\usepackage{amssymb}
\usepackage{latexsym}
\usepackage{graphicx}  
\usepackage{subfig}
\usepackage{amsthm}
\usepackage{xspace}

\newcommand{\up}[1]{\ensuremath{\mathrm{#1}}}
\DeclarePairedDelimiter{\norm}{\lVert}{\rVert}
\DeclareMathOperator{\de}{d\!}
\renewcommand{\vec}[1]{\bm{\mathrm{#1}}}
\newcommand{\transpose}{\mathsf{T}}
\newcommand{\abs}[1]{\left|#1\right|}

\newcommand{\flll}{\vphantom{l^{l^l}}}
\newcommand{\mtot}{m_\up{tot}}
\newcommand{\xl}{x_\up{L}}
\newcommand{\xr}{x_\up{R}}
\newcommand{\dxi}{\Delta x_i}
\newcommand{\dmi}{\Delta m_i}

\newcommand{\dmiph}{\Delta m_{i+1/2}}

\newcommand{\miph}{m_{i+1/2}}
\newcommand{\imh}{{i-1/2}}
\newcommand{\ip}{{i+1}}
\newcommand{\im}{{i-1}}
\newcommand{\iph}{{i+1/2}}

\newcommand{\ek}{e_\up{k}}
\newcommand{\ei}{e}
\newcommand{\pinf}{\Pi}
\newcommand{\oftm}{(t,\ m)}

\newcommand{\dedp}{\frac{\partial \ei}{\partial p}}
\newcommand{\dedv}{\frac{\partial \ei}{\partial V}}
\newcommand{\aw}{a}
\newcommand{\awq}{a^2}
\newcommand{\solutionop}{\mathcal{H}}
\newcommand{\cfl}{k_\up{CFL}}
\newcommand{\cflzero}{k_\up{CFL}^0}
\newcommand{\muscl}{MUSCL--Hancock\xspace}
\newcommand{\rk}{Runge--Kutta\xspace}

   \newcolumntype{C}{>{\centering\let\newline\\\arraybackslash}X}
   \newcolumntype{R}{>{\raggedleft\let\newline\\\arraybackslash}X}
   \newcolumntype{L}{>{\raggedright\let\newline\\\arraybackslash}X}
   \newcolumntype{Z}{>{\raggedright\let\newline\\\arraybackslash}X}

\newcommand{\diffsl}[1]{S_{#1,\up{L}}}
\newcommand{\diffal}[1]{A_{#1,\up{L}}}
\newcommand{\diffsr}[1]{S_{#1,\up{R}}}
\newcommand{\diffar}[1]{A_{#1,\up{R}}}

\newcommand{\reviewera}[1]{\textcolor{black}{#1}}
\newcommand{\reviewerb}[1]{\textcolor{black}{#1}}

\usepackage{lineno}

\begin{document}

\begin{flushleft}
{\usekomafont{title}\linespread{1.08}\selectfont An efficient implicit scheme for the multimaterial Euler equations in Lagrangian coordinates\par}

\vspace{4ex}

{\usekomafont{author}Simone Chiocchetti\textsuperscript{1,2,3,*} and Giovanni Russo\textsuperscript{4}}

\vspace{3ex}

{\usekomafont{date}%
\textsuperscript{1}Division of Mathematics, University of Cologne, Weyertal 86-90, 50931 Cologne, Germany\\[0.3ex]
\textsuperscript{2}Institute of Aerodynamics and Gasdynamics, University of Stuttgart, Pfaffenwaldring 21, 70569 Stuttgart, Germany\\[0.3ex]
\textsuperscript{3}Laboratory of Applied Mathematics, Department of Civil, Environmental and Mechanical Engineering,\\ University of Trento, Via Mesiano 77, 38123 Trento, Italy\\[0.3ex]
\textsuperscript{4}Department of Mathematics and Informatics, University of Catania, Viale Andrea Doria 6, 95125 Catania, Italy\\[0.8ex]
\textsuperscript{*}Corresponding author. E-mail: simone.chiocchetti@uni-koeln.de%
}
\end{flushleft}

\vspace{1ex}

\begin{abstract}
\noindent\textsf{\textbf{Abstract.}} Stratified fluids composed of a sequence of alternate layers show interesting
macroscopic properties, which may be quite different from those of the
individual constituent fluids. On a macroscopic
scale, such systems can be considered a sort of fluid metamaterial.
In many cases each fluid layer can be described by Euler equations 
following the stiffened gas equation of state.
The computation of detailed numerical solutions of such stratified material poses several 
challenges, first and foremost the issue of artificial smearing of material parameters across 
interface boundaries. Lagrangian schemes completely eliminate this issue, but at the cost of rather
stringent time step restrictions. 
In this work we introduce an implicit numerical method for the multimaterial Euler equations in 
Lagrangian coordinates.
The implicit discretization is aimed at bypassing the prohibitive time step restrictions present 
in flows with stratified media, where one of the materials is particularly dense, or rigid (or both). 
This is the case for flows of water-air mixtures, air-granular media, or similar high density ratio systems.
We will present the novel discretisation approach, which makes extensive use of the remarkable structure
of the governing equations in Lagrangian coordinates to find the solution by means of a single implicit
discrete wave equation for the pressure field, yielding a symmetric positive definite structure and 
thus a particularly efficient algorithm.
Additionally, we will introduce simple filtering strategies for counteracting the 
emergence of pressure or density 
oscillations typically
encountered in multimaterial flows, and will present results concerning the robustness, 
accuracy, and performance of the proposed method, including applications to stratified media with 
high density and stiffness ratios.


\end{abstract}


\section{Introduction}\label{sec:introduction}
\paragraph{Motivation} Metamaterials have attracted an enormous attention in science and technology in the last two decades. 
They are composite materials obtained by assembling a large number of small unit cells, in 1D, 2D or 3D patterns, in such a way that 
on a length scale much larger than the cell size, they behave as a homogeneous material with
properties that are different from the original constituent materials. In most cases, they exhibit
properties that cannot be found in nature, and which are not even intermediate between the
(mechanical or optical) properties of the constituents. 

An enormous literature on metamaterials is available. Here we mention a couple of papers on acoustic metamaterials 
\cite{cummer2016controlling,torrent2007acoustic}, which are a particular type of mechanical
metamaterials. In these papers the authors emphasize the remarkable properties and the anomalous behaviour of
such materials. 
Among the surprising properties of acoustic metamaterials, we recall the so-called {\em negative
mass\/}. This means that the system behaves like a peculiar mass-spring system which, if excited
with a periodic external forcing term with a suitable frequency, shows a displacement which is
in phase with the force, rather than in phase opposition, as one would expect
\cite{ma2016acoustic}. 
In \cite{lee2016origin} it is shown how to create a simple 1D acoustic metamaterial formed by a gas
tube, with compartments separated by elastic membranes. 

The analysis of the properties of acoustic metamaterials is in general based on linear or linearized
equations, so that it is assumed that the displacement, or any type of signal, is sufficiently small
so that linear theory can be used, or that non-linear effects can be computed by perturbation
methods, still assuming sufficiently small signals. We also point out that fully nonlinear innovative 
theories of metamaterials are object of current research, see for example \cite{peshkov2019}.

There are several cases in which non-linear effects are crucial and cannot be neglected. 
In a multilayer system formed by several pairs of layers of two different fluids, for example, the
non-linear response of the gases cannot be neglected, and shocks will in general form and propagate
in the multilayer medium, unless the initial state or the forcing terms are suitably small, or the periodicity
of the system indices an oscillator behaviour of the solutions, as it has been recently observed in several quasilinear hyperbolic systems when the waves propagates over a periodic medium, as in the case of non linear elasticity
in 1D \cite{leveque2003} and in 2D \cite{quezada_dispersion}
or shallow water over periodic bathymetry in 1D 
\cite{KLR2023} or in 2D  \cite{2021_solitary}. In all these cases, the amplitude of the wave is large enough to observe non-linear effects, still sufficiently small to avoid formation of shock waves. 
For large amplitude waves propagating on layered media, shocks will generally form. 

In \cite{phan2023numerical}, for example, the behaviour of a large number of pairs of layers is
studied numerically. Each pair is formed by two different fluids, each one treated as a stiffened
gas with its own standard density $\rho$, adiabatic exponent $\gamma$, and a constant $\Pi$ which
determines the stiffness of the fluid, and which is sometimes called attractive pressure
\cite{Chiapolino_Saurel_18}. 

In the paper it is shown that the a simple adiabatic homogeneous model is able to accurately
reproduce the behaviour of the multilayer system provided isentropic initial conditions are
assigned, and the integration time is short enough (or equivalently the deviation from global
equilibrium is small enough) so that no strong shocks form. When shocks form, the predictions of the
homogeneous model are no longer reliable, hence the use of a more accurate model is necessary. 

Another interesting feature of the multilayer system is that, depending on the fluid parameters, the
effective sound speed of the homogenised fluid may be a non monotonic function of the mass fraction
of the layers. As a result, the effective sound speed of the layer may be smaller than the sound
speed of both gases. 

\reviewerb{The connection itself of homogenised multiphase flow models to the detailed description of multilayered systems with sharp interfaces
coupling two or more fluids governed by the compressible Euler or Navier--Stokes equations is an additional point of interest.
For example, Bresch, Burtea and Lagoutière \cite{bresch2022continuous, bresch2023semidiscrete} discuss the link between such layered 
systems and multiphase flow models of Baer--Nunziato type. 
They provide formal proof, starting from a layered two-fluid system governed by the compressible Navier--Stokes equations 
    (which give no-slip coupling conditions between interface velocities), and 
    from a semi-discrete description of the same system, that homogenisation will yield a viscous, compressible, single-velocity, 
    two-pressure multiphase model of the Baer--Nunziato \cite{baer1986} flavour (with single velocity and two pressures as in \cite{saurel2009, pelanti2014}). 
    They also formally recover Kapila’s inviscid single-velocity, single-pressure model \cite{kapila2001} in the limit of vanishing viscosity.}
    
    \reviewerb{In this regard, in this paper we limit our scope to providing experimental evidence of the latter remark that multilayer systems governed by the compressible Euler equations yield 
    a homogenised behaviour that is captured at a coarse level by Kapila's model. The strategy used here is just incidental: homogenisation 
    effects are observed/captured by the scheme automatically, when space and time are under-resolved with respect to the natural fast/small scales of the detailed Euler--Euler flow.}
    
\reviewera{Tackling multiphase flows starting from the detailed description of the underlying two-fluid system was also the strategy 
adopted by Petrella, Abgrall, and Mishra in \cite{petrella2023montecarlo}, where the closure problem of multiphase flow models is avoided 
entirely by carrying out a set of front tracking simulations, embedded in a Monte Carlo loop, 
which yields solutions of the homogenised model as a result of the accumulated statistics. 
Like in this paper, the authors adopted a sharp description of material fronts, even if afterwards smoothing is obtained via the Monte Carlo averaging. Similarly, here we capture material interfaces sharply and then we can simply under-resolve
the timestepping to compare with the solutions of the homogenised models.}

\paragraph{Numerical methodology}
In \cite{phan2023numerical} the detailed numerical solution of the multilayer system was obtained by
a finite volume methods for solving the Euler equations that describe the fluids expressed in
Lagrangian form, where the mass coordinate and time are adopted as independent variables. This
choice appears to be a natural one, since it allows to have material boundaries between the two
fluids always located at cell edges, so strictly speaking in this configuration 
there is no need to use a diffuse interface multiphase model, such as for example Kapila's \cite{kapila2001} 
reduced Baer--Nunziato family model \cite{baer1986, pelanti2014, saurel2009}, 
or \cite{saurel1999b, RomenskiTwoPhase2007, RomenskiTwoPhase2010, romenski2016, thein2022, delrio2024}, 
and each cell completely belongs to either
the first or the second fluid. A positive side effect of such a separation is that very sharp material
interfaces and contacts are naturally obtained, directly at cell interfaces (which is very convenient 
in one space dimension)
or from a level set function \cite{sussman1994} in the more general case. 
Furthermore, working with ``pure'' states makes it
simpler to perform reconstructions and to design exact or approximate Riemann solvers. 
As mentioned, versatile diffuse interface treatments for this problem have shown promising results
for solid-fluid interfaces \cite{favrie2009, kemm2020} and the treatment of free surfaces \cite{dumbser2011b, gaburro2018, ferrari2024}.



Explicit schemes for the numerical solution of fluid equations suffer from classical
Courant-Friedrichs-Lewy time step restriction \cite{courant1928}, which can be written as $\Delta
t<C \Delta \xi/\lambda$, where $\lambda$ denotes the maximum characteristic speed of the system in
the computational domain, $C$ is a constant of order of magnitude 1, and $\Delta \xi$ denotes the
discretization of the independent space coordinate (which is typically space $x$ in Eulerian
coordinates and mass coordinate $m$ in the Lagrangian framework). 
The ability of constructing versatile and efficient high order massively 
parallel codes \cite{reinarz2020, kurz2024} 
is however a very attractive feature of such schemes.

Explicit integration also is sensitive, in Lagrangian an 
arbitrary Lagrangian-Eulerian \cite{hirt1974, duarte2004, berndt2011, barlow2016, loubere2010, ShashkovRemap5}, 
or \cite{dumbser2014, boscheri2015, gaburro2020}, 
to the formation of small cells that yield very stringent timestep restrictions, 
though this can be mitigated, for example, by using local timestepping methods \cite{dumbser2014}.

If sound speed in one of the two materials is much larger than in the other, say
$\lambda_1\gg\lambda_2$, then the most severe time restriction on the time step will be due to the
waves propagating in fluid 1. On the other hand, we expect that the faster material behaves almost
like an incompressible fluid compared to the other, therefore most of the effect of compressibility
will be evident in the softer fluid, for which stability and accuracy requirements on the time step
will be similar. For such a reason the system becomes {\em stiff\/}: stability requirements on the
time step in the first fluid are much more restrictive than accuracy requirements. 

Some correction strategies for explicit methods, made to deal with inaccuracies encountered 
by such schemes in the low-Mach regime, are also available \cite{barsukow2017, bruel2019, rieper2011}, 
but it would seem that an implicit or semi-implicit strategy is required to deal with excessive timestep 
restrictions characteristic of low Mach number flows.
In general, the use of some implicitness 
may mitigate the stiffness problem, thus improving the overall efficiency of the scheme. 

An extreme case of multifluid in which the densities and stiffness of the two fluids are very
different is a multilayer formed by a sequence of air and water layers. Here the ratio of the
densities is close to one thousand, while the ratio in the sound speed, in standard Eulerian coordinates,  is approximately four (see
Section \ref{sec:speeds}). 
In \cite{phan2023numerical} only explicit schemes with locally uniform spacing in the mass
Lagrangian coordinate have been adopted. The ratio of the sound speed between the two fluids in
Lagrangian coordinates is equal to the Eulerian one multiplied by the ratio in the density, which
causes the ratio of Lagrangian sound speed between water and air to be more than three thousand, making it
almost unfeasible to adopt an explicit scheme in Lagrangian coordinates with uniform mesh for the
numerical treatment of water-air multilayer systems. For this reason in \cite{phan2023numerical} the
authors considered only a moderate density ratio (either 2 or 10 in particular). 

A semi-implicit approach would indeed be a very natural choice for an efficient implicit FV scheme:
solving convection explicitly and treating acoustic waves implicitly, recovering pressure
and momentum, allows to construct very simple and efficient numerical discretizations \cite{munz2003, park2005, casullidumbser}. 
This approach has been widely adopted in the Eulerian framework to construct the so
called all-Mach number solvers for compressible Euler equations (see for example
\cite{degond2011all, cordier2012} for one of the first second order finite volume all Mach solvers for isentropic
equations, \cite{boscarino2018all} for a second order finite volume all-Mach solver for full Euler
equations on staggered grid, or \cite{avgerinos2019linearly} for a second order FV full Euler
all-Mach solver on unstaggered grid.) Second order accurate all Mach solvers for full 3D Euler
equations have been constructed and tested in \cite{boscheri2020second}, while third order finite
difference all Mach solvers for Euler equations have been developed in \cite{boscarino2019high} for
isentropic Euler equations, and in \cite{boscarino2022high} for full Euler equations. Using 
staggered grids, similar all-speed methods have been developed also in a hybrid Finite Volume-Finite Element 
framework \cite{busto2021}. Finally, an
all-Mach solver for isentropic two phase flows has been developed in \cite{lukavcova2022all} (see also \cite{lukakova2024}). 

A common feature of all such solvers is that only the fast waves are treated implicitly, 
and in most cases only linearly implicit, i.e.\ the implicit solver is applied to 
a linearization of the equations. High order in time is usually achieved using 
Implicit-Explicit time solvers, such as Runge-Kutta IMEX (see, for 
example \cite{ascher1997implicit}, or \cite{pareschi2005implicit}). A general 
technique to obtain high order accuracy when the system containing stiff and 
non-stiff terms is not written in additive nor in partitioned for is 
introduced in \cite{boscarino2016high}, and adopted, 
for example, in \cite{boscarino2019high} and \cite{boscarino2022high}.
It is worth noting that the versatility of such non-partitioned IMEX schemes
allows their use in complex frameworks like in \cite{boscheri2022}, where 
IMEX integration is used in conjunction with a black box solver derived from \cite{chiocchettimueller}.

Semi-implicit schemes have the advantage of being usually more efficient 
than fully implicit ones, 
for the same required accuracy. 
Nevertheless, very efficient implicit methods like \cite{bassi2015} are available, and 
implicit methods for fluid flow, even of Finite Volume type \cite{puppo2023} are currently
being actively developed in the literature.

Moreover, semi-implicit schemes tend to be 
less stable, and have been found to be more prone to spurious oscillations (see, for 
example, \cite{boscarino2015linearly}, where semi-implicit and fully 
implicit solvers are compared for a class of convection-diffusion problems). 
Nevertheless, a vast amount of successful schemes for low Mach and all Mach flows have been constructed
with this methodology, and in particular \cite{casullidumbser}, 
which provided the main guidelines
for the efficient formulation here introduced.

In this work, we instead aim at constructing a fully-implicit method along the steps of the recent 
work by Plessier et al. \cite{plessierdespres1} (see also \cite{delpino2024}), which is the main motivation for the present paper, 
with a different discretization strategy that leverages
the structure of the governing equations to build an efficient method that resembles a Lagrangian adaptation
of the scheme by Dumbser and Casulli \cite{casullidumbser}. Additionally, here we introduce 
several strategies for controlling spurious oscillations arising in the solution of 
two-material compressible flows and apply the new scheme to multi-layer systems following the
work by Phan et. al. \cite{phan2023numerical}, where explicit numerical methods were used.













The system of Euler equations for stiffened gas in Lagrangian coordinates offers a very natural opportunity for distilling a 
single scalar discrete wave equation for the pressure field, which can be solved iteratively by a
sequence of tridiagonal systems (predictor), leading to a preliminary (possibly oscillatory and non
conservative) solution which can be made conservative and non-oscillatory with a suitable
post-processing (corrector). Once a first order scheme is produced, singly-diagonally-implicit
Runge-Kutta (SDIRK) schemes can be adopted to rise time accuracy to second or higher order. 
An additional advantage of the approach is that it does not require sophisticated Riemann solvers 
to reproduce sharp contact and material interfaces. 

The method naturally generalizes to the treatment of layered multifluid, for which, as we already 
mentioned, the Lagrangian approach avoids the use of intermediate states. When dealing with material 
interfaces, the methods blends the use of conservative and primitive variables in the time 
advancement of the solution, using ideas from \cite{abgrall2001}. 
Strict conservation is locally lost, yet, the shock speed is correctly captured up to 
the accuracy of the method (see Section~\ref{sec:rp}).

Last, but not least, a whole section (Sec.\ \ref{sec:mesh}) is devoted to the construction 
of quasi-uniform grids, i.e.\ grids in the Lagrangian mass variable whose spacing varies smoothly 
among cells. This avoids the drawbacks of uniform cells when dealing with fluids with large mass 
density ratio: using uniform cells in the mass variable under-resolve the regions occupied by the 
low density fluid, or over-resolves the regions occupied by the high density fluid. On the other 
hand, using grid spacing which is approximately uniform in space causes an abrupt jump in the 
grid spacing in the mass variable, which may induce oscillations and sometimes unphysical solutions.

The plan of the paper is the following. After the introduction, 
Section \ref{sec:equation} describes the model under consideration, i.e.\ Euler system for 
a stiffened gas in Lagrangian coordinates. 
Section \ref{sec:method} is the core of the paper and describes in detail the construction 
of the method as outlined in the previous paragraphs (non-conservative, possibly oscillatory, 
fully implicit predictor, conservative non-oscillatory corrector). 
Section \ref{sec:mesh} deals with the construction of mass conservative well-graded 
Lagrangian meshes. 
Section \ref{sec:results} is devoted to illustrate the performance of the method on a 
sequence of classical tests for a single gas and for a multilayer systems. 
\reviewerb{In particular, Section \ref{sec:multilayer} is dedicated to applications to systems with many 
stratified layers, and compares and contrast the behaviour of the scheme with that obtained from the a homogenized multiphase model.} Finally, 
in Section \ref{sec:conclusions} we draw some conclusions and perspective about applications and extensions of the method. 






\section{Governing equations}\label{sec:equation}

We consider the one-dimensional Euler equations in Lagrangian (mass) coordinates.
The spatial mass abscissa $m$ (also referred to as Lagrangian abscissa in this manuscript) 
is defined at a point $x$ (in this manuscript also called
Eulerian coordinate) as the mass of the fluid (per unit area)
found between the left boundary ($m = 0$) and the Eulerian coordinate $x$, or formally
\begin{equation} \label{eqn:masscoord}
    m(x) = \int_{\xl}^{x} \rho(x^\prime)\,\de{x^\prime},
\end{equation}
where $\rho$ is the mass density of the fluid.
Symmetrically, the Eulerian coordinate $x$ can be recovered from the mass coordinate
as 
\begin{equation} \label{eqn:eulercoord}
    x(m) = \int_{0}^m V(m^\prime)\de{m^\prime},
\end{equation}
with $V = 1/\rho$ the specific volume of the fluid.

In the Lagrangian mass coordinate system, the conservation laws for mass, momentum, and energy read
\begin{subequations}
    \begin{align}
        &\partial_t \left(V\right) + \partial_m\left(-u\right) = 0, \label{eqn:mass}\\[1mm]
        &\partial_t \left(u\right) + \partial_m\left(p\right) = 0, \label{eqn:momentum}\\[1mm]
        &\partial_t \left(E\right) + \partial_m\left(u\,p\right) = 0, \label{eqn:energy}\\[1mm]
        &\partial_t \left(c\right) = 0, \label{eqn:colour}
    \end{align}
\end{subequations}
where we define $V = V\oftm$ to be the specific volume of the fluid (the reciprocal of the mass density $\rho$), 
the velocity field is $u = u\oftm$, and $E\oftm = \ek + \ei$ is the specific total energy density.
The system is closed with the definition of the specific kinetic energy $\ek = u^2/2$ 
and by providing an equation of state (EOS) that expresses the internal energy $\ei$ as 
a function of density (or specific volume) and pressure.
In this work we adopt the stiffened gas equation of state in the form 
\begin{equation}
    \label{eqn:eos}
    \ei = V\,\frac{p + \gamma\,\pinf}{\gamma - 1}.
\end{equation}
Note that, since we allow the presence of two separate (immiscible) fluids, the local 
adiabatic index $\gamma$ and the stiffness parameter $\pinf$ may vary between
two values $\gamma_1$ and $\gamma_2$, and $\pinf_1$ and $\pinf_2$. Owing to the immiscible
nature of the modelled fluids (which is maintained also at the discrete level), 
there is no necessity for mixture rules defining homogenised or averaged parameters
at interfaces. Nevertheless, as a piece of notational convenience, we can 
define a colour function $c$ which helps tracking the presence of one or the other fluid 
at a given Lagrangian location. The value of $c$ is initially assigned cell-by-cell in 
as either $c=1$ (meaning only fluid 1 is present) or $c=0$ (only fluid 2 is present).
This means that the parameters of the equation of state can be locally recovered 
by computing $\gamma = c\,\gamma_1 + (1 - c)\,\gamma_2$ and $\pinf = c\,\pinf_1 + (1 - c)\,\pinf_2$.
We should emphasize that such a notation is \textit{not} a mixture rule (or at least not a good one)
but rather a shorthand which is valid in this context only thanks to the fact that $c$ is allowed only 
to be valued \textit{exactly} $c = 0$ or \textit{exactly} $c = 1$. Moreover, in Lagrangian coordinates, 
Equation \eqref{eqn:colour} is not discretised at all, since it simply states that the EOS parameters
are immutably associated with each discrete cell.

See \cite{abgrall2001} for a discussion
of computationally favourable choices for averaged equation of state parameters.

\subsection{Characteristic speeds in Lagrangian coordinates}
\label{sec:speeds}
We briefly discuss the eigenvalues of the hyperbolic system of PDEs \eqref{eqn:mass}--\eqref{eqn:energy}, 
in regards to the CFL timestep restriction \cite{courant1928} they obey to when explicit timestepping
is adopted.
Having defined a vector of conserved variables $\vec{Q} = {\left(V,\ u,\ E\right)}^\transpose$, 
and a flux vector $\vec{F} = {\left(-u,\ p,\ u\,p\right)}^\transpose$, the system of conservation laws 
\eqref{eqn:mass}--\eqref{eqn:energy} can be cast as
\begin{equation}
    \partial_t \vec{Q} + \partial_m \vec{F} = \vec{0}, 
\end{equation}
or, in terms of the primitive variables $\vec{W} = {\left(V,\ u,\ p\right)}^\transpose$ and 
in its so-called quasi-linear form
\begin{equation}
    \partial_t \vec{W} + \vec{A}\,\partial_m \vec{W} = \vec{0},\quad \text{with}\quad \vec{A} = {\left(\frac{\partial{\vec{Q}}}{\partial{\vec{W}}}\right)}^{-1}\,
    \left(\frac{\partial \vec{F}}{\partial \vec{W}}\right).
\end{equation}
The system matrix is 
\begin{equation}
    \label{eq:A_Lag}
    \vec{A} = \begin{pmatrix*}[c]
    0 & -1 & 0\\
    0 & 0 & 1\\
    0 & a^2 & 0\\
    \end{pmatrix*},\quad \text{with}\quad a^2 = {\left({\dedp}\right)}^{-1}{\left(\dedv + p\right)} = \frac{\gamma\,\left(p + \pinf\right)}{V} = \rho\,\gamma\,\left(p + \pinf\right)
\end{equation}
and has eigenvalues $\vec{\lambda} = {\left(-a,\ 0,\ a\right)}^\transpose$.
Crucially, in mass coordinates, the square of the acoustic eigenvalues of the system are proportional to the mass density 
$\rho$. Since for explicit schemes the maximum stable timestep \cite{courant1928} is inversely proportional to 
the magnitude of the eigenvalues, such a dependence on $\rho$ implies that denser fluids will 
be more restrictive than lighter ones with regards to the maximum admissible timestep of explicit methods, if we assume a uniform mesh in the mass variable.
In contrast, the same analysis, when carried out for the Eulerian variant of the governing 
equations, would yield $\lambda_\up{E} = {\left(u-a_\up{E},\ u,\ u+a_\up{E}\right)}^\transpose$, 
with $a_\up{E}^2 = \gamma\,\left(p + \pinf\right)/\rho$, meaning that in the Eulerian 
framework denser fluids admit larger timesteps
for a fixed choice of the EOS parameters $\gamma$ and $\pinf$, assuming uniform mesh in the space variable.

Oftentimes, in simulations involving both a non-stiff fluid and a stiff one, (an example might be
air and water), 
the fluid characterised by higher values of $\gamma$ and $\pinf$ (that is, the stiff one), is also the 
denser of the two.
In the following paragraphs, we give a practical quantitative illustration of the discrepancy: at standard 
atmospheric conditions $p=10^5\,\up{Pa}$, 
$\gamma_\up{a} = 1.4$, 
$\pinf_\up{a} = 0\,\up{Pa}$, 
$\rho_\up{a} = 1.2\,\up{kg\,m^{-3}}$, 
an estimate for the speed of sound is $a_\up{E}^\up{a} \simeq 341.6\,\up{m\,s^{-1}}$, 
while at the same pressure the speed of sound for water can be estimated
to be $a_\up{E}^\up{w} \simeq 1482.3\,\up{m\,s^{-1}}$, with $\gamma_\up{w} = 7.3$, 
$\pinf_\up{w} = 3.0\times10^8\,\up{Pa}$, 
$\rho_\up{a} = 997\,\up{kg\,m^{-3}}$.
The same choice of material parameters, when adopted in the Lagrangian framework, 
yields $a^\up{a} \simeq 409.9\,\up{kg\,m^{-2}\,s^{-1}}$
and $a^\up{w} \simeq 1.478\times10^6\,\up{kg\,m^{-2}\,s^{-1}}$ (the units change due to 
the fact that the spatial coordinate $x$ is replaced by the mass coordinate $m$), 
meaning that in mass coordinates 
the water phase imposes a timestep more than 3600 times smaller than the one allowed by 
the air on the same mesh spacing, as opposed to the ratio being about 4.3 in the Eulerian
setting. 

Of course the physical speed of sound waves remains the same in both cases, but from the computational
standpoint, the timestep size of explicit Lagrangian schemes using mass coordinates is 
restricted by a much higher numerical signal speed.


\section{Numerical method}\label{sec:method}

\reviewerb{In this Section we present the scheme here developed and used, 
which is based on a Lagrangian staggered grid in mass coordinates, and makes use 
of a predictor step involving a smooth wave-equation discretisation (see Section \ref{sec:discwave}), with a finite difference 
type correction step. To manage spurious oscillations that may be generated as a result of these low-dissipation discretisation 
choices (Lagrangian coordinates, smooth predictor step, central differences), we also introduce ad hoc
filtering techniques that have been embedded in the iterative computation of the solution. These target 
oscillations of specific volume in Section \ref{sec:densityfiltering} and oscillations of pressure in Section \ref{sec:pressurefiltering}.
Leaving such spurious oscillations unchecked might lead to catastrophic failure of the computations when strong shockwaves are encountered, due to violation
of positivity in the pressure and in the specific volume, hence implementation of some form of oscillation control is crucial to the 
robustness of the scheme. Figure \ref{fig:issues} provides some visual examples of how spurious oscillations might manifest
even in mild problems, due to the intrinsic low dissipation at material interfaces given by the Lagrangian framework, as well as due to the presence
of non-homogeneous material interfaces. The reader is referred to the paper by Abgrall and Karni \cite{abgrall2001} for an extensive discussion of these issues in multifluid computation.}

\subsection{Staggered grid}\label{sec:grid}

The scheme developed in this work adopts a staggered grid, mainly with the goal of obtaining
a well-behaved and inexpensive discretisation for a wave equation yielding
a preliminary solution for the pressure field, whereby the pressure field $p_i$, the specific volume $V_i$ and the specific total energy $E_i$
can be seen as collocated at cell centers, (the conserved quantities can equivalently be interpreted as cell average values), and the velocity field $u_{i+1/2}$
is collocated at cell interfaces.
The left boundary is denoted with the Eulerian coordinate $\xl$, and by definition 
the corresponding mass coordinate is $m(\xl) = 0$, while the right boundary is $\xr$, 
and $m(\xr) = \mtot$ is the total mass of the system.

The computational domain is partitioned using a staggered grid.
The main grid is composed of
$N$ cells of index $i = 1,\ 2,\ \hdots,\ N$, having variable width $\dxi$ 
or mass content $\dmi$.
The cell boundaries are formally denoted $\miph$ regardless of the uniformity of the grid, 
and define a face-based dual grid, which is used for collocation of the velocity field $u$.
For notational convenience, we define at each cell boundary the mesh spacing for the dual 
grid as $\dmiph = \left(\Delta m_i + \Delta m_{i+1}\right)/2$, which 
corresponds to the distance between the cell centers of the two adjacent control volumes.


\subsection{Discretisation approach}\label{sec:discretisation}

\subsubsection{Wave equation}\label{sec:discwave}
The combination of the governing equations \eqref{eqn:mass}, \eqref{eqn:momentum}, and \eqref{eqn:energy}
into a single wave equation constitutes the core of the numerical method proposed in this work.
Although the derivation of such a wave equation is a trivial exercise, we carry out here the main
steps in order to clarify the scope of applicability of the proposed framework.

An important preliminary note is that the employed wave equation relies on the validity
of the strong form of the governing equations, implying that
density, velocity, and pressure, should be smooth fields, rather than discontinuous as they will 
instead be in practice. Nonetheless, it will be shown with numerical examples that this fact 
does not hinder the applicability of the numerical method here presented.

Further corrections to the discretisation scheme, accounting for the restricted validity of this 
assumption will be detailed in the subsequent paragraphs.

The derivation begins by applying the chain rule to Equation \eqref{eqn:energy} to obtain

\begin{equation}
    \partial_t\left(\ei + u^2/2\right) 
    + \partial_m(up) = \dedp\,\partial_t p + \dedv\,\partial_t V + u\,\partial_t u + u\,\partial_m p  + p\,\partial_m u = 0, 
\end{equation}
which, due to the momentum equation \eqref{eqn:momentum} simplifies to 
\begin{equation}\label{eqn:wave2}
    \dedp\,\partial_t p + \dedv\,\partial_t V  + p\,\partial_m u = 0.
\end{equation}
Inserting the mass equation \eqref{eqn:mass}
into \eqref{eqn:wave2}, and making use of Eq.\eqref{eq:A_Lag}, after simple algebraic rearrangement, 
the wave equation reads
\begin{equation} \label{eqn:wave}
    \partial_t p + \aw^2\,\partial_m u = 0.
\end{equation}
In the case of the stiffened gas EOS \eqref{eqn:eos}, the explicit expressions of the partial derivatives of the internal energy are
\begin{equation}
    \dedp = \frac{V}{\gamma - 1}, \qquad \dedv = \frac{p + \gamma\,\pinf}{\gamma - 1}
\end{equation}
and thus the wavespeed is $\aw = \sqrt{\gamma\,\left(p + \pinf\right)/V}$.
Thanks to the Lagrangian coordinates, the wave equation assumes the very simple
form \eqref{eqn:wave} and, most important, it does not feature any nonlinear convective terms, which would be present in 
its Eulerian analogue.

\subsubsection{Discrete wave equation}
The continuous wave equation \eqref{eqn:wave} can be used as the basis for the construction
of a \textit{discrete} system of algebraic equations which can be solved rather efficiently as a 
sequence of implicit linear problems for the pressure field only.

Discrete values $p_i^n$, $p_i^{n+1}$ (cell $i$, time level $t = t^n$ and $t = t^{n+1}$ respectively) 
for the pressure field are collocated at the cell centers of the main grid ($m = m_i$). Consequently 
a rather natural finite-difference discretisation of \eqref{eqn:wave} is
\begin{equation} \label{eqn:wavedisc0}
    \frac{p_i^{n+1} - p_i^n}{\Delta t} + {(\awq)}_i^{n+1}\,\frac{u_\iph^{n+1} - u_\imh^{n+1}}{\Delta m_i} = 0.
\end{equation}
We remark that conveniently the computation of the wavespeed ${(\awq)}_i^{n+1}$ is required \textit{only} at cell centers, 
where the discrete pressure field is directly collocated, requiring no interpolation of pressure
or internal energy, and guaranteeing that no artificial mixture-states have to be generated, 
as would be the case if $\awq$ had to be evaluated at cell interfaces. This is a welcome feature
since simple averaging rules might generate non-physical, highly inaccurate, or oscillatory results (see 
for example \cite{abgrall2001}).

Following the strategy used for the fully discrete semi-implicit methodology 
introduced in \cite{casullidumbser}, 
we introduce a finite difference discretisation for the velocity field update 
at the grid interface location
\begin{equation} \label{eqn:velocityupdate}
    u_{\iph}^{n+1} = u_\iph^n - \frac{\Delta t}{\Delta m_{i+1/2}}\,\left(p_\ip^{n+1} - p_\im^{n+1}\right), 
\end{equation}
which can then be inserted in \eqref{eqn:wavedisc0}
giving
\begin{equation}
    \label{eqn:nlsystem}
    \frac{p_i^{n+1} - p_i^n}{\Delta t} + \frac{(a^2)_i^{n+1}}{\Delta m_i}\,
    \left[
        u_{i+1/2}^n - \frac{\Delta t}{\Delta m_{i+1/2}}\,\left(p_{i+1}^{n+1} - p_i^{n+1}\right)
        -u_{i-1/2}^n + \frac{\Delta t}{\Delta m_{i-1/2}}\,\left(p_{i}^{n+1} - p_{i-1}^{n+1}\right)
    \right] = 0,
\end{equation}
The nonlinear system of algebraic equations   \eqref{eqn:nlsystem} for $p^{n+1}$ is solved via a simple
fixed point iteration which generates a sequence of 
approximations $\widetilde{(a^2)}_i^{n+1,r} = a^2(V^{n+1,r-1}_i,\ p^{n+1,r-1}_i)$
for the values of the nonlinear function $a^2(V,\ p)$,
converging towards the value at the new time level $(a^2)_i^{n+1} = a^2(V^{n+1}_i,\ p^{n+1}_i)$.

As an initial guess for the first iteration $r = 1$, we use the pressure at the 
previous time level $p_i^{n+1,0} = p_i^n$, 
while the specific volume is obtained as
\begin{equation}
    V_i^{n+1,0} = V_i^n - \frac{\Delta t}{\Delta m_i}\,\left(f^V_\iph - f^V_\imh\right),
\end{equation}
with the numerical flux $f^V_\iph$ given by be the simple Rusanov \cite{rusanov1961} approximate Riemann solver 
(with mild, flow velocity-based, signal speed estimates) as
\begin{equation} \label{eqn:rusanovguess}
    f^V_\iph = \frac{1}{2}\,\left(-u_i - u_{i+1}\right) - \frac{1}{2}\,\max\left(\rho_i\,\abs{u_i},\ \rho_{i+1}\,\abs{u_{i+1}}\right)\,\left(V_{i+1} - V_{i}\right).
\end{equation}
The scheme is not particularly sensitive to this initial guess, as evidenced by the fact that the sharp
contact discontinuities typical of Lagrangian schemes employing complete Riemann solvers are recovered
by the proposed scheme despite not using such contact-resolving Riemann solvers at all.
Simpler initial guesses, such as $V_i^{n+1,0} = V_i^n$ can be considered, with no effect on the accuracy or 
robustness of the scheme. The one given in \eqref{eqn:rusanovguess} 
is the one used for all 
tests included in this paper, out of an abundance of caution, 
but no solid arguments for why it should be preferred to simpler (or trivial) initial
guesses have been found in this work.



The linear system presents a tridiagonal structure
\begin{equation} \label{eqn:linearpressure}
k_{i-1}\,p_{i-1}^{n+1} + 
k_i\,p_{i}^{n+1} + 
k_{i+1}\,p_{i+1}^{n+1} = b_i
\end{equation}
with the coefficients being
\begin{equation}
\begin{aligned}
    & k_{i-1} = -\widetilde{(a^2)}_i^{n+1}\,\frac{\Delta t}{\Delta m_i}\,\frac{\Delta t}{\Delta m_{i-1/2}}, \\
    & k_i = 1 + \widetilde{(a^2)}_i^{n+1}\,\frac{\Delta t}{\Delta m_i}\,\left(\frac{\Delta t}{\Delta m_{i-1/2}}+\frac{\Delta t}{\Delta m_{i+1/2}}\right),\\
    & k_{i+1} = -\widetilde{(a^2)}_i^{n+1}\,\frac{\Delta t}{\Delta m_i}\,\frac{\Delta t}{\Delta m_{i+1/2}},
\end{aligned}
\end{equation}
and the known right hand side
\begin{equation}
    b_i = p_i^n - \widetilde{(a^2)}_i^{n+1}\,\frac{\Delta t}{\Delta m_i}\,\left(u_{i+1/2}^n - u_{i-1/2}^n\right),
\end{equation}
where the iteration index $r$ for the pressure $p_i^{n+1,r}$ and the wavespeed $\widetilde{(a^2)}_i^{n+1,r}$ 
have been omitted for the sake of brevity.

As soon as the discrete pressure field $p_i^{n+1}$ at cell centers is recovered by solving the linear system \eqref{eqn:linearpressure}, 
the velocity at the cell boundary $u_{i+1/2}^{n+1}$ can also immediately be computed using the central difference \eqref{eqn:velocityupdate}, 
and subsequently the specific volume is updated in an analogous manner
\begin{equation} \label{eqn:volumeupdate}
    (V_i^{n+1})_\up{c} = V_i^n - \frac{\Delta t}{\Delta m_i}\,\left(-u_{i+1/2} + u_{i-1/2}\right).
\end{equation}

The quantity given by \eqref{eqn:volumeupdate} is denoted $(V_i^{n+1})_\up{c}$ to indicate that
it is obtained from the direct discretisation of the PDE via a conservative update formula, based on 
a staggered central difference. The updated value of the specific volume $V_i^{n+1}$ is computed
after having found a suitable set of diffusive fluxes, in order to introduce numerical dissipation and control the 
oscillations generated by the central
update formula.
The determination of such fluxes is detailed in the next subsection.

Noted that the state variables have to be updated at each iteration in a rather nonlinear manner, 
inclusive of numerical stabilisation via diffusion and data-dependent filtering.  
The iteration loop is split into two nested parts: a first one (an inner loop) in which a solution 
for $(V_i^{n+1})_\up{c}$, $u_{i+1/2}^{n+1}$, and $p_i^{n+1}$ is sought, that is, the solution of the
discrete wave equation \eqref{eqn:nlsystem} \textit{per se}, and a second one (an outer loop)
which begins with running to convergence an instance of the inner loop and subsequently applies
diffusion and filtering operators to the discrete volume and pressure.
The final result of the outer loop consists of the conservative state vector composed of 
$V_i^{n+1}$, $u_{i+1/2}^{n+1}$, and $E_i^{n+1}$.
For each of the two loops, the iteration is stopped as soon as the difference between two successive updates 
falls below a prescribed tolerance.


The discretisation approach presents similarities to 
 the ones proposed in \cite{bustohybrid} and \cite{casullidumbser} 
 for the Eulerian variant of the governing equations: namely the use of a staggered grid is key
 to obtaining a simple scalar tridiagonal linear system for the predictor wave equation.

\subsubsection{Numerical diffusion for the specific volume} \label{sec:densityfiltering}
The update formula \eqref{eqn:volumeupdate} for the specific volume is a second order central difference, and does not include any artificial
viscosity or Riemann solver dissipation. For this reason,
undesirable spurious oscillations can be introduced by it. In order to inhibit the development of 
such numerical artifacts, we introduce a \textit{stabilisation/filtering procedure} composed of two steps.
First, the candidate solution $\left(V_i^{n+1}\right)_\up{c}$ given by \eqref{eqn:volumeupdate} is filtered in order to detect
cells affected by excessive oscillations, and for each cell an alternative non-oscillatory value 
$V_i^{n+1,\ast}$ is computed. Second, since simply replacing in each cell the candidate solution
$\left(V_i^{n+1}\right)_\up{c}$ with $V_i^{n+1,\ast}$ might lead to unacceptable mass conservation errors, each 
value of the specific volume $\left(V_i^{n+1}\right)_\up{c}$ is updated with a conservative numerical diffusion flux, 
yielding
a corrected discrete volume field $V_i^{n+1}$ satisfying mass conservation. The set of conservative fluxes is determined in such
a way that its effect is as close as possible to that of the nonconservative cell-by-cell replacement of $\left(V_i^{n+1}\right)_\up{c}$ with $V_i^{n+1,\ast}$. An important difference of course is 
that the flux form
\begin{equation}\label{eqn:diffusionflux}
    V_i^{n+1} = \left(V_i^{n+1}\right)_\up{c} - \frac{\Delta t}{\Delta m_i}\,\left(f_{i+1/2}^{V,\up{d}} - f_{i-1/2}^{V,\up{d}}\right)
\end{equation}
ensures that mass is globally conserved for any choice of the diffusive fluxes $f_{i+1/2}^{V,\up{d}}$, remarking that, in the same way,  $\left(V_{i}^{n+1}\right)_c$ is also computed by a comparable Finite Volume-type formula \eqref{eqn:volumeupdate}, which satisfies conservation exactly as formula \eqref{eqn:diffusionflux}.
In order to achieve the desired diffusive effect, the fluxes are
set to have the form 

\begin{equation} \label{eq:consdiffusionflux}
f_{i+1/2}^{V,\up{d}} = -S_{i+1/2}\left[\left(V_{i+1}^{n+1}\right)_\up{c} - \left(V_{i}^{n+1}\right)_\up{c}\right],
\end{equation}
with $S_{i+1/2}$ a non-negative coefficient analogous to the signal speed estimates of the Rusanov or Local Lax Friedrichs approximate fluxes.

\paragraph{Computation of the target filtered volume}
The first step of the correction procedure applied to the specific volume consists in the determination of 
a filtered 
specific volume $V_i^{n+1,\ast}$. The quantity can also be termed
\textit{target} volume, since it is not used directly to replace $\left(V_{i}^{n+1}\right)_\up{c}$, that is, the volume obtained by the conservative central difference update formula \eqref{eqn:volumeupdate}, but
rather to determine the coefficients $S_{i+1/2}$ of the conservative diffusion fluxes \eqref{eq:consdiffusionflux}, aiming at obtaining a result that approximates 
the target volume $V_i^{n+1,\ast}$ as closely as possible, while retaining exact conservation of mass.

Towards this goal, in each cell $i$, \textit{if and only if the cell is a local extremum of the specific volume}, three candidates are made available:
1. the volume in the adjacent cell on the left $(V_i^{n+1,\ast})_1 = \left(V_{i-1}^{n+1}\right)_\up{c}$, 2. the volume in the right 
adjacent cell
 $(V_i^{n+1,\ast})_2 = \left(V_{i+1}^{n+1}\right)_\up{c}$, and 3. a third value obtained by a least-squares  reconstruction 
 involving a five cell stencil, of which only four are actually used: the two right neighbours and the two left neighbours of cell $i$. Contrary to other polynomial
 reconstruction operators \cite{barthfrederickson, DumbserKaeser06b, cravero2018cweno} used in high order Finite Volume method, the cell average of the central cell 
 $\left(V_{i}^{n+1}\right)_\up{c}$ is completely ignored, since it coincides with the value potentially requiring a replacement. 
 A local extremum is defined operatively as any cell where the left and right slopes switch sign, that is where
 \begin{equation}
 \left[(V_{i+1}^{n+1})_\up{c} - (V_{i}^{n+1})_\up{c})\right]\,\left[(V_{i-1}^{n+1})_\up{c} - (V_{i}^{n+1})_\up{c}\right] > \epsilon,
 \end{equation}
 having set $\epsilon = 10^{-14}$, a small tolerance preventing random switching-type behaviour in 
 completely flat regions.
 
 For the reconstruction procedure that has to be carried out at extrema to compute $(V_i^{n+1,\ast})_3$, formally we define 
 the vector of the stencil data for cell $i$ to be
\begin{equation}
          \vec{V}_i = 
     \left(
     (V_{i-2}^{n+1})_\up{c},\ 
     (V_{i-1}^{n+1})_\up{c},\ 
     (V_{i+1}^{n+1})_\up{c},\ 
     (V_{i+2}^{n+1})_\up{c}
     \right)^\transpose.
\end{equation}
 The least-square reconstructed candidate value $(V_i^{n+1,\ast})_3$ can be then computed by 
 solving the normal equations of the overdetermined system given by imposing conservation over each one of the four cells in the
 stencil with data $\vec{V}_i$. 
 The four reconstruction equations are obtained by imposing that the cell average $(V_{j}^{n+1})_\up{c}$ of each cell $j$ in 
 the stencil spanning from cell $i-2$ to cell $i+2$ (cell $i$ excluded) be the same as that of 
 a stencil-supported quadratic polynomial 
 \begin{equation}
 P_i(m) = \sum_{k=0}^{2}\hat{P}_{i,k}\,\left(m - m_{i-2}\right)^k = \hat{P}_{i,0} + 
 \hat{P}_{i,1}\,\left(m - m_{i-2}\right) + \hat{P}_{i,2}\,\left(m - m_{i-2}\right)^2,
 \end{equation}
 averaged over the same cell $j$.
The reconstruction equations are formally expressed as
 \begin{equation} \label{eq:reconstruction}
 \frac{1}{\Delta m_j} \int_{m_j}^{m_j + \Delta m_j} P_i(m)\,\de{m} = (V_{j}^{n+1})_\up{c}, \qquad j\in \left\{i-2,\ i-1,\ i+1,\ i+2\right\}
 \end{equation}
  recalling that $m_j$ is the leftmost point of cell $j$ and $\Delta m_j = m_{j+1} - m_j$.
  Since the four equation system \eqref{eq:reconstruction} is overdetermined, it is solved
  in the least squares sense, yielding a vector of three coefficients $\hat{\vec{P}}_i$.
  The polynomial $P_i(m)$ is then integrated over cell $i$ to obtain 
  the reconstructed central cell average $(V_i^{n+1,\ast})_3$: in practice it is sufficient 
  to take the dot product of the coefficients $\hat{\vec{P}}_i$ with an averaging operator $\vec{a}_i$, whose three components are the integrals of the three basis monomials $1$, $m - m_{i-2}$, and $(m-m_{i-2})^2$ over cell $i$, divided by its size $\Delta m_i$.
 In our implementation, the solution to the normal equations and the averaging of 
 the polynomial $P_i(m)$ over cell $i$ are computed as
 \begin{equation}
     (V_i^{n+1,\ast})_3 = \hat{\vec{P}}_i\cdot\vec{a}_i = \left[\left(\vec{A}_i^\transpose\,\vec{A}_i\right)^{-1}\,\vec{A}_i^\transpose\,\vec{V}_i\right]\cdot\vec{a}_i
 \end{equation}
 with the definitions 
 \begin{equation}
     \vec{A}_i = \begin{pmatrix*}
     1 & \Delta m_{i-2}/2 & \Delta m_{i-2}^2/3 + \Delta m_{i-2}\, m_{i-2} + m_{i-2}^2\\[1mm]
     1 & \Delta m_{i-1}/2 & \Delta m_{i-1}^2/3 + \Delta m_{i-1}\, m_{i-1} + m_{i-1}^2\\[1mm]
     1 & \Delta m_{i+1}/2 & \Delta m_{i+1}^2/3 + \Delta m_{i+1}\, m_{i+1} + m_{i+1}^2\\[1mm]
     1 & \Delta m_{i+2}/2 & \Delta m_{i+2}^2/3 + \Delta m_{i+2}\, m_{i+2} + m_{i+2}^2\\[1mm]
     \end{pmatrix*}, \qquad
     \vec{a}_i = 
     \begin{pmatrix*}
     1\\[1mm]
     \Delta m_i/2 + m_{i}\\[1mm]
     \Delta m_i^2/3 + \Delta m_i\,m_{i} + m_i^2\\[1mm]
     \end{pmatrix*}.
 \end{equation}

 Finally, the filtered specific volume $V_i^{n+1,\ast}$ is a weighted average of the three candidates
 \begin{equation}
  V_i^{n+1,\ast} = \sum_{k = 1}^3 w_k\, (V_i^{n+1,\ast})_k,\quad\text{ with }\quad
     w_k = \frac{\widetilde{w_k}}{\sum_{k=1}^3 \widetilde{w_k}}, \qquad \widetilde{w_k} = \left(\frac{1}{d_k + \epsilon}\right)^{8}, \qquad \epsilon = 10^{-14}
 \end{equation}
 where the weights $w_k$ are computed similarly to the classic WENO nonlinear coefficients
 \cite{shu1997, LPR:99, jiang1996}, but with a CWENO-type \cite{cravero2019, cravero2018cweno} approach, using the aggressive choice of parameters given in \cite{dumbserweno, dumbsercweno}, 
 and with $d_k$ being a simple indicator proportional to the distance of each of the candidate values for
 the filtered specific volume and the value to be corrected itself.
 Formally, for what concerns the expression of the indicator $d_k,\ k\in{1,\ 2,\ 3}$, we set
 \begin{equation}
     d_1 = (V_i^{n+1,\ast})_1 - \left(V_i^{n+1}\right)_\up{c}, \qquad
     d_2 = (V_i^{n+1,\ast})_2 - \left(V_i^{n+1}\right)_\up{c}, \qquad
     d_3 = K\,\left[(V_i^{n+1,\ast})_3 - \left(V_i^{n+1}\right)_\up{c}\right],
 \end{equation}
 which states that the weighted average will be biased towards correcting $\left(V_i^{n+1}\right)_\up{c}$
 as little as possible: the indicator $d_k$ is simply a measure of the distance between each candidate value $(V_i^{n+1,\ast})_k$ and the cell average $(V_i^{n+1})_\up{c}$ obtained by the conservative update formula \eqref{eqn:volumeupdate}. Additionally, the least-squares reconstructed value $(V_i^{n+1,\ast})_3$ is slightly penalised (by a multiplicative factor $K=4$ in its indicator $d_3$) 
 with respect to 
 the other two candidates since at discontinuities the other two are clearly more reliable, while
 in smooth regions, where the data are well approximated by the reconstruction polynomial and the least-squares option is thus more accurate and preferable, the match with ${\left(V_i^{n+1}\right)}_\up{c}$ will be such that ${(V_i^{n+1,\ast})}_3$ will have by far the largest weight regardless of the penalisation factor $K$.

\begin{figure}[!bp]
\centering
\includegraphics[width=0.495\textwidth]{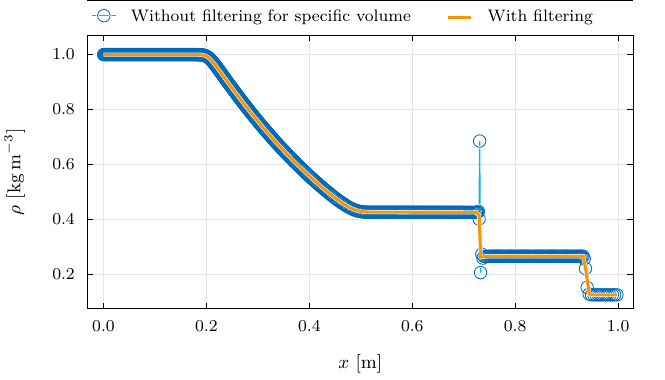}%
\includegraphics[width=0.495\textwidth]{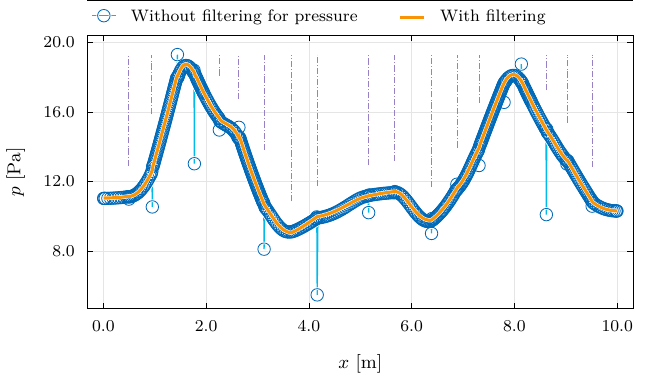}%
\caption{Motivation for the filtering techniques introduced in this work. On the left: central fluxes 
for the density update introduce spurious oscillations in correspondence of contact discontinuities. 
The oscillations are eliminated by numerical diffusion-like filtering (Section~\ref{sec:densityfiltering}). On the right: pressure oscillations
appear at material interfaces, across which the parameters of the equation of state are discontinuous (the interface
positions are highlighted with dashed lines). The oscillations disappear with filtering (Section~\ref{sec:pressurefiltering}). Both simulations were run with 
SDIRK2 time-stepping and $N=1000$ cells. For the left panel $\cfl=5.0$ and for the right panel $\cfl=1000$.}
\label{fig:issues}
\end{figure}

\paragraph{Conservative redistribution of specific volume correction}
Expanding on the conservative diffusive update formula given in \eqref{eqn:diffusionflux}, 
we set the diffusive flux to be of the form \eqref{eq:consdiffusionflux}, here replicated for the reader's convenience
\begin{equation} \label{eqn:diffusionfluxrus}
f_{i+1/2}^{V,\up{d}} = -S_{i+1/2}\,\left[(V_{i+1}^{n+1})_\up{c} - (V_i^{n+1})_\up{c}\right].
\end{equation} 
The numerical factor (one half) that would be present in the Rusanov flux is omitted here in \eqref{eqn:diffusionfluxrus}, 
since the values of $S_{i+1/2}$ are not derived from some bounding estimates for the signal speeds of the PDE system, but rather computed as the result of an optimisation procedure aiming towards obtaining
a close approximation of the (nonconservative) corrected values $V_i^{n+1,\ast}$, while retaining exact conservation.
Hence the conservative diffusive correction for the specific volume reads
\begin{equation}\label{eqn:diffusionconservative}
    V_i^{n+1} = \left(V_i^{n+1}\right)_\up{c} + \frac{\Delta t}{\Delta m_i}\,S_{i-1/2}\,\left[\left(V_{i-1}^{n+1}\right)_\up{c} - \left(V_{i}^{n+1}\right)_\up{c}\right] 
    + \frac{\Delta t}{\Delta m_i}\,S_{i+1/2}\,\left[\left(V_{i+1}^{n+1}\right)_\up{c} - \left(V_{i}^{n+1}\right)_\up{c}\right].
\end{equation}

In the following paragraphs, we outline the procedure by which a set of coefficients $S_{i+1/2}$ can be assigned to each
interface, so to mimic the output of the nonconservative filtering procedure given in the beginning of Section~\ref{sec:densityfiltering} by means of a 
strictly conservative formulation.

In each cell, we would like that the corrected specific volume $V_i^{n+1}$ obtained from the
conservative formula \eqref{eqn:diffusionconservative} be as close as possible to the
non-conservative filtered value $V_i^{n+1,\ast}$.
Finding the cell boundary coefficients $S_{i+1/2}$ that allow the best approximation of 
the target $V_i^{n+1} \rightarrow V_i^{n+1,\ast}$ over the $N$ cells of the computational domain
constitutes a rather large {optimisation} problem 
in the $N+1$ unknowns $S_{i+1/2}$. 

On the contrary, within 
each cell, if the effect of the conservative diffusion flux on the two neighbours is ignored, the target can be 
can be trivially imposed
as an exact constraint $V_i^{n+1} = V_i^{n+1,\ast}$ leading to 
\begin{equation} \label{eqn:diffusionconstraint}
\left(V_i^{n+1}\right)_\up{c} + \diffal{i}\,\diffsl{i} + \diffar{i}\,\diffsr{i} = V_i^{n+1,\ast},
\end{equation}
with
\begin{equation}
\diffal{i} = \frac{\Delta t}{\Delta m_i}\,\left[\left(V_{i-1}^{n+1}\right)_\up{c} - \left(V_{i}^{n+1}\right)_\up{c}\right], \quad
\diffar{i} = \frac{\Delta t}{\Delta m_i}\,\left[\left(V_{i+1}^{n+1}\right)_\up{c} - \left(V_{i}^{n+1}\right)_\up{c}\right].
\end{equation}

Note that \eqref{eqn:diffusionconstraint} has multiple solutions and in order to obtain a unique pair 
$\diffsl{i}, \diffsr{i}$, an additional constraint has to be imposed, specifically we choose to select the solution that has minimum square norm $\norm{\vec{S}_i}^2 = \diffsl{i}^2 + \diffsr{i}^2$, as explained in 
\eqref{eqn:minimumnorm}.
The definition of a unique intercell face value
\begin{equation}
S_{i+1/2} = \frac{1}{2} \, \left(\diffsl{i} + \diffsr{i}\right)
\end{equation}
guarantees that the resulting update \eqref{eqn:diffusionconservative} due to artificial diffusion will be conservative.


Below we explain in detail the procedure adopted to compute $\diffsl{i}$ and $\diffal{i}$.
In order to enforce equality \eqref{eqn:diffusionconstraint}, one 
one only has to express either $\diffsl{i}$ or $\diffsr{i}$ as a linear function of the other.
For example, if the coefficient $\diffal{i}$ is not zero (meaning that there is a jump in 
the specific volume $V$ between cell $i$ and cell $i-1$), then one can set
\begin{equation} \label{eqn:diffusionsl}
    \diffsl{i} = \frac{V_i^{n+1,\ast} - V_i^{n+1}}{\diffal{i}} - \frac{\diffar{i}}{\diffal{i}}\,\diffsr{i},
\end{equation}
or equivalently, if $\diffar{i} > 0$, i.e. when a nonzero jump is present at the interface $i+1/2$, one can impose
\begin{equation} \label{eqn:diffusionsr}
    \diffsr{i} = \frac{V_i^{n+1,\ast} - V_i^{n+1}}{\diffar{i}} - \frac{\diffal{i}}{\diffar{i}}\,\diffsl{i}.
\end{equation}
Within cell $i$, a local optimality condition can be then formulated by 
finding the two component vector $\vec{S}_i = (\diffsl{i},\ \diffsr{i})^\transpose$ which has minimum norm.
This is immediately obtained starting for example from \eqref{eqn:diffusionsl}, computing the norm of $\vec{S}_i$
\begin{equation} \label{eqn:minimumnorm}
    \norm{\vec{S}_i}^2 = \diffsl{i}^2 + \diffsr{i}^2 = B_0^2 - 2\,B_0\,B_1\,\diffsr{i} + \left(1 + B_0^2\right)\,\diffsr{i},
\end{equation}
and setting its derivative with respect to $\diffsr{i}$ to zero, obtaining
\begin{equation}
    \diffsr{i} = \frac{B_0\,B_1}{1 + B_1^2}, \qquad \diffsl{i} = B_0 - B_1\,\diffsr{i},
\end{equation}
with the definitions (omitting the index $i$ for the sake of notational compactness)
\begin{equation}
    B_0 = \frac{V_i^{n+1,\ast} - \left(V_i^{n+1}\right)_\up{c}}{\diffal{i}}, \qquad B_1 = \frac{\diffar{i}}{\diffal{i}}.
\end{equation}
Provided that Equation \eqref{eqn:diffusionsr} is valid, 
the same procedure 
can be carried out starting from it instead of \eqref{eqn:diffusionsl}, leading to the same results.
In order to avoid division by zero and minimize roundoff errors, in practice we choose \eqref{eqn:diffusionsl}
if $\abs{\diffal{i}} \ge \abs{\diffar{i}}$ and \eqref{eqn:diffusionsr} otherwise.
If \eqref{eqn:diffusionsr} is chosen, then the sought coefficients are
\begin{equation}
    \diffsl{i} = \frac{C_0\,C_1}{1 + C_1^2}, \qquad \diffsr{i} = C_0 - C_1\,\diffsl{i},
\end{equation}
having defined (again omitting the implied index $i$)
\begin{equation}
    C_0 = \frac{V_i^{n+1,\ast} - \left(V_i^{n+1}\right)_\up{c}}{\diffar{i}}, \qquad C_1 = \frac{\diffal{i}}{\diffar{i}}.
\end{equation}
In our implementation, towards the goal of preventing symmetry-breaking effects from 
floating point roundoff errors as much as possible, whenever both
$\diffal{i}$ and $\diffar{i}$ are strictly positive (in practice greater than $10^{-8}$) and differ
from each other by less than the same threshold of $10^{-8}$, then both \eqref{eqn:diffusionsl}
and \eqref{eqn:diffusionsr} are used to generate two values of $\diffsl{i}$ and $\diffsr{i}$ in cell
$i$, which are then arithmetically averaged.

If both $\diffal{i}$ and $\diffar{i}$ fall below the $10^{-8}$ threshold, then $\diffsl{i} = 0,$ and
$\diffsr{i} = 0$ for cell $i$, which corresponds to no diffusion. 
Note that in principle there is no guarantee that the conservative redistribution of these diffusive terms
will preserve the non oscillatory character of the (non-conservative) corrected values of the specific 
volume, but in practice, due to the iterative nature of the scheme, the correction tends to 
behave well and does not introduce the oscillations back.


\subsubsection{Update of the velocity field}
The velocity field $u_{i+1/2}^{n+1}$ is updated at the edge locations according to \eqref{eqn:velocityupdate}, and 
the auxiliary cell-center values $u_i^{n+1}$ are computed immediately after by interpolation, that is by defining an interpolated cell center velocity
\begin{equation}
u_i^{n+1} = \frac{\Delta m_{i+1/2}\,u_{i-1/2} + \Delta m_{i-1/2}\,u_{i+1/2}}{\Delta m_{i-1/2} + \Delta m_{i+1/2}}.
\end{equation}
It should be stressed that $u_i^{n+1}$ is used just for the computation of specific total energy
whenever the quasi-conservative pressure corrections detailed in Section~\ref{sec:pressurefiltering}
are applied, and the quantity to be integrated in time is $u_{i+1/2}$, since interpolating backwards
from cell centers to edge values would introduce Lax-Friedrichs-like dissipation terms \cite{lax1954},
quickly degrading the accuracy of the method. 
If more aggressive numerical diffusion is required, as we found to be necessary for the two of the shock tube problems of Toro 
shown in Section~\ref{sec:rp}, we refer the reader to the procedure detailed in Appendix~\ref{appendix:extra_diffusion}.

\subsubsection{Conservative update of energy, numerical diffusion and filtering}
\paragraph{Conservative update of total energy}
A candidate specific energy updated value is
which is based on a conservative flux form and thus works also in presence of shockwaves
\begin{equation} \label{eqn:pressureconservative}
    (E_i^{n+1})_\up{c} = E_i^n - \frac{\Delta t}{\Delta m_i}\left(
    ({F_E})_{i+1/2}^{n+1} - ({F_E})_{i-1/2}^{n+1}
    \right),
\end{equation}
with the fluxes being 
\begin{equation} \label{eq:energyflux}
\begin{aligned}
    &({F_E})_{i-1/2}^{n+1} = u_{i-1/2}^{n+1}\,p_{i-1/2}^{n+1} - \frac{1}{2}\,({s_E})_{i-1/2}\,\left(p_{i-1/2,+}^{n+1}   - p_{i-1/2,-}^{n+1}\right),\\
    &({F_E})_{i+1/2}^{n+1} = u_{i+1/2}^{n+1}\,p_{i+1/2}^{n+1} - \frac{1}{2}\,({s_E})_{i+1/2}\,\left(p_{i+1/2,+}^{n+1} -  p_{i+1/2,-}^{n+1} \right).
\end{aligned}
\end{equation}
The energy fluxes \eqref{eq:energyflux} are the sum of a central part and a diffusive part, just like the Rusanov (or local Lax--Friedrichs) \cite{rusanov1961}
and Ducros \cite{ducros2000} fluxes.
Note however that the jumps $p_{i-1/2,+}^{n+1}   - p_{i-1/2,-}^{n+1}$ are not the jumps of the conserved variables (the specific energy $E$) but rather
they are pressure jumps. This is motivated by the fact that diffusive part of \eqref{eq:energyflux} is intended
to introduce as little numerical diffusion as possible and in particular at material interfaces we can expect 
the pressure field to be smooth. On the other hand, even in complete absence of shock waves, the specific
energy $E$ might jump from one cell to another in a discontinuous manner, for example at material interfaces where the parameters of the equation of state change.
The central part simply requires the computation of a cell-edge interpolated pressure
\begin{equation}
    p_{i+1/2}^{n+1} = \frac{\Delta m_{i+1}\,p_{i}^{n+1} + \Delta m_{i}\,p_{i+1}^{n+1}}{\Delta m_{i} + \Delta m_{i+1}},
\end{equation}
and the diffusive part is a function of the minmod-reconstructed pressure jumps $p_{i+1/2}^+ - p_{i+1/2}^-$ at each cell interface, also computed starting from 
the discrete wave equation solution $p_i^{n+1}$. The boundary extrapolated values $p_{i+1/2}^-$ and $p_{i+1/2}^+$ are defined
as
\begin{equation}
p_{i+1/2,-}^{n+1} = p_{i}^{n+1} + \delta p_{i}/2,\qquad p_{i+1/2,+}^{n+1} = p_{i+1}^{n+1} - \delta p_{i+1}/2,
\end{equation}
and the pressure increments due to the reconstruction procedure are 
\begin{equation}
    \delta p_{i} = \Delta m_i \, \up{minmod} \left(\frac{p_i - p_{i-1}}{\Delta m_{i-1/2}},\ \frac{p_{i+1} - p_i}{\Delta m_{i+1/2}}\right),
\end{equation}
with the minmod function defined by
\begin{equation}
\up{minmod}\left(\delta_\up{L},\ \delta_\up{R}\right) = \frac{1}{2}\,\left(\flll\up{sign}\left(\delta_\up{L}\right) + \up{sign}\left(\delta_\up{R}\right)\right)\,\min\left(\flll\abs{\delta_\up{L}},\ \abs{\delta_\up{R}}\right)
\end{equation}
The coefficients $(s_E)_{i-1/2}$ and $(s_E)_{i+1/2}$ are Rusanov-type signal speed estimates computed as
\begin{equation}
    (s_E)_{i+1/2} = \max\left(\frac{V_i^{n+1}}{\gamma_i - 1},\ \frac{V_{i+1}^{n+1}}{\gamma_{i+1} - 1}\right)\,\max\left(\frac{u_i^{n+1}}{V_{i}^{n+1}},\ \frac{u_{i+1}^{n+1}}{V_{i+1}^{n+1}}\right),
\end{equation}
where the first scaling coefficient, proportional to $V/(\gamma - 1)$, is intended to convert pressure jumps into jumps of specific internal energy $E$. Since such scaling factor introduces a dimensional constant proprtional to $V$, the correct dimensionality of $(s_E)_{i+1/2}$ (a velocity) is restored by taking the maximum of the momenta in the cells $i$ and $i+1$.
If the two fluids have similar densities and equation of state and the velocity is continuous across an interface, the factor $(s_E)_{i+1/2}$ is approximately $(s_E)_{i+1/2} \sim {u_{i+1/2}^{n+1}}/(\gamma_{i+1/2} - 1)$, whereas if a significant jump in density is present, the coefficient $(s_E)_{i+1/2}$ incorporates a scaling factor
proportional to the ratio between the maximum and the minimum density, effectively increasing the dissipation at shocks and material interfaces.

It should be remarked that such estimates, in order to introduce as little numerical dissipation as possible, are based on the flow velocity and not on the speed of sound or the eigenvalues of 
the governing hyperbolic system, as they
would be for a standard explicit scheme. Moreover, due to the use of pressure 
jumps in lieu of specific energy jumps, the signal speeds $(s_E)_{i+1/2}$ include a dimensional scaling coefficient, based on the derivative of specific energy $E$ with respect to pressure, so that dimensionality
and scaling of the involved quantities are accounted for.

\paragraph{Removal of pressure oscillations/merging different candidates for the pressure field} \label{sec:pressurefiltering}

With reference to the second panel of Figure~\ref{fig:issues}, 
and taking into account for the observations given in \cite{abgrall2001}, we introduce 
a quasi-conservative correction which is intended to eliminate the pressure spikes that arise
at material interfaces, where the parameters of the equation of state jump between different values.
We point out that such small non-conservative corrections, as discussed in \cite{gaburro2024}, 
do not lead to highly undesirable effects like shockwaves propagating with the wrong speed, 
precisely due to their small pointwise nature. This statement is tested in the paper 
by applying the method to a variety of classic benchmarks featuring strong shocks.
Towards this goal, we simply compute the specific energy $E_i^{n+1}$ in a cell by cell fashion as a convex
combination of its conservative updated value $(E_i^{n+1})_\up{c}$ and a second value $(E_i^{n+1})_\up{p}$
obtained solely based on predictor information, which is nonconservative in principle.
\begin{equation}
    E_i^{n+1} = \left(1 - w_i\right)\,(E_i^{n+1})_\up{c} + w_i\,(E_i^{n+1})_\up{p}
\end{equation}
with $(E_i^{n+1})_\up{p}$ denoting the specific total energy that can be computed as a pointwise function
\begin{equation}
(E_i^{n+1})_\up{p} = (E_i^{n+1})_\up{p}\left(p_i^{n+1},\ u_{i-1/2}^{n+1},\ u_{i+1/2}^{n+1},\ V_i^{n+1}\right) = V_i^{n+1}\,\frac{p_i^{n+1} + 
    \gamma_i\,\Pi_i}{\gamma_i - 1} + \frac{(u_i^{n+1})^2}{2}.
\end{equation}
The blending weight $w_i$ is computed as 
\begin{equation} \label{eqn:pressureweight}
    w_i = \min\left(1,\ \frac{1}{D}\,\left|(p_i^{n+1})_\up{c} - \frac{\Delta m_{i+1/2}\,(p_{i-1}^{n+1})_\up{c} + \Delta m_{i-1/2}\,(p_{i+1}^{n+1})_\up{c}}{\Delta m_{i-1/2} + \Delta m_{i+1/2}}\right|\right)^r.
\end{equation}
The motivation for this choice of blending weight is given in the following paragraphs.
Here,  $(p_i^{n+1})_\up{c}$ denotes the pressure field induced by the conservative energy update, that is
\begin{equation}
    (p_i^{n+1})_\up{c} = \frac{\gamma_i - 1}{V_i^{n+1}}\left[\left(E_{i}^{n+1}\right)_\up{c} - \frac{(u_i^{n+1})^2}{2}\right] - \gamma_i\,\Pi_i,
\end{equation}
and $D$ being a scaling factor defined as 
\begin{equation}
    D = \max\left(p_{i-2}^{n+1},\ p_{i-1}^{n+1},\ p_{i+1}^{n+1},\ p_{i+2}^{n+1}\right) - 
\min\left(p_{i-2}^{n+1},\ p_{i-1}^{n+1},\ p_{i+1}^{n+1},\ p_{i+2}^{n+1}\right) + \epsilon.
\end{equation}
Finally we set the exponent $r = 8$ in \eqref{eqn:pressureweight} to be a fixed parameter used to smoothly approximate a binary switch between $w_i=0$ and $w_i=1$, and set $\epsilon=10^{-14}$, a small constant aimed at avoiding division by zero.
The purpose of \eqref{eqn:pressureweight} is to compare the pressure $(p_i^{n+1})_\up{c}$ given by the conservative update of energy \eqref{eqn:pressureconservative} in 
cell $i$ with the pressure obtained by linearly interpolating the two neighbouring values, thus measuring
the deviation from a simple linear distribution.
Wherever such discrepancy is high, relative to the scaling factor $D$ based on the pressure field $p_i$ (the one derived solely from 
the wave equation), a pressure spike has been detected and the weighting factor $w_i$ will tend towards unity, 
while in smooth regions of the pressure field the correction is essentially null, since the match between 
the pressure $(p_i^{n+1})_\up{c}$ and the neighbour-interpolated value will be much smaller than the feature size
indicated by $D$, thereby driving $w_i$ towards zero, due to the high exponent $r=8$, approximating a rather sharp switch.
This means that the heuristic \eqref{eqn:pressureweight} is designed in such a way that the pressure field 
is replaced by the smooth predictor solution whenever quantified as oscillatory, while safeguarding the 
conservativity of the method as much as possible since in most flow regions the nonconservative contribution is negligible.



\subsection{Time integration}\label{sec:dirk}
The time integration approach adopted in this work is very simple and leverages the Stiffly-accurate
Diagonally implicit Runge--Kutta schemes originally introduced by \cite{crouzeixsdirk, alexandersdirk} (see also \cite{ascher1997implicit, pareschi2005implicit} and the Review by Kennedy and Carpenter \cite{rkreview}).
Formally, we define an implicit solution operator $\mathcal{H}_\up{I}(\vec{Q},\ \Delta t)$ that maps 
an initial state vector $\vec{Q} = \left(V_i,\ u_{i+1/2},\ E_{i}\right)$, 
with $i=1,\ 2,\ \hdots,\ N$ (and including $i=0$ for the left boundary value of the velocity field), 
to an updated vector $\mathcal{H}_\up{I}(\vec{Q}^n,\ \Delta t)$, given by our implicit discretisation scheme.
With reference to the notation adopted in the above sub-sections, the updated values 
for a timestep of size $\Delta t$ will be 
\begin{equation}
\mathcal{H}_\up{I}(\vec{Q}^n,\ \Delta t) = \left(V^{n+1}_i,\ u^{n+1}_{i+1/2},\ E_i^{n+1}\right).
\end{equation}
The approach proposed in this paper is tested with two time integration schemes, one of second order (labelled SDIRK2)
with two stages, defined by the Butcher array
\begin{equation}
    \begin{tabular}{c|cc}
        $\gamma$ & $\gamma$ & $0$ \\
        $1$      & $a_{21}$ & $\gamma$\\
        \midrule
                 & $a_{21}$ & $\gamma$\\
    \end{tabular}
    \qquad 
    \text{with}
    \qquad
    \begin{aligned}
        & \gamma = 1 - 1/\sqrt{2}, \\
        & a_{21} = 1 - \gamma, \\
    \end{aligned}
\end{equation}
and implemented in practice as a sequence of two Euler steps
\begin{equation}
\begin{aligned}
    &\vec{Q}_1 = \solutionop_{\up{I}}\left(\vec{Q}_1^\ast,\ \gamma\,\Delta t\right), \\
    &\vec{Q}^{n+1} = \solutionop_{\up{I}}\left(\vec{Q}_2^\ast,\ \gamma\,\Delta t\right), \\
\end{aligned}
    \qquad 
    \text{with}
    \qquad
\begin{aligned}
    &\vec{Q}_1^\ast = \vec{Q}^n,\\
    &\vec{Q}_2^\ast = c_{20}\,\vec{Q}^n + c_{21}\,\vec{Q}_1.\\
\end{aligned}
\end{equation}
The coefficients $c_{20}$ and $c_{21}$ are
$c_{20} = 1 - a_{21}/\gamma$ and $c_{21} = a_{21}/\gamma$.

A second time integration scheme is the three-stage, third order scheme defined by
\begin{equation}
    \begin{tabular}{c|ccc}
        $\gamma$  & $\gamma$ &      $0$ &      $0$\\
        $k_c$     & $a_{21}$ & $\gamma$ &      $0$\\
        $1$       & $a_{31}$ & $a_{32}$ & $\gamma$\\
        \midrule
                  & $a_{31}$ & $a_{32}$ & $\gamma$\\
    \end{tabular}
    \qquad 
    \text{with}
    \qquad
    \begin{aligned}
        &\gamma = 0.435866521508459,\\
        &a_{21} = k_c - \gamma,\\
        &a_{31} = 1 - k_d - \gamma,\\
        &a_{32} = k_d,\\
    \end{aligned}
    \qquad
    \begin{aligned}
        &k_a = 1 - 4\,\gamma + 2\,\gamma^2,\\
        &k_b = 3\,\gamma\,\left(2 - 3\,\gamma + \gamma^2\right) - 1,\\
        &k_c = (2/3 - 3\,\gamma + 2\,\gamma^2)/k_a,\\
        &k_d = -3\,k_a^2/(4\,k_b),\\
    \end{aligned}
\end{equation}
for which the implementation reads
\begin{equation} \label{eqn:rkvariant}
\begin{aligned}
    &\vec{Q}_1 = \solutionop_{\up{I}}\left(\vec{Q}_1^\ast,\ \gamma\,\Delta t\right), \\
    &\vec{Q}_2 = \solutionop_{\up{I}}\left(\vec{Q}_2^\ast,\ \gamma\,\Delta t\right), \\
    &\vec{Q}^{n+1} = \solutionop_{\up{I}}\left(\vec{Q}_3^\ast,\ \gamma\,\Delta t\right), \\
\end{aligned}
    \qquad 
    \text{with}
    \qquad
\begin{aligned}
    &\vec{Q}_1^\ast = \vec{Q}^n,\\
    &\vec{Q}_2^\ast = c_{20}\,\vec{Q}^n + c_{21}\,\vec{Q}_1,\\
    &\vec{Q}_3^\ast = c_{30}\,\vec{Q}^n + c_{31}\,\vec{Q}_1 + c_{32}\,\vec{Q}_2,\\
\end{aligned}
\end{equation}
The variation \eqref{eqn:rkvariant} is provided so that in practice one 
can conveniently implement the Runge--Kutta integrator purely as a sequence of implicit Euler steps.
The coefficients $c_{ij}$ are
\begin{equation}
    c_{21} = a_{21}/\gamma,\qquad c_{20} = 1 - c_{21},\qquad c_{32} = a_{32}/\gamma,\qquad
    c_{31} = (a_{31} - c_{32}\,a_{21})/\gamma,\qquad c_{30} = 1 - c_{31} - c_{32}.
\end{equation}


\subsection{Explicit reference schemes}
The explicit Lagrangian schemes used for reference in the convergence and efficiency studies are the
based on the well-known method of Munz \cite{munz94}, with the addition of an exact two-material 
Riemann solver two Strong-Stability-Preserving Runge-Kutta time-stepping schemes \cite{gottliebssprk}.
In the notation adopted in this paper, the scheme of Munz, including a low-dissipation TVD reconstruction 
operator (see \cite{vanleer1979, kurganovgminmod}), 
mapping an initial condition $\vec{Q}^n$ to its evolved state $\vec{Q}^{n+1}$ over a time interval $\Delta t$.
With  notation, the first one of the two explicit time-stepping schemes is the two stage 
method \cite{gottliebssprk} (labelled SSPRK2)
having Butcher table
\begin{equation}
    \begin{tabular}{c|ccc}
        $0$  & $0$ &      $0$\\
        $1$     & $a_{21}$ & $0$\\
        \midrule
                & $b_{1}$ & $b_{2}$\\
    \end{tabular}
    \qquad
    \text{with}
    \qquad
    \begin{aligned}
        & a_{21} = 1,\\
        & b_1 = b_2 = 1/2,\\
    \end{aligned}
\end{equation}
which has a very simple implementation given by
\begin{equation}
\begin{aligned}
    \vec{Q}_1 &= \solutionop_\up{E}\left(\vec{Q}^n,\ \Delta t\right),\\
    \vec{Q}^{n+1} &= (1 - c_1)\,\vec{Q}^n + c_1\,\solutionop_\up{E}\left(\vec{Q}_2,\ \Delta t\right),\\
\end{aligned}
\qquad
\text{with}
\qquad
c_1 = 1/2.\\
\end{equation}

Furthermore, we compared our results with those produced by the same second-order TVD spatial reconstruction
and exact Riemann solver, but with a three-stage third order SSP Runge--Kutta scheme \cite{gottliebssprk}, 
labelled in this paper SSPRK3, defined by the Butcher table
\begin{equation}
    \begin{tabular}{c|ccc}
        $0$  & $0$ &      $0$ &      $0$\\
        $1$     & $a_{21}$ & $0$ &      $0$\\
        $1/2$ & $a_{31}$ & $a_{32}$ & $0$\\
        \midrule
                & $b_{1}$ & $b_{2}$ & $b_{3}$\\
    \end{tabular}
    \qquad
    \text{with}
    \qquad
    \begin{aligned}
        & a_{21} = 1,\\
        & a_{31} = a_{32} = 1/4,\\
        & b_1 = b_2 = 1/6,\\
        & b_3 = 2/3,
    \end{aligned}
\end{equation}
which in practice is implemented as
\begin{equation}
\begin{aligned}
    \vec{Q}_1 &= \solutionop_\up{E}\left(\vec{Q}^n,\ \Delta t\right),\\
    \vec{Q}_2 &= (1 - c_1)\,\vec{Q}^n + c_1\,\solutionop_\up{E}\left(\vec{Q}_1,\ \Delta t\right),\\
    \vec{Q}^{n+1} &= (1 - c_2)\,\vec{Q}^n + c_2\,\solutionop_\up{E}\left(\vec{Q}_2,\ \Delta t\right),\\
\end{aligned}
\qquad
\text{with}
\qquad
\begin{aligned}
    & c_1 = 1/4,\\
    & c_2 = 2/3.
\end{aligned}
\end{equation}


\section{Generating mass-constrained meshes with smoothly variable non-uniform spacing in mass coordinates}\label{sec:mesh}

We discuss the task of generating a computational mesh (in our simple one-dimensional context, this means a sequence of 
mesh spacing values $\Delta m_i$), with the aim of applying the implicit Lagrangian scheme to flows in stratified systems,
characterized by several layers, alternating regions of constant width having two different densities $\rho_1$ and $\rho_2$.

A first natural choice consists in selecting a uniform fixed spacing $\Delta m$, and
a second one specifies two values for the mesh spacing $\Delta m_1$ and $\Delta m_2$, for example
computed in such a way that $\Delta m_1 = \rho_1\,\Delta x$ and $\Delta m_2 = \rho_2\,\Delta x$ corresponding
to the same step $\Delta x$ in Eulerian coordinates.

The first choice, a uniform spacing $\Delta m$, when used in stratified media with alternating densities, 
suffers from two significant drawbacks. First, when the constant width layers have very different densities, 
the resolution in the lighter regions can be much lower, meaning that either excessively fine grids will
be required for the denser fluid layers, or that the flow in lighter layers will be severely under-resolved.

In the canonical example of water-air-like density ratio $\rho_2/\rho_1 = 1000$, if a single grid point is associated with each one 
of 10 air-like layers (which is totally inadequate for any flow condition), then 10000 additional cells will be 
required for keeping the uniform spacing in the 10 water-like layers.

The main issue with the second type of mesh, having two mass-spacing values $\Delta m_1$ and $\Delta m_2$
jumping across layer boundaries where the density switches from $\rho_1$ to $\rho_2$, 
is that the artificial diffusion of the numerical 
scheme (and in general the modified equation associated with it), will change discontinuously at boundaries, 
potentially giving rise to spurious oscillations and in general deteriorating the quality of the solution.

Hence it is clear that some strategy for defining a smoothly varying mesh is required for applications
involving stratified multilayer media.

A simple enough strategy for computing the Lagrangian mass-coordinate spacing values $\Delta m_i$ is given in the next paragraphs.

\begin{figure}[!bp]
\centering
\includegraphics[width=0.99\textwidth]{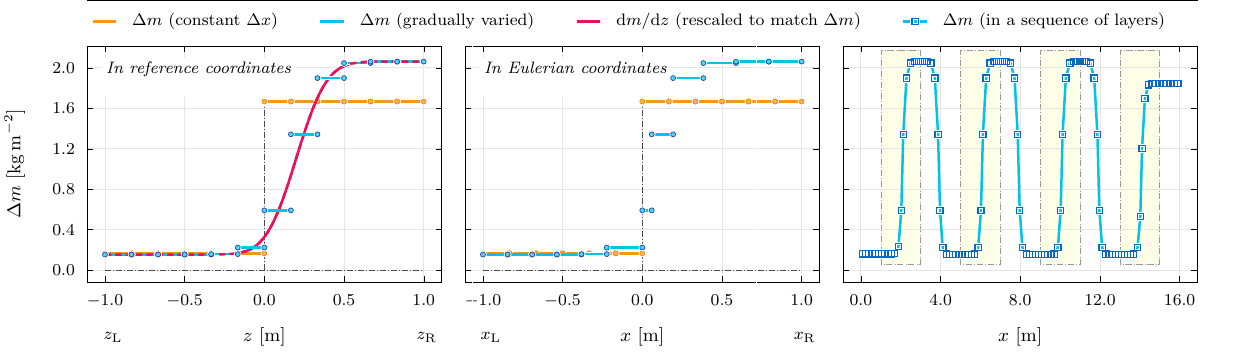}%
\caption{Illustration of the graded mesh generation procedure. In the left panel, we show
the cell mass content $\Delta m$ as a function of the auxiliary (reference) coordinate $z$, along with 
the mesh spacing that would be induced by fixed uniform value of $\Delta x$, and the mass density
function $\de m/\de z$ from which $\Delta m$ is integrated (graphically rescaled to match the plot of $\Delta m$). In the central panel, 
we show how the smooth transient in reference coordinates translates to Eulerian coordinates (that is, as a function of $x$).
In the right panel we show how several copies of the half-layer-couples are connected in a sequence of 
four layer couples. The regions corresponding to the central panel are highlighted. Note the presence of a boundary effect on the mesh spacing due to the fact that the first and last half-layer of the domain has a constant mesh spacing, unlike the internal ones.}
\label{fig:gradedmesh}
\end{figure}

\subsection{Mesh for a step density distribution}

In this section, we outline the procedure adopted for constructing a well graded Lagrangian mesh for a sub-domain 
composed of two regions of constant density. We refer to this sub-system as a half-layer pair, since it can be 
obtained by splitting each constant density layer of a stratified medium in two halves and focusing on a given interface 
between two layers.
Since the final mesh for the multilayer system will be composed of a series of such half-layer couples, suitably mirrored in an alternated sequence, studying this simple configuration of two half-layers will be sufficient to define
the final mesh in full. In particular, it will be enough to define $n$ values $\Delta m_i$ in each layer pair to specify the mesh in full.

With reference to Figure~\ref{fig:gradedmesh}, 
consider a pair of half-layers having density $\rho_\up{L}$ and 
$\rho_\up{R}$ in the left layer and in the right layer respectively. An auxiliary coordinate system $z$ (which can 
be thought of as an ``enumerating'' space where distance varies by a fixed amount per cell) is 
constructed, such that $z=0$ at the interface between the two layers, $z=z_\up{L}=x_\up{L}$ at the center of the left
layer (corresponding by definition to the left extremum of left half-layer) and $z=z_\up{R}=x_\up{R}$ at the center of the right
layer (corresponding by definition to the right extremum of right half-layer). Note that, while the $x$ and $z$ coordinates match at the interface between the two layers ($x=z=0$), and at the boundaries of the half-layer pair ($-1\,\up{m}$ and $+1\,\up{m}$ in this example), $z$ does not coincide at all with $x$ within each layer: 
in Figure~\ref{fig:gradedmesh} it is clear that $z$ increases exactly by the constant $\Delta z$ in each cell, while 
in the Eulerian reference frame each cell has a different width.

Within each half-layer pair, we will obtain a sequence of discrete mesh spacing values $\Delta m_i$ by defining and integrating a smooth \textit{mesh density function} $\de{m}/\de{z}(z)$, which will be constructed in such a way that
the width and mass content of each layer are described exactly, which is the main nontrivial constraint in generating 
Lagrangian mass-coordinate meshes for two-density systems if such grids are not simply taken to have layer-wise constant spacing as in \cite{phan2023numerical}.
Smoothness of the mesh at the junction points (where mirrored half-layer couples are glued together) is ensured by the fact that the mesh spacing is almost constant near both boundaries of a half-layer pair, 
that is, the mesh density assumes a sigmoid shape.
In the example of Figure~\ref{fig:gradedmesh}, the reference frame runs from $z_\up{L} = -1\,\up{m}$ to $z_\up{R} = 1\,\up{m}$.
In this reference frame we define each one of the $n$ mesh spacing as the integral of a density-like function, or formally
\begin{equation}
    \Delta m_i = \int_{z_1^i}^{z_2^i} \frac{\de {m}}{\de{z}}(z)\,\de{z} ,\quad \text{ with } z_1^i = {z_i - \Delta z/2},\ z_2^i = {z_i + \Delta z/2},\qquad i = 1,\ 2,\ \hdots,\ n
\end{equation}
and with $\Delta z$ specified as $\Delta z = (z_\up{R} - z_\up{L})/n$.
The number $n$ of cells of the half-layer pair can be chosen arbitrarily, and the procedure will automatically
generate a density-like function $\delta m$ which satisfies two mandatory requirements for any choice of $n$: first, the total mass contained
in the left half-layer must be equal to $M_\up{L} = \rho_\up{L}\,\left(-z_\up{L}\right)$, and likewise 
in the right half layer $M_\up{R} = \rho_\up{R}\,z_\up{R}$. 
Second, at the same time, the Eulerian mesh spacing values $\Delta x_i$ will sum to $-z_\up{L}$ on the left 
half-layer and sum to $z_\up{R}$ on the right half layer.
This means that the smoothly varying mesh is constructed in such a way that both the Eulerian coordinate $x$ 
and the Lagrangian mass coordinate $m$ of 
each material interface are the ones initially specified by the problem setup, meaning that the mesh is 
forced to preserve at the same time the width and mass content of each layer.
Formally the mass conservation constraints are imposed by requiring that 
\begin{equation}\label{eq:meshmasscons}
\int_{z_\up{L}}^0 \frac{\de {m}}{\de{z}}(z)\,\de{z} = M_\up{L} = -\rho_\up{L}\,z_\up{L}, \qquad \int_0^{z_\up{R}} \frac{\de {m}}{\de{z}}(z)\,\de{z} = M_\up{R} = \rho_\up{R}\,z_\up{R}, 
\end{equation}
meaning, the density function integrated over each half-layer 
must equal the total mass expected for that segment, 
so to ensure that both mass conservation and the positions of 
material interfaces are strictly maintained.

The density-like function is chosen so that it will smoothly interpolate the discrete mesh spacing $\Delta m_i$ 
between a left value $\Delta m_\up{L}$ and
a right value $\Delta m_\up{R}$. Specifically, the following form is adopted 
\begin{equation} \label{eq:meshdensity}
    \frac{\de{m}}{\de{z}}(z) = \Delta m_\up{L} + \frac{1}{2}\left(\Delta m_\up{R} - \Delta m_\up{L}\right)\,\left[1 + \up{erf}\left(\frac{z - z_0}{L_\beta}\right)\right].
\end{equation}
As a clarification note, we recall that $\mathrm{erf}(x) = \left({2}/{\sqrt{\pi}}\right) \int_0^x \exp{\left(-t^2\right)} \, \de t$, 
which gives a function normalized between -1 and 1. 
The formula is parametrised by $z_0$ (a coordinate shift) and $L_\beta$, a characteristic length controlling the width of the
transition zone between two constant regions in which the mesh spacing will be $\Delta m_\up{L}$ and 
$\Delta m_\up{R}$. In practice the length $L_\beta$ is computed as $L_\beta = \beta\,\min(-z_\up{L},\ z_\up{R})$, 
so the parameter to be fixed is a nondimensional value $\beta$.

For any fixed choice of $z_0$ and $\beta$, the 
mesh spacing at the left and right  $\Delta m_\up{L}$ and 
$\Delta m_\up{R}$ will be determined by imposing the 
mass conservation constraint \eqref{eq:meshmasscons} on \eqref{eq:meshdensity}, which yields
\begin{equation}
    \begin{pmatrix*}
    \Delta m_\up{L}\\
    \Delta m_\up{R}
    \end{pmatrix*} = {\begin{pmatrix*}
    -z_\up{L} + k_\up{L} & -z_\up{L} - k_\up{L}\\
    z_\up{R} - k_\up{L} & z_\up{R} + k_\up{R}\\
    \end{pmatrix*}}^{-1}\,\begin{pmatrix*}
    -\rho_\up{L}\,z_\up{L}\\
    \rho_\up{R}\,z_\up{R}\\
    \end{pmatrix*},
\end{equation}
having defined the auxiliary constants
\begin{equation}
    \begin{aligned}
        & k_\up{L} = \frac{L_\beta}{\sqrt{\pi}}\,\left\{
            \exp\left[-\left(\frac{-z_\up{L} + z_0}{L_\beta}\right)^2\right] - \exp\left[-\left(\frac{z_0}{L_\beta}\right)^2\right]\right\} + 
            \left(-z_\up{L} + z_0\right)\,\up{erf}\left(\frac{-z_\up{L} + z_0}{L_\beta}\right) - z_0\,\up{erf}\left(\frac{z_0}{L_\beta}\right),\\
        & k_\up{R} = \frac{L_\beta}{\sqrt{\pi}}\,\left\{
            \exp\left[-\left(\frac{z_\up{R} - z_0}{L_\beta}\right)^2\right] - \exp\left[-\left(\frac{z_0}{L_\beta}\right)^2\right]\right\} + 
            \left(z_\up{R} - z_0\right)\,\up{erf}\left(\frac{z_\up{R} - z_0}{L_\beta}\right) - z_0\,\up{erf}\left(\frac{z_0}{L_\beta}\right).\\
    \end{aligned}
\end{equation}

The mass content $\mathcal{M}(z_1,\ z_2)$ associated with a generic interval delimited by $z_1$ and $z_2$ is 
\begin{equation}
\begin{aligned}
    \mathcal{M}(z_1,\ z_2) = 
        \frac{L_\beta\,\left(\Delta m_\up{R} - \Delta m_\up{L}\right)}{2\,\sqrt{\pi}}\,
        \left[
            \exp\left(-\zeta_2^2\right) - \exp\left(-\zeta_1^2\right)\right]
         + \frac{1}{2}\,\left(\Delta m_\up{R} + \Delta m_\up{L}\right)\,\left(z_2 - z_1\right)&\ + \\
          +\ \frac{1}{2}\,\left(\Delta m_\up{R} - \Delta m_\up{L}\right)\,\left[
         \left(z_2 - z_0\right)\,\up{erf}\left(\zeta_2\right) - \left(z_1 - z_0\right)\,\up{erf}\left(\zeta_1\right)\right]&,
\end{aligned}
\end{equation}
with 
$\zeta_1 = (z_1 - z_0)/L_\beta$
and
$\zeta_2 = (z_2 - z_0)/L_\beta$.
Then the mesh spacing for each cell $i$ is directly computed by evaluating $\Delta m_i = \mathcal{M}(z_i-\Delta z/2,\ z_i + \Delta z/2)$.
Since the values of $\beta$ and $z_0$ which yield well graded meshes, or the range of selectable values itself, are not 
immediately apparent, a simple procedure to compute these parameters automatically is discussed in Appendix \ref{app:mesh}.

\section{Numerical results}\label{sec:results}

\reviewerb{
This section collects several numerical experiments aimed at validating the novel elements of
the numerical solver proposed in this paper, and discussing their stability and convergence characteristics from 
an experimental point of view, as well as applications of the designed scheme to multi-layered, finely stratified two-fluid systems.
The reader interested in the application results for two-fluid multi-layered systems 
may jump directly to Section \ref{sec:multilayer}, where these problems are considered, skipping Sections \ref{sec:rp}--\ref{sec:bubble}.}

\reviewerb{In Sections \ref{sec:rp}--\ref{sec:bubble}, we aim to demonstrate that the proposed numerical scheme 
is robust and accurate, and experimentally assess its convergence properties. 
It also demonstrates the degree to which the scheme can handle sharp material interfaces and 
strong shockwaves and generates or does not generate oscillations as a result of their presence. 
Such experimental validation is necessary due to the fact that the scheme makes extensive use of a predictor stage 
based on a smooth wave equation, hence its behaviour in the presence of shockwaves cannot be taken for granted, 
and the filtering techniques for oscillations (Section \ref{sec:densityfiltering} and \ref{sec:pressurefiltering}) require validation 
due to their critical importance in preventing spurious oscillations and potential crashes.}

\subsection{Single material Riemann problems}\label{sec:rp}
We begin the validation of the proposed numerical methodology by showing the behaviour of the method with
varying mesh sizes and CFL values. To this purpose we select the simple Sod shock tube problem.

The method is then benchmarked against a much more stringent battery of tests, consisting of the Lax 
shock tube and the six problems given in the Toro's 2009 book \cite{torobook}.
These Riemann problems are designed to assess the accuracy and robustness of a numerical method with
respect to shockwaves (both fast and slow moving) and strong rarefactions.

The general form of the initial condition is composed of two constant states separated by a sharp 
discontinuity located at the spatial coordinate $x_\up{d}$
\begin{equation}
    \rho(t=0,\ x),\ u(t=0,\ x),\ p(t=0,\ x) = 
    \left\{
    \begin{aligned}
        \rho_\up{L},\ u_\up{L},\ p_\up{L} &\qquad\text{ if } x \le x_\up{d},\\
        \rho_\up{R},\ u_\up{R},\ p_\up{R} &\qquad\text{ if } x > x_\up{d},\\
    \end{aligned}
    \right.
\end{equation}
where the parametric quantities $\rho_\up{L},\ u_\up{L},\ p_\up{L},\ \rho_\up{R},\ u_\up{R},\ p_\up{R}$
are listed in Table \ref{tab:rp}, together with the final simulation times $t_\up{end}$.

\begin{table}[!tp]
\caption{Initial conditions for the single material Riemann problems.
The table lists the left and right states, the domain extrema $x_\up{L}$ and $x_\up{R}$, the initial position
of the discontinuity $x_\up{d}$, and the final time $t_\up{end}$.
In all cases, the material parameters are $\gamma_\up{L} = \gamma_\up{R} = 1.4$ and $\Pi_\up{L} = \Pi_\up{R} = 0.0$.
As the rest of the paper, the table uses SI units.
}
\label{tab:rp}
\begin{tabularx}{\textwidth}{lrRRRRRRR}
\toprule
              & Sod     & Lax     & Toro1     & Toro2     & Toro3       & Toro4         & Toro5          & Toro6 \\
\midrule
$\rho_\up{L}$ & $1.0$   & $0.445$ & $1.0$   & $1.0$  & $1.0$    & $5.99924$  & $1.0$       & $1.0$ \\
$u_\up{L}$    & $0.0$   & $0.698$ & $0.75$  & $-2.0$ & $0.0$    & $19.5975$  & $-19.59745$ & $2.0$ \\
$p_\up{L}$    & $1.0$   & $3.528$ & $1.0$   & $0.4$  & $1000.0$ & $460.894$  & $1000.0$    & $0.1$ \\
\midrule
$\rho_\up{R}$ & $0.125$ & $0.5$   & $0.125$ & $1.0$  & $1.0$    & $5.99242$  & $1.0$       & $1.0$ \\
$u_\up{R}$    & $0.0$   & $0.0$   & $0.0$   & $2.0$  & $0.0$    & $-6.19633$ & $-19.59745$ & $-2.0$\\
$p_\up{R}$    & $0.1$   & $0.571$ & $0.1$   & $0.4$  & $0.01$   & $46.095$   & $0.01$      & $0.1$ \\
\midrule
$x_\up{L}$  & $0.0$   & $0.0$   & $0.0$   & $0.0$  & $0.0$    & $-0.2$     & $0.0$       & $-2.1$\\
$x_\up{d}$    & $0.5$   & $0.5$   & $0.4$   & $0.5$  & $0.5$    & $0.5$      & $0.6$       & $0.0$   \\
$x_\up{R}$  & $1.0$   & $1.0$   & $1.0$   & $1.0$  & $1.0$    & $1.2$      & $1.0$       & $2.1$ \\
\midrule
$t_\up{end}$  & $0.25$  & $0.1$   & $0.2$   & $0.15$ & $0.012$  & $0.035$    & $0.012$     & $0.8$ \\
\bottomrule
\end{tabularx}
\end{table}

\begin{figure}[!b]
\centering
\includegraphics[width=0.99\textwidth]{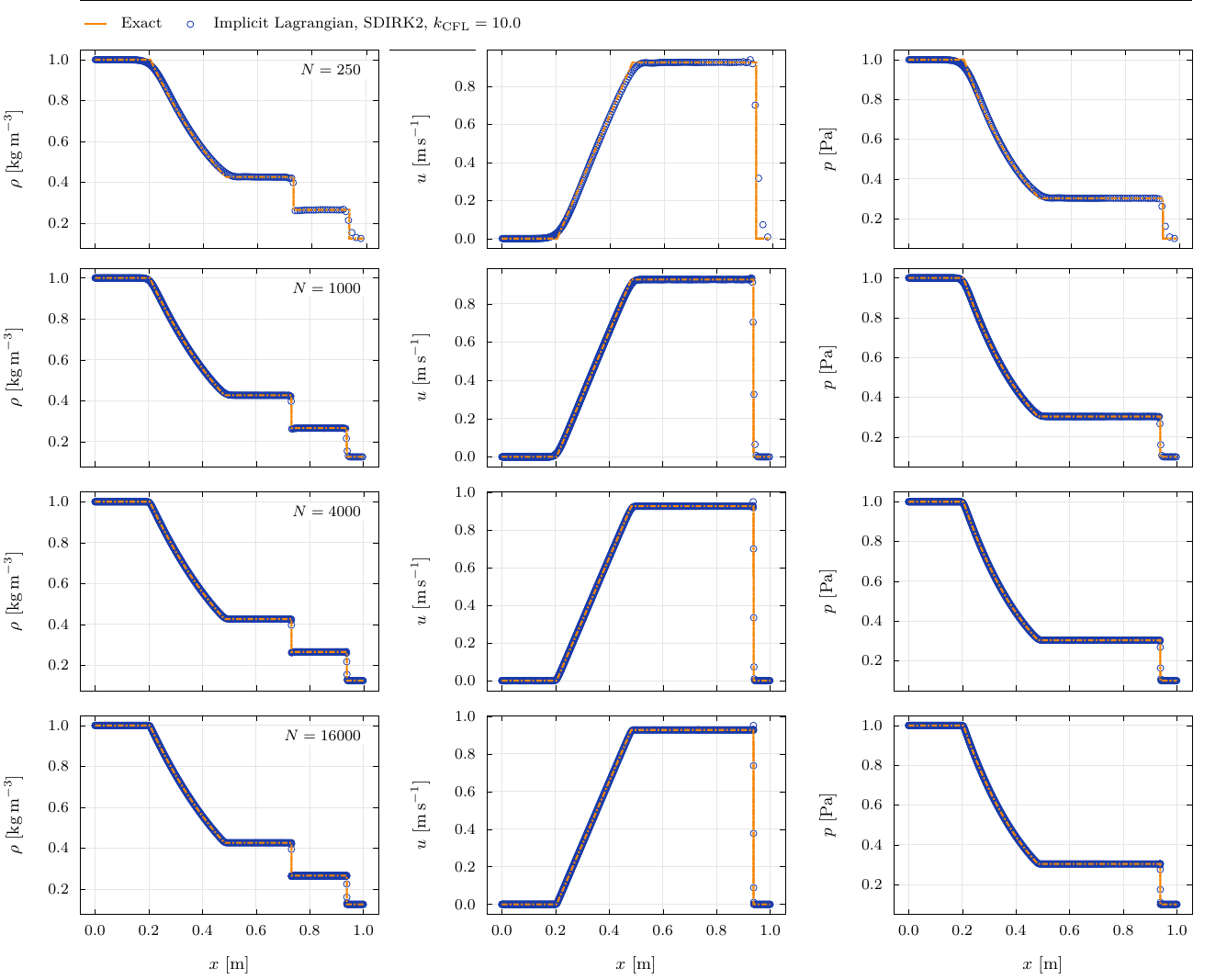}%
\caption{Behaviour of the proposed implicit Lagrangian numerical method at different uniform mesh sizes, from $N=250$ to $N=16000$ cells.
The number of points found at contact (0 to 1) and at shocks (2-3) 
is largely independent of mesh resolution, like it would be for an explicit Lagrangian scheme.
Major oscillations are absent, despite the relatively large CFL number ($\cfl=10.0$) and the presence of mild shockwaves.
The simulations employ a second order diagonally implicit \rk time integrator (SDIRK2).}
\label{fig:meshsweep}
\end{figure}

\begin{figure}[!b]
\centering
\includegraphics[width=0.99\textwidth]{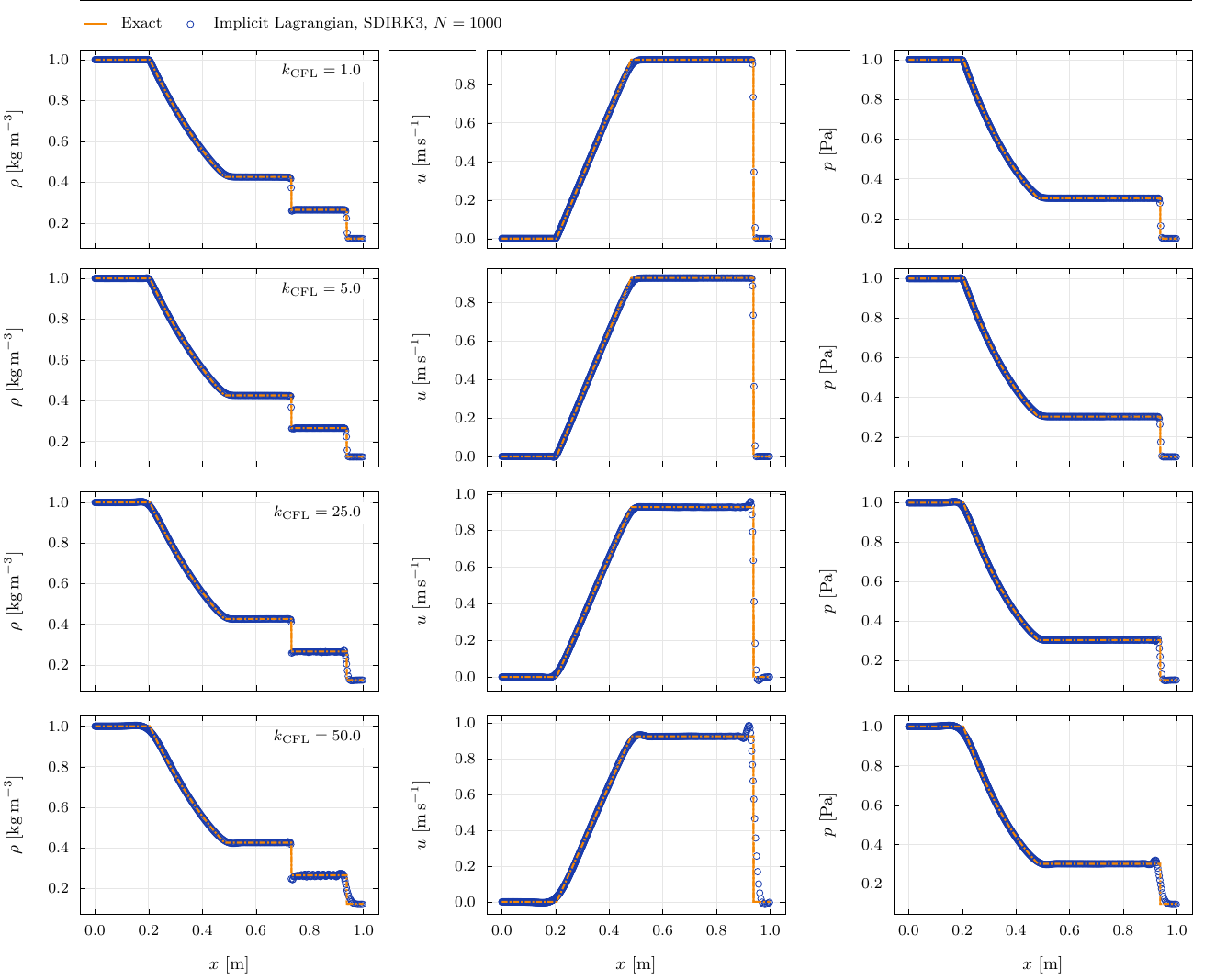}%
\caption{Behaviour of the proposed implicit Lagrangian numerical method at different CFL values, ranging from $\cfl = 1.0$ to $\cfl = 50.0$.
The first timestep of the simulations at $\cfl=25.0$ and $\cfl=50.0$ uses $\cflzero=\cfl/10$. Some oscillations 
and additional artificial diffusion appear for the higher CFL numbers
but the overall structure of the solution is well preserved. The simulations 
employ a third order diagonally implicit \rk time integrator (SDIRK3) and a constant-mass mesh of 
$N=1000$ cells.
The run at $\cfl=50.0$ is carried out in 16 timesteps.}
\label{fig:cflsweepfirst}
\end{figure}

\begin{figure}[!b]
\centering
\includegraphics[width=0.33\textwidth]{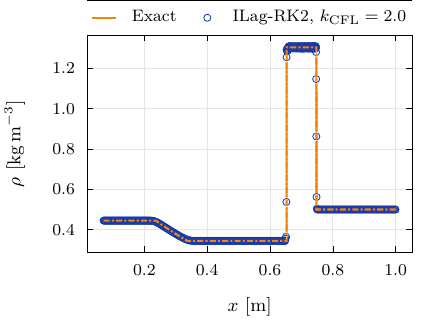}%
\includegraphics[width=0.33\textwidth]{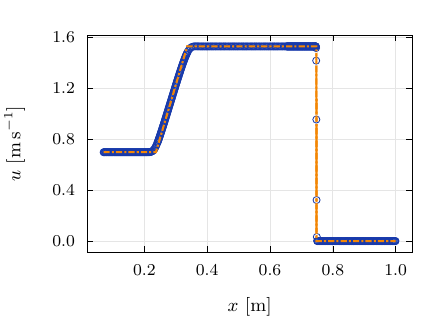}%
\includegraphics[width=0.33\textwidth]{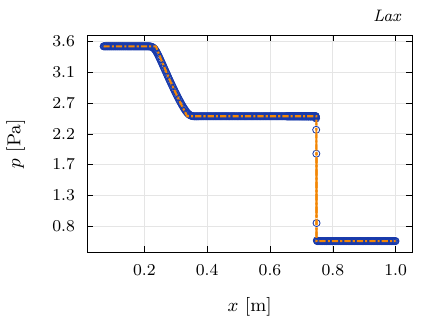}\\[4mm]
\includegraphics[width=0.33\textwidth]{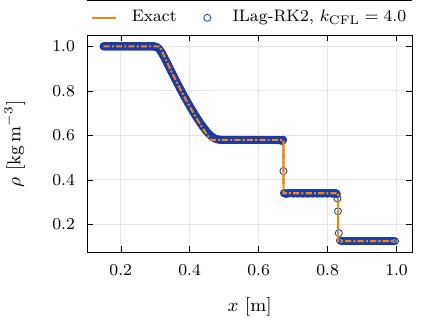}%
\includegraphics[width=0.33\textwidth]{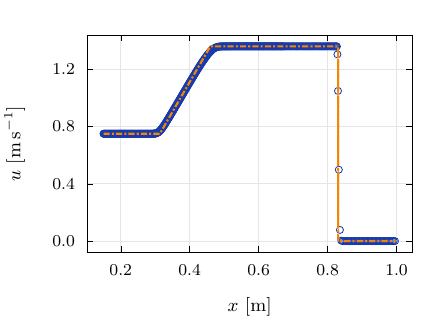}%
\includegraphics[width=0.33\textwidth]{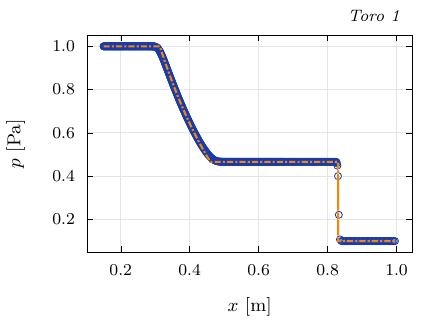}%
\caption{Numerical results for the Lax shock tube problem and for the modified Sod problem with transsonic rarefaction (Toro 1). 
The simulations employ a second order diagonally implicit \rk time integrator (SDIRK2) and a constant-mass mesh of 
$N=1000$ cells.}
\label{fig:rp1}
\end{figure}

\begin{figure}[!b]
\centering
\includegraphics[width=0.33\textwidth]{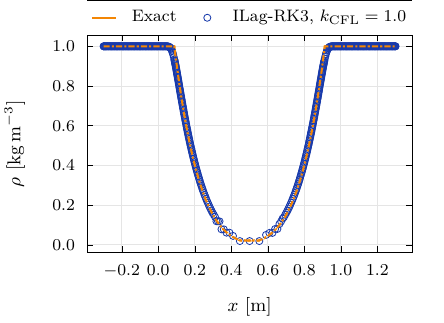}%
\includegraphics[width=0.33\textwidth]{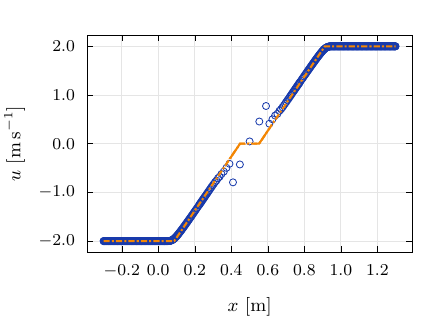}%
\includegraphics[width=0.33\textwidth]{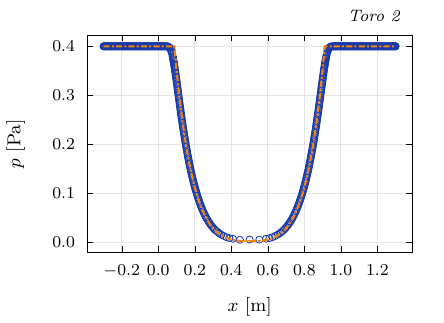}\\[4mm]
\includegraphics[width=0.33\textwidth]{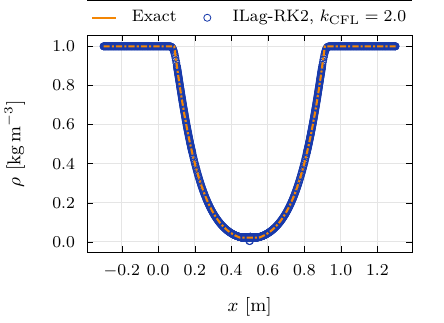}%
\includegraphics[width=0.33\textwidth]{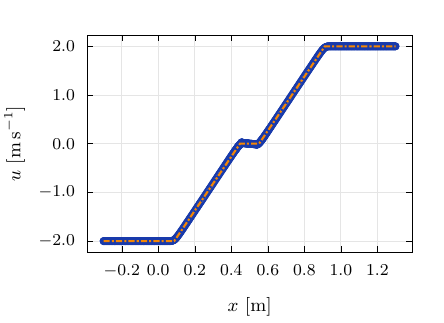}%
\includegraphics[width=0.33\textwidth]{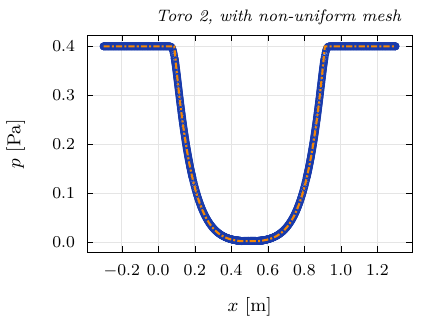}%
\caption{Numerical results for the second problem of Toro (also known as 123 problem) involving two strong rarefactions
generating a near-vacuum central state. The uniform mesh of $N=1000$ cells is expands with the rarefactions, 
causing inaccurate results in the central portion of the domain. In the bottom panels, we show the results obtained using
a finer non-uniform mesh of $N=2000$ cells, which accommodates the expansion effect.
The simulations employ a second order diagonally implicit \rk time integrator (SDIRK2).}
\label{fig:rp2}
\end{figure}

\begin{figure}[!b]
\centering
\includegraphics[width=0.33\textwidth]{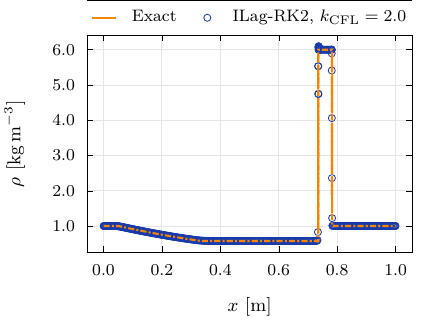}%
\includegraphics[width=0.33\textwidth]{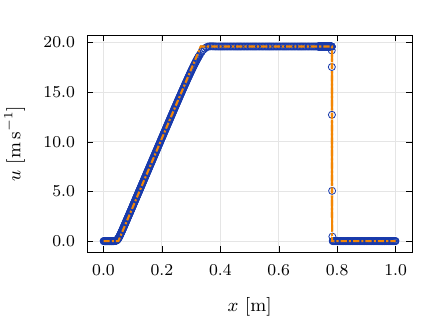}%
\includegraphics[width=0.33\textwidth]{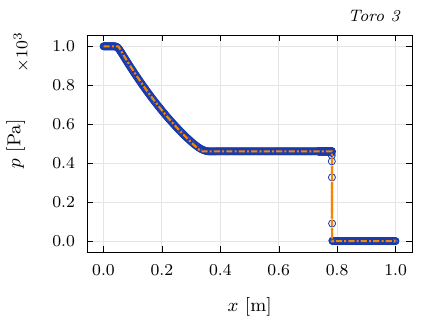}\\[4mm]
\includegraphics[width=0.33\textwidth]{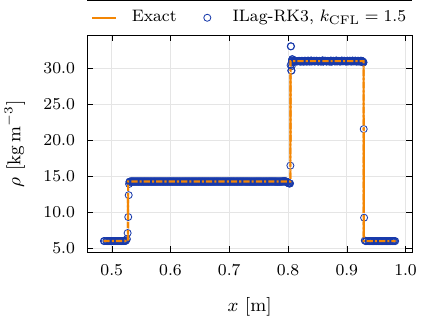}%
\includegraphics[width=0.33\textwidth]{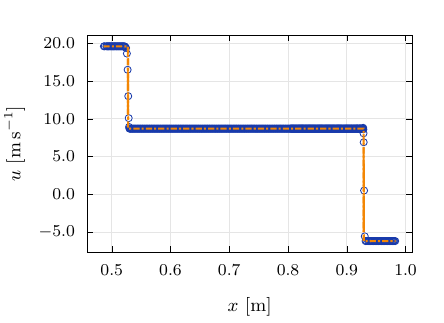}%
\includegraphics[width=0.33\textwidth]{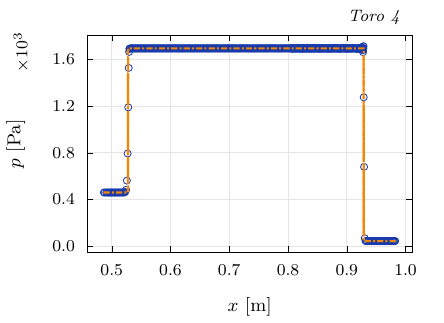}\\[4mm]
\includegraphics[width=0.33\textwidth]{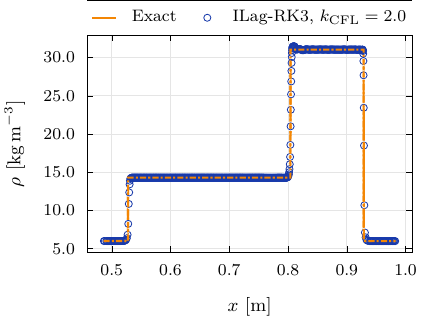}%
\includegraphics[width=0.33\textwidth]{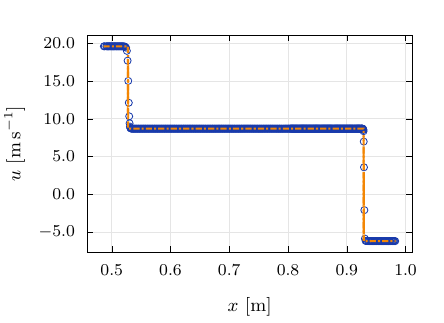}%
\includegraphics[width=0.33\textwidth]{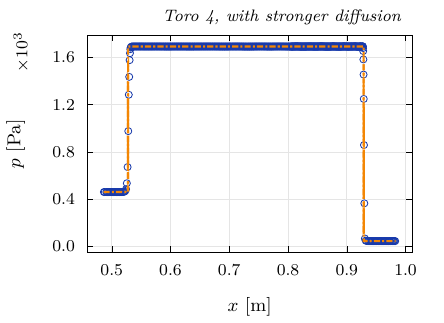}%
\caption{Numerical results for the third (first row) and fourth (second and third row) problems of Toro, obtained via 
the proposed implicit Lagrangian method with second or 
third order diagonally implicit \rk time integration (SDIRK2 or SDIRK3), labelled ILag-RK2 and ILag-RK3 respectively.
All simulations use a constant-mass mesh of $N=1000$ cells. For the Toro4 problem, the contact discontinuity 
presents some spurious oscillations and in the bottom panels we show the results obtained with a much more
diffusive variant of the regularisation operator for specific volume.}
\label{fig:rp3}
\end{figure}

\begin{figure}[!b]
\centering
\includegraphics[width=0.33\textwidth]{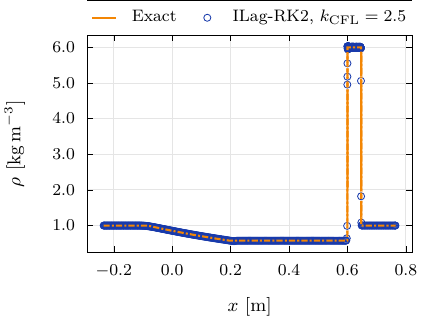}%
\includegraphics[width=0.33\textwidth]{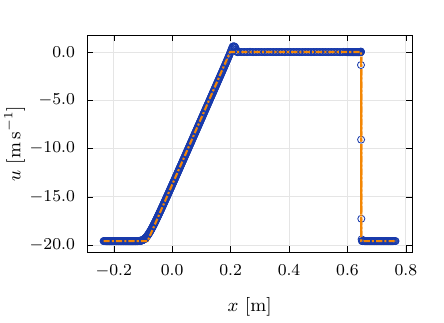}%
\includegraphics[width=0.33\textwidth]{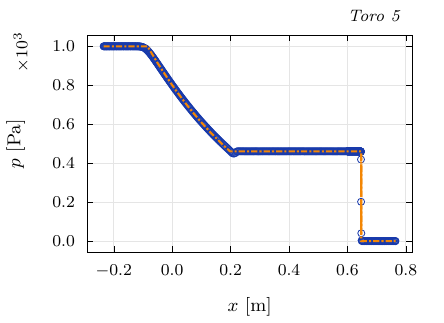}\\[4mm]
\includegraphics[width=0.33\textwidth]{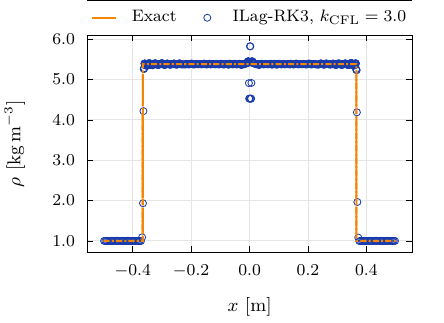}%
\includegraphics[width=0.33\textwidth]{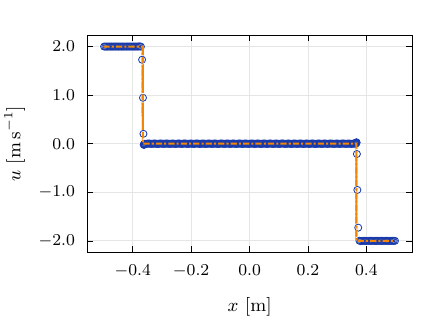}%
\includegraphics[width=0.33\textwidth]{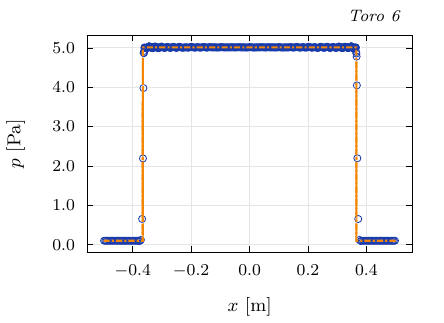}\\[4mm]
\includegraphics[width=0.33\textwidth]{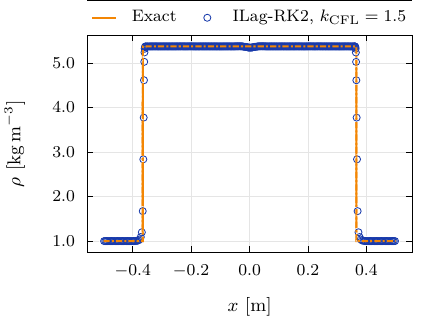}%
\includegraphics[width=0.33\textwidth]{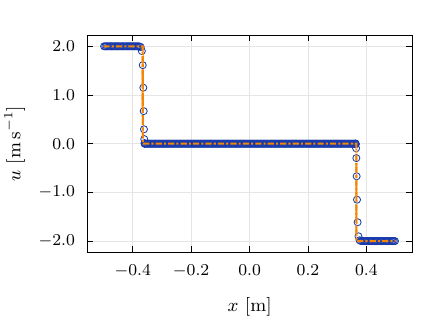}%
\includegraphics[width=0.33\textwidth]{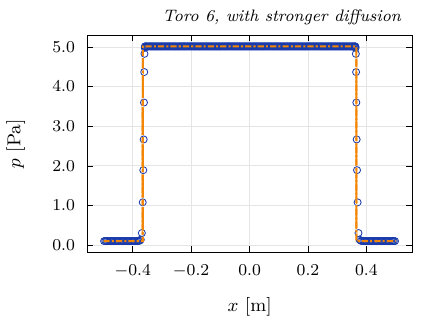}%
\caption{Numerical results for the fifth (first row) and sixth (second and third row) problems of Toro, obtained via 
the proposed implicit Lagrangian method with second or 
third order diagonally implicit \rk time integration (SDIRK2 or SDIRK3), labelled ILag-RK2 and ILag-RK3 respectively.
All simulations use a constant-mass mesh of $N=1000$ cells. For the Toro6 problem, oscillations due to 
well known wall heating/cooling effecs are present. The issue can be mitigated by adopting
a more aggressive numerical diffusion for the specific volume, as shown in the bottom panels.}
\label{fig:rp4}
\end{figure}

In Figure~\ref{fig:meshsweep} we show how the implicit Lagrangian method
here proposed behaves on the simple Sod shock tube problem at different 
mesh resolutions ranging from $N=250$ to $N=16000$ cells of equal mass.
The Courant number is fixed by setting $\cfl = 10.0$.
The main observations are that the contact is correctly captured without 
smearing (as expected from a Lagrangian scheme) and that the shockwave
is always spread over two mesh points despite the increased Courant number
with respect to what standard explicit time-`stepping would prescribe.

Analogously, in Figure~\ref{fig:cflsweepfirst} we fix the mesh
size to $N=1000$ uniform mass cells and let the Courant number vary between 
$\cfl = 1.0$ and $\cfl = 50.0$.
For the tests at $\cfl = 25.0$ and $\cfl = 50.0$, 
we ramp up the timestep starting from $\cflzero=\cfl/10$ and increasing linearly over ten timesteps. 
This means that 
the computation of Sod's problem in the last row of Figure~\ref{fig:cflsweepfirst} has been carried out
over 16 timesteps, for a mesh size of $N=1000$ uniform mass cells.
This is a slight variation on the practical recommendation given in Toro's book \cite{torobook}
about the choice of step size for the early stages of a numerical simulation involving discountinuous initial
data. In this case, instead of choosing a reduced size for the first, say, five timesteps, we gradually increase the timestep 
from a small initial value $\cflzero$ to the full value $\cfl$ in a linear fashion, so that
the time integrator does not see abrupt jumps in the step size.
As the timestep size increases, we see that both 
nonlinear waves are subject to additional numerical diffusion, with the contact remaining sharp and some 
oscillations can be observed especially for what concerns the velocity of the fluid at the shock front.
Nonetheless, the computations did not fail due to positivity violations.

Next, we carry out Lax's problem and the battery of benchmarks that can be found in Toro's 2009 book \cite{torobook}.
These test problems are aimed at assessing the robustness and accuracy of numerical methods for 
the compressible Euler equations in a variety of scenarios. The problems are labelled Toro1 to Toro6.
Toro1 is a modified version of Sod's problem, where a sonic point is present in the left moving rarefaction.
In Figure~\ref{fig:rp1} we report the results of the classic Lax problem and of the Toro1 sonic rarefaction problem, 
using the proposed implicit Lagrangian scheme with SDIRK2 time integration over a mesh of 1000 constant mass cells.
Contacts are sharp and no sonic glitch or other spurious artifacts are observed, with the Courant numbers 
($\cfl=2.0$ and $\cfl=4.0$ respectively) exceeding those prescribed by explicit schemes, despite the presence of
shockwaves.

In Figure~\ref{fig:rp2}, we plot the numerical results obtained for the so-called
123-problem (labelled Toro2 here). The problem features two strong symmetrical rarefaction waves giving 
rise to a very low density central state. 
{In the upper panels we present the numerical results obtained with a uniform mesh in the mass coordinate
with 1000 cells. It appears that the low density central region is poorly resolved, therefore the numerical solution is not accurate there. The lower panels are obtained adopting a slightly finer grid of 2000 elements which are
initially denser towards the center of the computational domain, thus providing a much more accurate solution. }

The robustness of the proposed discretization is then tested on the Toro3 and Toro4 Riemann problems.
In the first row of Figure~\ref{fig:rp3} we can observe that the contact wave of the Toro3 problem 
is affected by some small spurious oscillations at $\cfl=2$, and even more so in the second row 
for what concerns the Toro4 problem at $\cfl=1.5$. In this latter test, the oscillations are pronounced
enough to warrant the use of stronger numerical dissipation on the density/specific volume, which 
successfully controls the size of the oscillations, but at the price of losing exactly sharp contact
discontinuities, as can be observed in the third row of Figure~\ref{fig:rp3}.

Similar observations can be given with regards to the results of the Toro5 and Toro6 tests shown in 
Figure~\ref{fig:rp4}: at $\cfl=2.5$, some mild oscillations are deteriorating the contact wave and the 
base of the rarefaction in the Toro5 problem, and strong density artifacts can be 
observed in the center of Toro6 problem at $\cfl=3.0$, which is nothing but two colliding shocks, where usually wall
heating effects are observed. In \cite{casullidumbser}, a semi-implicit scheme running at $\cfl=0.5$
was shown to produce a wall cooling effect instead. As done for the Toro4 problem, additional numerical 
diffusion in the density, and a reduction of the Courant number to $\cfl=1.5$ can reduce these oscillations, 
while of course at the same time smearing all waves a little bit more.

\subsection{Convergence results}\label{sec:resultsconvergence}

\begin{table}[!t]
\caption{
Convergence results for the proposed implicit method, using the first order diffusion operator for the total energy.
In the table are shown the $L^1$, $L^2$, and $L^\infty$ norms of the error on the specific total energy $E$ and the corresponding convergence rates, 
for two different implicit time integration schemes (SDIRK2 and SDIRK3) and with varying choices 
of $k_\up{CFL}$. Note that the convergence rates for the higher CFL numbers approach 2, despite the underlying spatial discretization being of first order, 
since at high CFL numbers the temporal discretisation error was found to be the dominant one in this case.}
\label{tab:convergencefirst}
\begin{tabularx}{\textwidth}{lrRRRRRR}
\toprule
                   & $N$       & $\mathcal{E}_{L^1}$   & $\mathcal{E}_{L^2}$   & $\mathcal{E}_{L^\infty}$ & $\mathcal{O}_{L^1}$ & $\mathcal{O}_{L^2}$ & $\mathcal{O}_{L^\infty}$ \\

\midrule
SDIRK3             & $    100$ & $4.32 \times 10^{-3}$ & $6.06 \times 10^{-3}$ & $2.26 \times 10^{-2}$    &  $-$      &    $-$      &    $-$    \\ 
$k_\up{CFL}=1.0$   & $    400$ & $9.62 \times 10^{-4}$ & $1.70 \times 10^{-3}$ & $8.52 \times 10^{-3}$    &  $1.08$   &    $0.92$   &    $0.70$ \\ 
                   & $   1\,600$ & $2.34 \times 10^{-4}$ & $4.44 \times 10^{-4}$ & $2.59 \times 10^{-3}$    &  $1.02$   &    $0.97$   &    $0.86$ \\ 
                   & $   6\,400$ & $5.80 \times 10^{-5}$ & $1.13 \times 10^{-4}$ & $6.79 \times 10^{-4}$    &  $1.01$   &    $0.99$   &    $0.96$ \\ 
                   & $  25\,600$ & $1.45 \times 10^{-5}$ & $2.83 \times 10^{-5}$ & $1.72 \times 10^{-4}$    &  $1.00$   &    $1.00$   &    $0.99$ \\[1mm]
                   \midrule
SDIRK2             & $    100$ & $4.03 \times 10^{-3}$ & $5.44 \times 10^{-3}$ & $1.95 \times 10^{-2}$    &  $-$      &    $-$      &    $-$    \\ 
$k_\up{CFL}=1.0$   & $    400$ & $9.34 \times 10^{-4}$ & $1.62 \times 10^{-3}$ & $7.91 \times 10^{-3}$    &  $1.05$   &    $0.88$   &    $0.65$ \\ 
                   & $   1\,600$ & $2.32 \times 10^{-4}$ & $4.38 \times 10^{-4}$ & $2.53 \times 10^{-3}$    &  $1.01$   &    $0.94$   &    $0.82$ \\ 
                   & $   6\,400$ & $5.79 \times 10^{-5}$ & $1.12 \times 10^{-4}$ & $6.75 \times 10^{-4}$    &  $1.00$   &    $0.98$   &    $0.95$ \\ 
                   & $  25\,600$ & $1.45 \times 10^{-5}$ & $2.82 \times 10^{-5}$ & $1.72 \times 10^{-4}$    &  $1.00$   &    $1.00$   &    $0.99$ \\[1mm]
                   \midrule
SDIRK2             & $    100$ & $2.30 \times 10^{-2}$ & $2.55 \times 10^{-2}$ & $4.92 \times 10^{-2}$    &  $-$      &    $-$      &    $-$    \\ 
$k_\up{CFL}=20.0$  & $    400$ & $7.21 \times 10^{-3}$ & $9.94 \times 10^{-3}$ & $2.94 \times 10^{-2}$    &  $0.84$   &    $0.68$   &    $0.37$ \\ 
                   & $   1\,600$ & $8.29 \times 10^{-4}$ & $1.68 \times 10^{-3}$ & $8.73 \times 10^{-3}$    &  $1.56$   &    $1.28$   &    $0.88$ \\ 
                   & $   6\,400$ & $6.10 \times 10^{-5}$ & $1.50 \times 10^{-4}$ & $9.77 \times 10^{-4}$    &  $1.88$   &    $1.75$   &    $1.58$ \\ 
                   & $  25\,600$ & $1.17 \times 10^{-5}$ & $2.25 \times 10^{-5}$ & $1.43 \times 10^{-4}$    &  $1.19$   &    $1.36$   &    $1.39$ \\ 
                   & $ 102\,400$ & $3.44 \times 10^{-6}$ & $6.59 \times 10^{-6}$ & $3.95 \times 10^{-5}$    &  $0.88$   &    $0.89$   &    $0.93$ \\[1mm]
                   \midrule
SDIRK2             & $    800$ & $1.69 \times 10^{-2}$ & $2.08 \times 10^{-2}$ & $4.58 \times 10^{-2}$    &  $-$      &    $-$      &    $-$    \\ 
$k_\up{CFL}=100.0$ & $   3\,200$ & $3.87 \times 10^{-3}$ & $6.11 \times 10^{-3}$ & $2.25 \times 10^{-2}$    &  $1.06$   &    $0.88$   &    $0.51$ \\ 
                   & $  12\,800$ & $3.87 \times 10^{-4}$ & $9.08 \times 10^{-4}$ & $5.54 \times 10^{-3}$    &  $1.66$   &    $1.38$   &    $1.01$ \\ 
                   & $  51\,200$ & $2.55 \times 10^{-5}$ & $7.45 \times 10^{-5}$ & $6.03 \times 10^{-4}$    &  $1.96$   &    $1.80$   &    $1.60$ \\ 
                   & $ 204\,800$ & $1.66 \times 10^{-6}$ & $4.15 \times 10^{-6}$ & $2.75 \times 10^{-5}$    &  $1.97$   &    $2.08$   &    $2.23$ \\[1mm]
\bottomrule 
\end{tabularx}
\end{table}

\begin{table}[!t]
\caption{Convergence results for the proposed implicit method, using the second order diffusion operator for the total energy.
In the table are shown the $L^1$, $L^2$, and $L^\infty$ norms of the error on the specific total energy $E$ and the corresponding convergence rates, 
for two different implicit time integration schemes (SDIRK2 and SDIRK3) and with varying choices of $k_\up{CFL}$.
The error norms relative to  
standard explicit Lagrangian scheme with exact Riemann solver and second order generalized minmod reconstruction and 
SSPRK2 or SSPRK3 
time integration ($k_\up{CFL} = 0.8$) are shown for comparison.
}
\label{tab:convergencesecond}
\begin{tabularx}{\textwidth}{lrRRRrrr}
\toprule
                     & $N$      & $\mathcal{E}_{L^1}$   & $\mathcal{E}_{L^2}$   & $\mathcal{E}_{L^\infty}$ & $\mathcal{O}_{L^1}$ & $\mathcal{O}_{L^2}$ & $\mathcal{O}_{L^\infty}$ \\

\midrule
SDIRK3               & $100$    & $1.12 \times 10^{-3}$ & $1.62 \times 10^{-3}$ & $5.49 \times 10^{-3}$    & $-$                 & $-$                 & $-$    \\
$k_\up{CFL}=1.0$     & $400$    & $8.19 \times 10^{-5}$ & $1.43 \times 10^{-4}$ & $6.57 \times 10^{-4}$    & $1.89$              & $1.75$              & $1.53$ \\
                     & $1\,600$   & $5.94 \times 10^{-6}$ & $1.26 \times 10^{-5}$ & $1.02 \times 10^{-4}$    & $1.89$              & $1.75$              & $1.34$ \\
                     & $6\,400$   & $3.88 \times 10^{-7}$ & $9.34 \times 10^{-7}$ & $9.08 \times 10^{-6}$    & $1.97$              & $1.88$              & $1.74$ \\
                     & $25\,600$  & $2.42 \times 10^{-8}$ & $6.25 \times 10^{-8}$ & $6.94 \times 10^{-7}$    & $2.00$              & $1.95$              & $1.85$ \\
Ref. Explicit SSPRK3 & $25\,600$  & $3.25 \times 10^{-8}$ & $9.32 \times 10^{-8}$ & $8.20 \times 10^{-7}$    & $-$                 & $-$                 & $-$    \\[1mm]
\midrule
SDIRK2               & $100$    & $1.13 \times 10^{-3}$ & $1.48 \times 10^{-3}$ & $3.34 \times 10^{-3}$    & $-$                 & $-$                 & $-$    \\
$k_\up{CFL}=1.0$     & $400$    & $9.52 \times 10^{-5}$ & $1.99 \times 10^{-4}$ & $1.30 \times 10^{-3}$    & $1.79$              & $1.45$              & $0.68$ \\
                     & $1\,600$   & $6.61 \times 10^{-6}$ & $1.77 \times 10^{-5}$ & $1.67 \times 10^{-4}$    & $1.92$              & $1.75$              & $1.48$ \\
                     & $6\,400$   & $4.32 \times 10^{-7}$ & $1.22 \times 10^{-6}$ & $1.31 \times 10^{-5}$    & $1.97$              & $1.93$              & $1.84$ \\
                     & $25\,600$  & $2.71 \times 10^{-8}$ & $8.04 \times 10^{-8}$ & $9.06 \times 10^{-7}$    & $2.00$              & $1.96$              & $1.93$ \\
Ref. Explicit SSPRK2 & $25\,600$  & $5.95 \times 10^{-8}$ & $1.81 \times 10^{-7}$ & $1.59 \times 10^{-6}$    & $-$                 & $-$                 & $-$    \\[1mm]
\midrule
SDIRK2               & $100$    & $2.52 \times 10^{-2}$ & $2.85 \times 10^{-2}$ & $5.36 \times 10^{-2}$    & $-$                 & $-$                 & $-$    \\
$k_\up{CFL}=20.0$    & $400$    & $8.02 \times 10^{-3}$ & $1.11 \times 10^{-2}$ & $3.23 \times 10^{-2}$    & $0.83$              & $0.68$              & $0.36$ \\
                     & $1\,600$   & $9.65 \times 10^{-4}$ & $2.01 \times 10^{-3}$ & $1.04 \times 10^{-2}$    & $1.53$              & $1.23$              & $0.82$ \\
                     & $6\,400$   & $7.10 \times 10^{-5}$ & $2.01 \times 10^{-4}$ & $1.53 \times 10^{-3}$    & $1.88$              & $1.66$              & $1.38$ \\
                     & $25\,600$  & $4.47 \times 10^{-6}$ & $1.36 \times 10^{-5}$ & $1.16 \times 10^{-4}$    & $1.99$              & $1.94$              & $1.86$ \\
                     & $102\,400$ & $2.81 \times 10^{-7}$ & $8.56 \times 10^{-7}$ & $7.32 \times 10^{-6}$    & $2.00$              & $2.00$              & $1.99$ \\[1mm]
\midrule
SDIRK2               & $800$    & $1.74 \times 10^{-2}$ & $2.13 \times 10^{-2}$ & $4.66 \times 10^{-2}$    & $-$                 & $-$                 & $-$    \\
$k_\up{CFL}=100.0$   & $3\,200$   & $3.97 \times 10^{-3}$ & $6.28 \times 10^{-3}$ & $2.30 \times 10^{-2}$    & $1.07$              & $0.88$              & $0.51$ \\
                     & $12\,800$  & $4.03 \times 10^{-4}$ & $9.51 \times 10^{-4}$ & $5.79 \times 10^{-3}$    & $1.65$              & $1.36$              & $0.99$ \\
                     & $51\,200$  & $2.79 \times 10^{-5}$ & $8.31 \times 10^{-5}$ & $6.75 \times 10^{-4}$    & $1.93$              & $1.76$              & $1.55$ \\
                     & $204\,800$ & $1.75 \times 10^{-6}$ & $5.32 \times 10^{-6}$ & $4.55 \times 10^{-5}$    & $2.00$              & $1.98$              & $1.95$ \\
\bottomrule
\end{tabularx}
\end{table}

\begin{figure}[!b]
\centering
\includegraphics[width=0.99\textwidth]{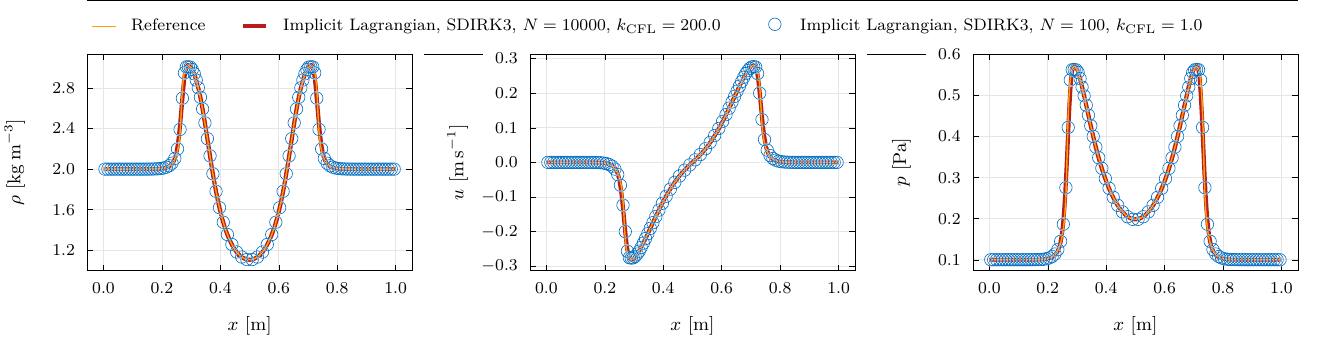}%
\caption{Numerical solutions at time $t=0.2\,\up{s}$ for the smooth convergence benchmark. In the three panels 
are shown the profiles of density, velocity, and pressure computed by the proposed implicit Lagrangian method
using two different choices of mesh and timestep size, compared with a reference solution obtained with an 
Eulerian \muscl scheme on a fine mesh of 200\,000 uniform cells. The reference solution is accurately reproduced with 
$N=100$ uniform mass cells, and on finer grids with high CFL number ($N=10000$ cells, $\cfl = \cflzero = 200$, 10 timesteps total).}
\label{fig:convergenceslices}
\end{figure}
We validate the convergence properties and order of accuracy of the proposed 
method by comparing its results with those of a reference solution obtained on a very fine grid 
by means of a standard explicit Finite Volume method. The test problem is chosen in such a way that
the solution can be considered smooth until the final time, but also making sure that the effects
of the nonlinearity of the governing equations are highlighted by the test, hence making sure
that the presented convergence results are not only valid for small amplitude waves.
\begin{equation}
    \rho(t=0,\ x) = \rho_0,\qquad u(t=0,\ x) = u_0,\qquad p(t = 0,\ x) = p_0\,\left\{1 + 10\, \exp\left[-(x - x_0)^2/\ell^2\right]\right\},
\end{equation}
with $\rho_0 = 2.0\,\up{kg\,m^{-3}}$, $u_0 = 0.0\,\up{m\,s^{-1}}$, $p_0 = 0.1\,\up{Pa}$, $\ell = 0.075\,\up{m}$,
constant (single-material) parameters for the stiffened gas equation of state
$\gamma_1 = \gamma_2 = \gamma = 1.4$, $\Pi_1 = \Pi_2 = \Pi = 0.5\,\up{Pa}$. Note that, while we express these quantities with SI units, 
the parameters and reference states for this problem have been chosen in such a way to highlight the nonlinear behaviour of the PDE system, in particular, obtaining significant deviations from the quasi-linear regime of small perturbations and obtain nonlinear wave steepening.

The simulation domain is $x \in {[0,\ 1]}$ and the final time
$t_\up{end} = 0.2\,\up{s}$ is chosen such that the data can be seen as smooth when sampled at the 
considered mesh resolutions (ranging between $N = 100$ and $N = 204800$ cells), while at the same time 
significant wave steepening effects are observed.

\begin{figure}[!b]
\centering
\includegraphics[width=0.495\textwidth]{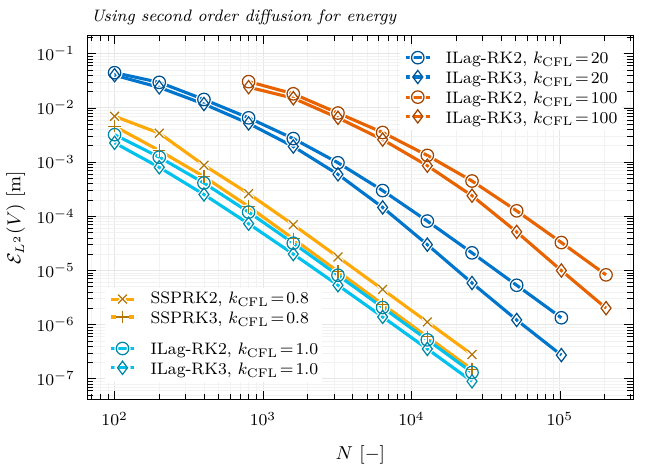}%
\includegraphics[width=0.495\textwidth]{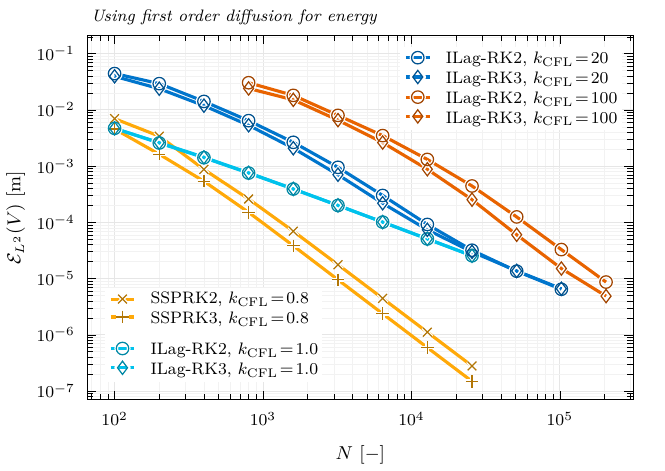}%
\caption{Bilogarithmic plots of the $L^2$ error norms given by the proposed implicit Lagrangian method
for the specific volume variable $V$ using different time integration schemes (second and third order SDIRK, 
labelled ILag-RK2 and ILag-RK3 respectively) and different timestep sizes (from $\cfl=\cflzero=1$ to $\cfl=\cflzero=100$).
The error norms given by a standard explicit second order Lagrangian scheme using SSP \rk time integration (SSPRK2 and SSPRK3) 
are included for comparison.
On the left, the implicit schemes use the second order energy diffusion operator and 
second order convergence rates are obtained at all CFL numbers for sufficiently fine meshes. On coarse meshes with
high CFL number, lower convergence rates are observed because the simulations are completed in too few timesteps:
for example at $\cfl=100$ the simulation with 800 grid points consists of two timesteps. At unit CFL, the errors 
of the implicit Lagrangian method are lower than those of the explicit reference scheme.
In the right panel, we show the results obtained without the low-diffusion correction on energy: at low CFL, or for 
very coarse meshes the order of convergence degenerates to 1 (as expected), but for higher CFL numbers
such an effect is negligible and observed errors are similar to those obtained with the second order
numerical diffusion operator.
}
\label{fig:convergence}
\end{figure}

The results of Figure~\ref{fig:convergenceslices} can be used to quantify the wave steepening effects 
and how the solution is captured at $\cfl=1.0$ with $N=100$ cells and at $\cfl=200.0$ with $N=10000$ cells, 
in which case the solution is marched to the final time $t_\up{end} = 0.2\,\up{s}$ in 10 timesteps.
We also report the corresponding error norms and experimental orders of convergence in Table~\ref{tab:convergencefirst}, 
where no piecewise linear reconstruction was used, yielding a scheme of first order spatial accuracy, 
and in Table~\ref{tab:convergencesecond}, where a minmod-limited piecewise linear reconstruction is 
employed in order to compute the artificial diffusion for energy.
A graphical rendition of the same data can be found in Figure~\ref{fig:convergence}.
When the piecewise linear reconstruction is used for the computation 
of artificial diffusion on energy, second order of convergence
is experimentally verified, given sufficiently fine meshes that 
allow to resolve all the features of the flow as smooth.
If the data reconstruction is not used, one would expect to achieve
only first order of accuracy, and indeed at $\cfl=1$, this is what we find.
When using very large timesteps ($\cfl=20$, $\cfl=100$) however, convergence
slopes closer to 2 are found.
Note that in Table~\ref{tab:convergencefirst} and Figure~\ref{fig:convergence}, 
at sufficiently high CFL numbers, the scheme converges with accuracy higher than 1, despite the fact that
the spatial discretization in these cases is only first order accurate. This is simply due to the fact that with large timesteps, 
the temporal discretization error has become the dominant one in this test problem.
Moreover, it is often the case that application to systems with many layers end up using very large cell counts simply to represent 
with some detail the layered structure of the medium, precisely leading to situations in which a significant benefit can be gained from
the use of higher order temporal integrators.

Based on the results of this study, we also carried out an efficiency study
to estimate the computational expense of the proposed method. 
Generally, when using implicit schemes, one expects to spend 
more computational time per Runge--Kutta stage with respect to an
explicit method, but such higher per-stage cost is offset by the fact that
the implicit scheme can carry out much larger timesteps.
In our case, we compare the proposed implicit Lagrangian method
with a reference explicit Lagrangian method, which makes use of 
low dissipation reconstruction (a low-diffusion 
generalized MinMod limiter, see \cite{kurganovgminmod}), together
with the exact Riemann solver, in order to match the
resolution of the proposed scheme.
In our naive Fortran implementation, which does \textit{not}
explicitly make use of vector instructions, 
the implicit method requires less computational time than the 
explicit one even when the timestep size is matched.
This is due to the fact that the explicit method
relies on the rather expensive exact Riemann solver, while
the implicit scheme is indeed iterative but its implicit 
part is solved by means of the Thomas algorithm and for 
this reason does not carry the same cost usually associated
with implicit methods. Extension to (nonlinear) equations of 
state more complex than the stiffened gas EOS can be 
carried out with the aid of the method by Casulli and Zanolli 
\cite{casulli2012}.
The use of the exact Riemann solver as a comparison is
justified by the fact that for this smooth test problem, 
the reference explicit scheme at $\cfl=0.8$ yields
numerical errors which are still slightly higher than those obtained
with the Implicit Lagrangian method in the Courant 
range between 
$\cfl=0.5$ and $\cfl=2.0$. 
The result shows that the method 
can achieve minimum error not only at $\cfl=1$ but in 
a range around it, and more over such minimum error is 
lower than that given by the reference scheme using a
little-dissipative reconstruction limiter and the exact 
Riemann solver. 

The implicit scheme is then favorable as a practical alternative
to the explicit one, given that no extra computational cost
is associated with its use, and at the same time the user
can choose to use much larger timesteps if as needed.
In light of the results of Section~\ref{sec:rp}, we can 
however see that an explicit method might still be 
preferable for problems involving very strong shockwaves.

\begin{figure}[!b]
\centering
\includegraphics[width=0.495\textwidth]{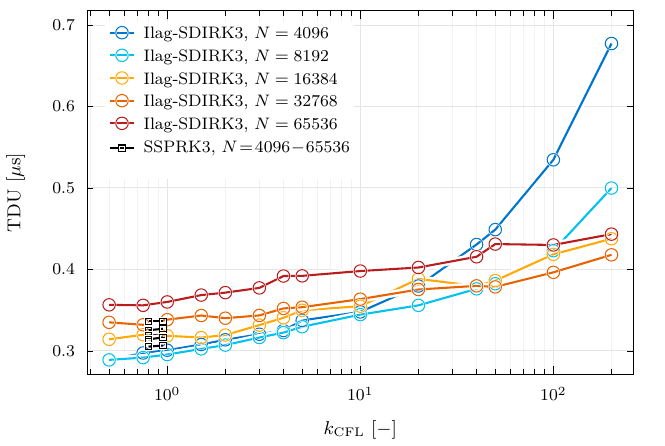}%
\includegraphics[width=0.495\textwidth]{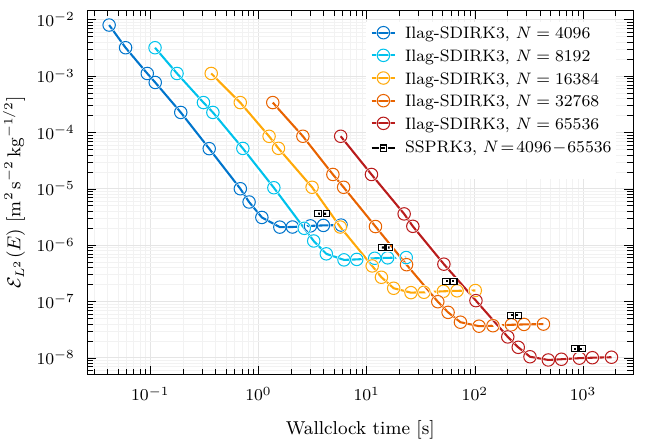}%
\caption{Performance measurements. In the left panel, computational times per degree of freedom update per Runge--Kutta stage (TDU) as a 
function of the CFL number, for mesh sizes ranging from $N=4096$ to $N=65536$, compared with the 
corresponding timing of the explicit scheme. While the CFL number ranges more than two orders of magnitude, the 
time per stage grows at most by a factor of 3, and is almost constant for sufficiently fine meshes. In the right 
panel, the $L^2$ error norms for the total energy are plotted as a function of the wallclock time, 
which depends on $\cfl$, ranging here between $\cfl=0.5$ (more timestep, longer 
wallclock time) and $\cfl=200$ (fewer timeteps, shorter computations). 
 For all mesh sizes, the proposed method achieves better results
than the explicit reference scheme using an exact Riemann solver (at $\cfl=0.8$ and $\cfl=0.95$). The plots also show that the same error norms can be obtained with significantly shorter computational times, selecting for example $\cfl=2.0$.}
\label{fig:performance}
\end{figure}

\subsection{Two material Riemann problems}\label{sec:abgrallkarni}
A series of classic two-material Riemann problems are used to evaluate the
performance of the proposed implicit Lagrangian method. 
The results here shown extend the validation of the proposed
scheme to systems where the material properties are allowed to jump
across cell interfaces. To this purpose we reproduce the results
of Abgrall and Karni, which involve material interfaces with 
significant differences in density and compressibility between adjacent layers.

\begin{table}[!tp]
\caption{Initial conditions for the two-material Riemann problems.
The table lists the left and right states, the material parameters for the two regions and the final time $t_\up{end}$.
In all cases, the domain extrema are $x_\up{L} = 0.0$ and $x_\up{R} = 1.0$ and the initial position of the 
discontinuity is $x_\up{d} = 0.5$. As the rest of the paper, the table uses SI units.}
\label{tab:rpmultimaterial}
\begin{tabularx}{\textwidth}{lRRRRRRrrrrr}
\toprule
   & $\rho_\up{L}$ & $u_\up{L}$ & $p_\up{L}$ & $\rho_\up{R}$ & $u_\up{R}$ & $p_\up{R}$  & $\gamma_\up{L}$  & $\Pi_\up{L}$ & $\gamma_\up{R}$ & $\Pi_\up{R}$  & $t_\up{end}$ \\ 
\midrule
AK1 & $1.0    $ & $0.0        $ & $1.0      $ & $0.125  $ & $0.0        $ & $0.1   $     & $1.6$ & $0.0$               & $1.2  $ & $0.0$ & $0.25$\\
AK2 & $1.0    $ & $0.0        $ & $500.0    $ & $1.0    $ & $0.0        $ & $0.2   $     & $1.4$ & $0.0$               & $1.6  $ & $0.0$ & $0.015$\\
AK3 & $1000.0 $ & $0.0        $ & $10^9   $ & $50.0   $ & $0.0        $   & $10^5  $     & $4.4$ & $6\!\times\!10^8$ & $1.4  $ & $0.0$ & $1.5\!\times\!10^{-4}$ \\
\bottomrule
\end{tabularx}
\end{table}

The initial conditions are simply given by
\begin{equation}
    \rho(t=0,\ x),\ u(t=0,\ x),\ p(t=0,\ x),\ \gamma(x),\ \Pi(x) = 
    \left\{
    \begin{aligned}
        \rho_\up{L},\ u_\up{L},\ p_\up{L},\ \gamma_\up{L},\ \Pi_\up{L} &\qquad\text{ if } x \le x_\up{d},\\
        \rho_\up{R},\ u_\up{R},\ p_\up{R},\ \gamma_\up{R},\ \Pi_\up{R} &\qquad\text{ if } x > x_\up{d},\\
    \end{aligned}
    \right.
\end{equation}
for three separate choices of initial data and material parameters, as listed 
in Table~\ref{tab:rpmultimaterial}.
The first test, labelled AK1, is a variant of Sod's problem and is tested at $\cfl = 10.0$, while the second (AK2) is a more violent shock tube, for which we ran 
the simulation at $\cfl = 2.0$. The third test (AK3) features a strong pressure jump and the CFL number was set to $\cfl=4.0$.

\begin{figure}[!b]
\centering
\includegraphics[width=0.33\textwidth]{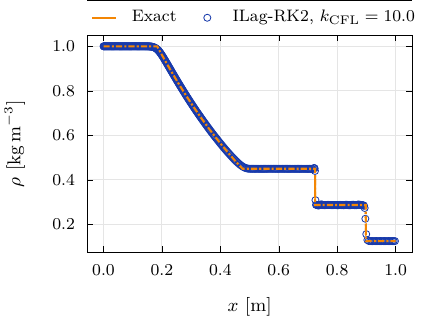}%
\includegraphics[width=0.33\textwidth]{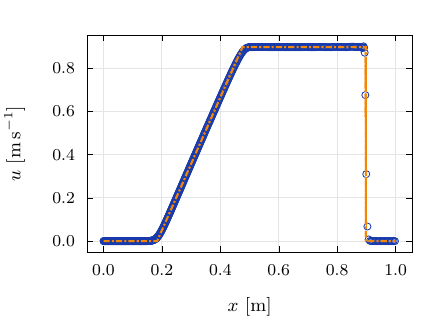}%
\includegraphics[width=0.33\textwidth]{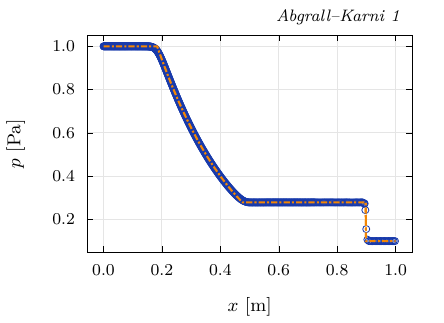}\\[4mm]
\includegraphics[width=0.33\textwidth]{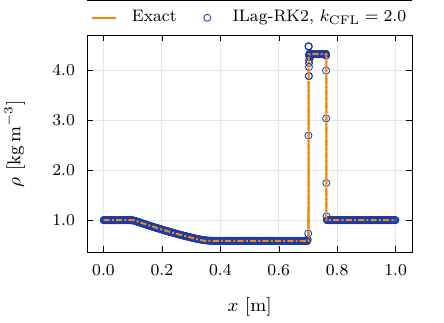}%
\includegraphics[width=0.33\textwidth]{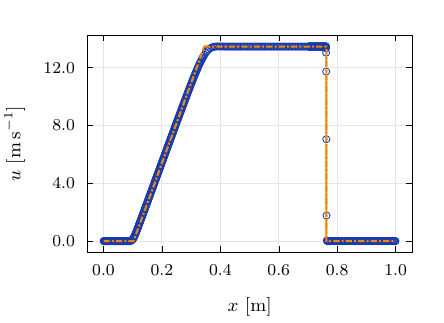}%
\includegraphics[width=0.33\textwidth]{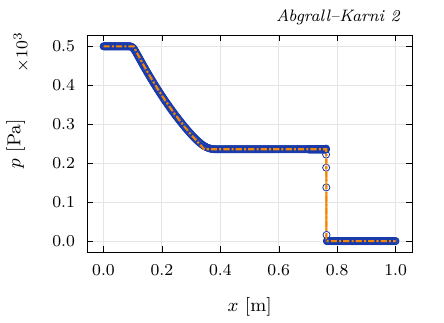}\\[4mm]
\includegraphics[width=0.33\textwidth]{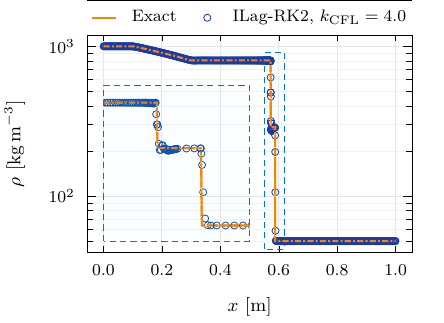}%
\includegraphics[width=0.33\textwidth]{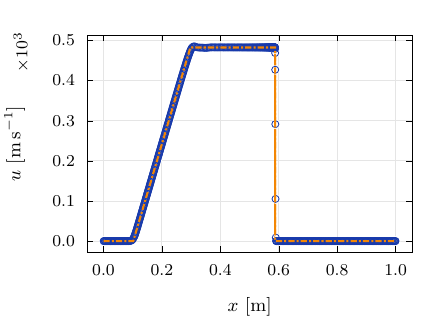}%
\includegraphics[width=0.33\textwidth]{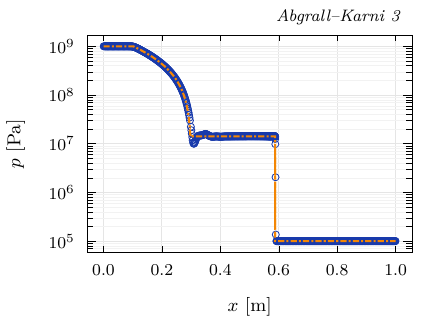}%
\caption{Numerical results for the multi-material Riemann problems of Abgrall and Karni, obtained
using the proposed implicit Lagrangian method with second order diagonally implicit \rk time integration (SDIRK2). 
The first two problems use a uniform-mass mesh of $N=1000$ cells, while the third, due to the large density ratio, 
adopts a smoothed quasi-uniform mesh in $\Delta x$, also composed of $N=1000$ cells. In the leftmost panel of the third row, a magnification window has been
added to better display the behaviour of the contact discontinuity found in the immediate vicinity of the shock front.}
\label{fig:rp5}
\end{figure}

The results of the tests labelled AK1 and AK2 use a uniform mesh of $N=1000$ cells (in mass space) are reported in Figure~\ref{fig:rp5}, 
which also contains the results of test AK3, where a quasi-uniform mesh of $N=1000$ cells (in physical space) is used. This choice was made
to provide a better distribution of the resolution of the scheme, since the density ratio of the test would yield a mesh with very few cells in the lighter 
right side, if a mass-uniform mesh were employed.

The results match the analytical solutions of the two-material 
Riemann problems and show limited spurious oscillations, 
despite the fact that the used CFL numbers always
exceed unity, and in particular the density and pressure oscillations (shown in log-scale in the third row of Figure~\ref{fig:rp5}) 
have acceptable amplitudes. In the last row of Figure~\ref{fig:rp5}, a magnification square has been added to better show 
small discontinous contact feature of the density profile in the vicinity of the shock front. The comparison with the analytical solutions does not show particularly severe issues.


\subsection{One-dimensional shock--bubble interaction}\label{sec:bubble}

\begin{figure}[!b]
\centering
\includegraphics[width=0.99\textwidth]{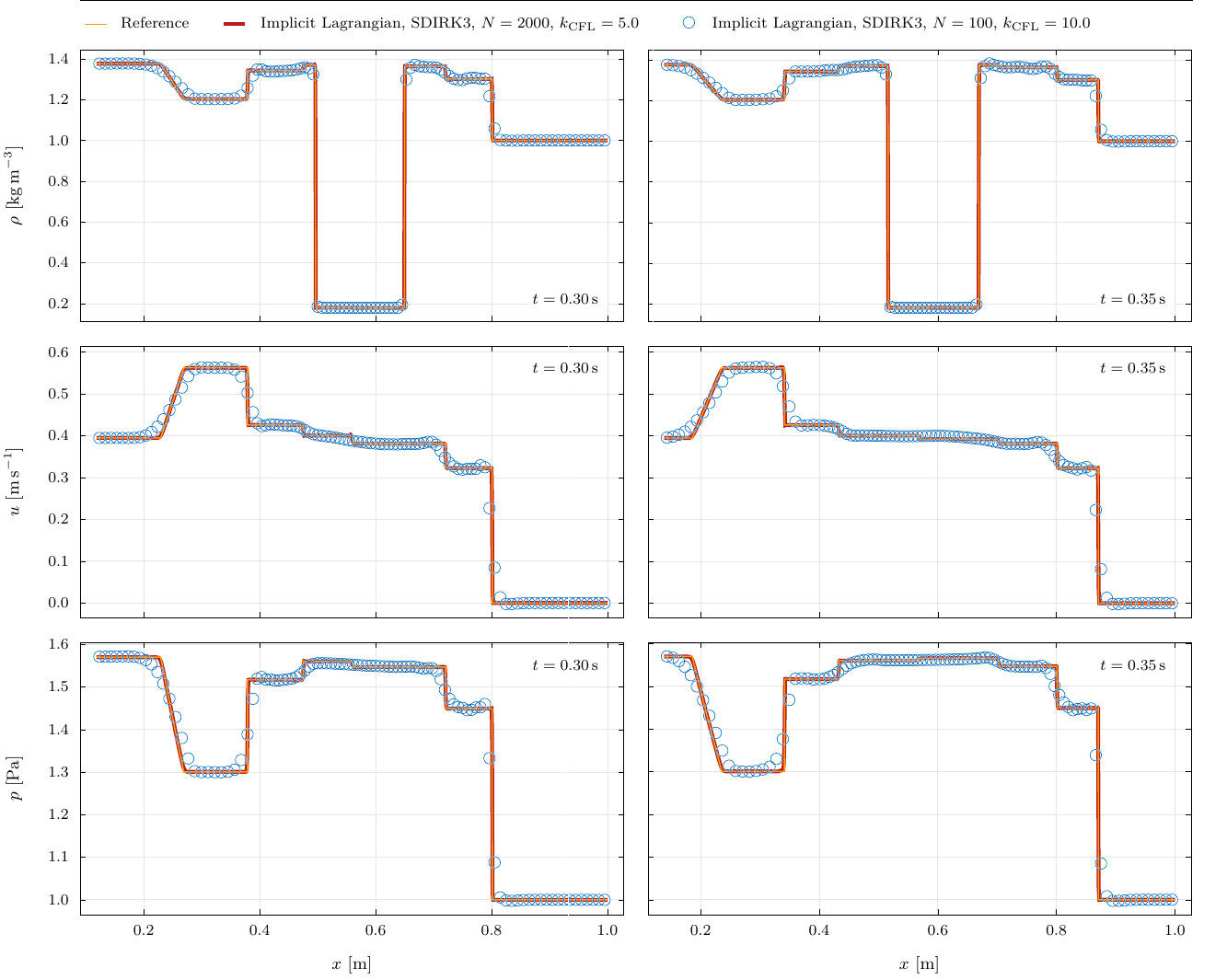}%
\caption{One-dimensional shock--bubble interaction problem of Quirk and Karni \cite{quirkkarni} at time $t=0.30\,\up{s}$ and 
$t=0.35\,\up{s}$. 
The solution obtained with the proposed implicit Lagrangian scheme 
and third order diagonally implicit \rk time-`stepping, using the second-order energy diffusion operator, 
on $N=2000$ uniformly spaced mesh points
(i.e. with non-uniform mass $\Delta m$) at $\cfl=5.0$, visually overlaps the reference given 
by an explicit Eulerian \muscl method
with 50000 uniform cells. No significant artefacts are visible in the under-resolved solution using 
$N=100$ points at $\cfl=10$, while adequately capturing all flow features.
 }
\label{fig:bubble}
\end{figure}

Further testing of the implicit Lagrangian scheme has been carried out
by reproducing the results of Quirk and Karni \cite{quirkkarni} 
on a one-dimensional shock--bubble interaction problem (Fig. \ref{fig:bubble}).
The solution obtained with the proposed implicit Lagrangian scheme 
and third order diagonally implicit \rk time-stepping, using the second-order energy diffusion operator, 
on $N=2000$ uniformly spaced mesh points
(i.e.\ with non-uniform mass $\Delta m$) at $\cfl=5.0$, visually overlaps the reference given 
by an explicit Eulerian \muscl method
with 50000 uniform cells. No significant artifacts are visible in the under-resolved solution using 
$N=100$ points at $\cfl=10$, while adequately capturing all flow features.

\begin{equation}
    \rho(t=0,\ x),\ u(t=0,\ x),\ p(t=0,\ x), \gamma(x) = 
    \left\{
    \begin{aligned}
        \rho_\up{L},\ u_\up{L},\ p_\up{L},\ \gamma_1 &\qquad\text{ if } x \le x_0,\\
        \rho_\up{R},\ u_\up{R},\ p_\up{R},\ \gamma_1 &\qquad\text{ if } x > x_0 \text{ and } x \le x_1,\\
        \rho_\up{B},\ u_\up{B},\ p_\up{B},\ \gamma_2 &\qquad\text{ if } x > x_1 \text{ and } x \le x_2,\\
        \rho_\up{R},\ u_\up{R},\ p_\up{R},\ \gamma_1 &\qquad\text{ if } x > x_2,\\
    \end{aligned}
    \right.
\end{equation}
The states are specified by
$\rho_\up{L} = 1.3765\,\up{kg\,m^{-3}}$, 
$\rho_\up{R} = 1.0\,\up{kg\,m^{-3}}$, 
$\rho_\up{B} = 0.138\,\up{kg\,m^{-3}}$, 
$u_\up{L} = 0.2948\,\up{m\,s^{-1}}$, 
$u_\up{R} = u_\up{B} = 0.0\,\up{m\,s^{-1}}$, 
$p_\up{L} = 1.57\,\up{Pa}$, and 
$p_\up{R} = p_\up{B} = 1.0\,\up{Pa}$.
The computational domain is delimited by $x_\up{L} = 0.0\,\up{m}$ and $x_\up{R} = 1.0\,\up{m}$.
The initial shock position is $x_0 = 0.25\,\up{m}$ and the initial position of the one-dimensional bubble
is delimited by $x_1 = 0.4\,\up{m}$ and $x_2 = 0.6\,\up{m}$.
The test adopts the ideal gas law with $\gamma_1 = 1.4$ and $\gamma_2 = 1.67$, hence 
in practice we also set $\Pi_1 = \Pi_2 = 0.0\,\up{Pa}$.

In Figure~\ref{fig:bubble} we show the solutions given by the proposed scheme with third order
SDIRK time-`stepping, at times $t = 0.30\,\up{s}$ and $t = 0.35\,\up{s}$.
The test is run on two meshes, composed of $N=2000$ or $N=100$ cells with constant $\Delta x = (x_\up{R} - x_\up{L})/N$ 
(and consequently discontinuous $\Delta m$), 
using $\cfl=5.0$ and $\cfl=10.0$ respectively.
The results match the reference solution given by an explicit, Eulerian, second order, path-conservative \cite{pares2006, castro2006}
Finite Volume MUSCL--Hancock method applied to the Kapila system, on a uniform mesh composed of 50\,000 control volumes.

\subsection{Application to stratified systems}\label{sec:multilayer}

\begin{table}[!t]
\caption{Initial conditions and material parameters for the stratified medium problems. The number of layer pairs $N_1$ and the final solution times $t_\up{end}$ are also reported.}
\label{tab:rpmultilayer}
\begin{tabularx}{\textwidth}{lRRRRRRRRRRR}
\toprule
    & $\rho_{1}$ & $\rho_{2}$ & $\gamma_{1}$ & $\gamma_{2}$ & $\Pi_{1}$         & $\Pi_{2}$ & $N_\up{l}$ & $t_\up{end}$ \\
\midrule
SM1 & $20.0$ & $10.0$ & $4.4$ & $1.4$ & $100.0$ & $0.0$ & $20$  & $2.0$\\
SM2 & $10^4$ & $10.0$ & $4.4$ & $1.4$ & $100.0$ & $0.0$ & $20$  & $2.5$\\
SM3 & $10^4$ & $10.0$ & $4.4$ & $1.4$ & $10^8$  & $0.0$ & $20$  & $1.0$\\
SM4 & $20.0$ & $10.0$ & $4.4$ & $1.4$ & $10^4$  & $0.0$ & $200$ & $2.0$\\
\bottomrule
\end{tabularx}
\end{table}




\begin{figure}[!b]
\centering
\includegraphics[width=0.99\textwidth]{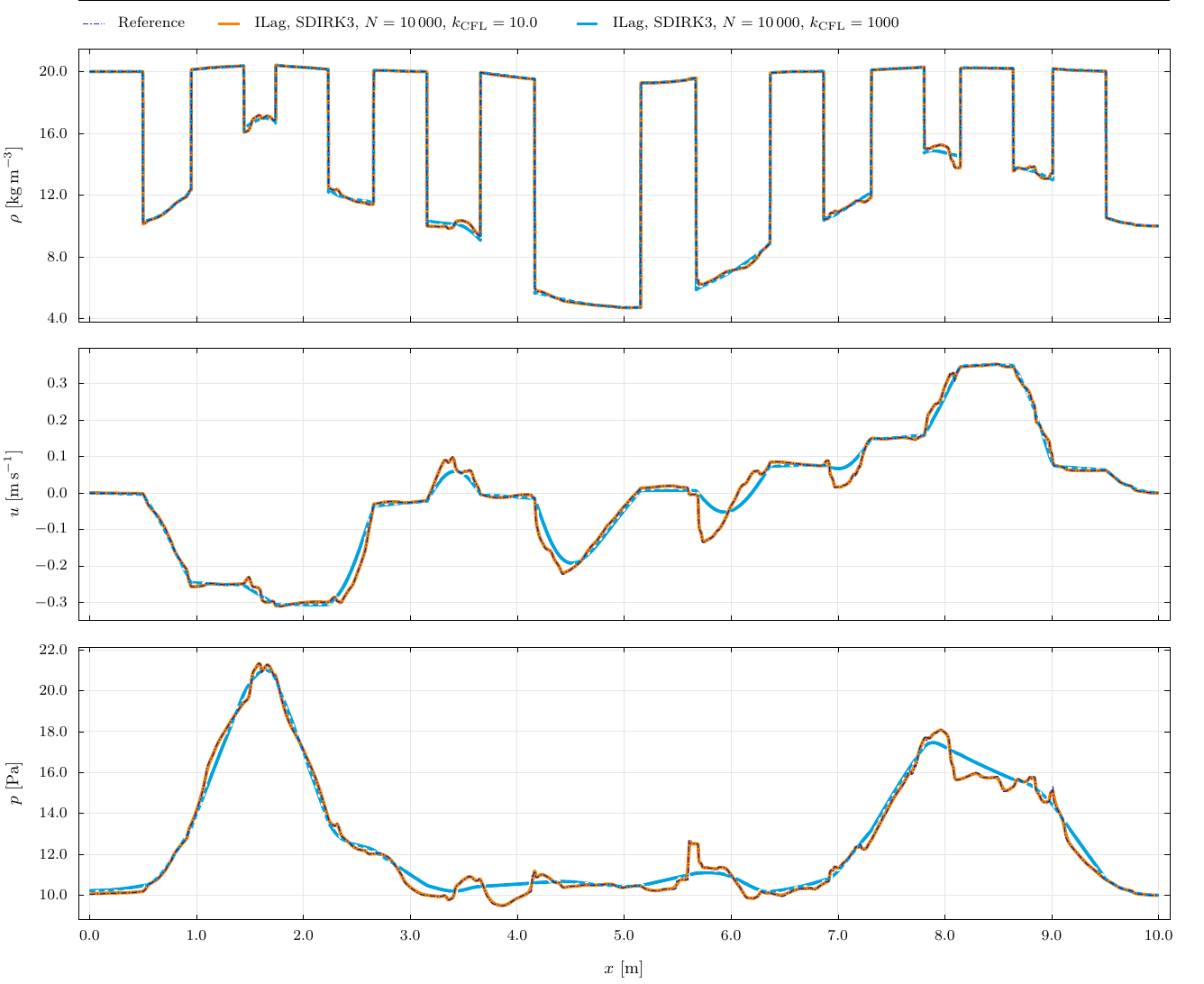}
\caption{Numerical results for the stratified medium problem SM1. For each flow variable, we
show the solutions given by the proposed implicit Lagrangian method on a grid of $N=10\,000$ cells
(500 cells per layer) at $\cfl=10.0$ and $\cfl=1000$.
The first solution visually coincides with the reference, and with higher CFL the main structures
are captured, while high frequency details are filtered out by the coarse time-stepping.
The reference solution has been obtained by solving the Kapila two-phase flow model model 
with an explicit Eulerian method, on a fine grid of half a million cells.}
\label{fig:multilayer1}
\end{figure}

\begin{figure}[!b]
\centering
\includegraphics[width=0.99\textwidth]{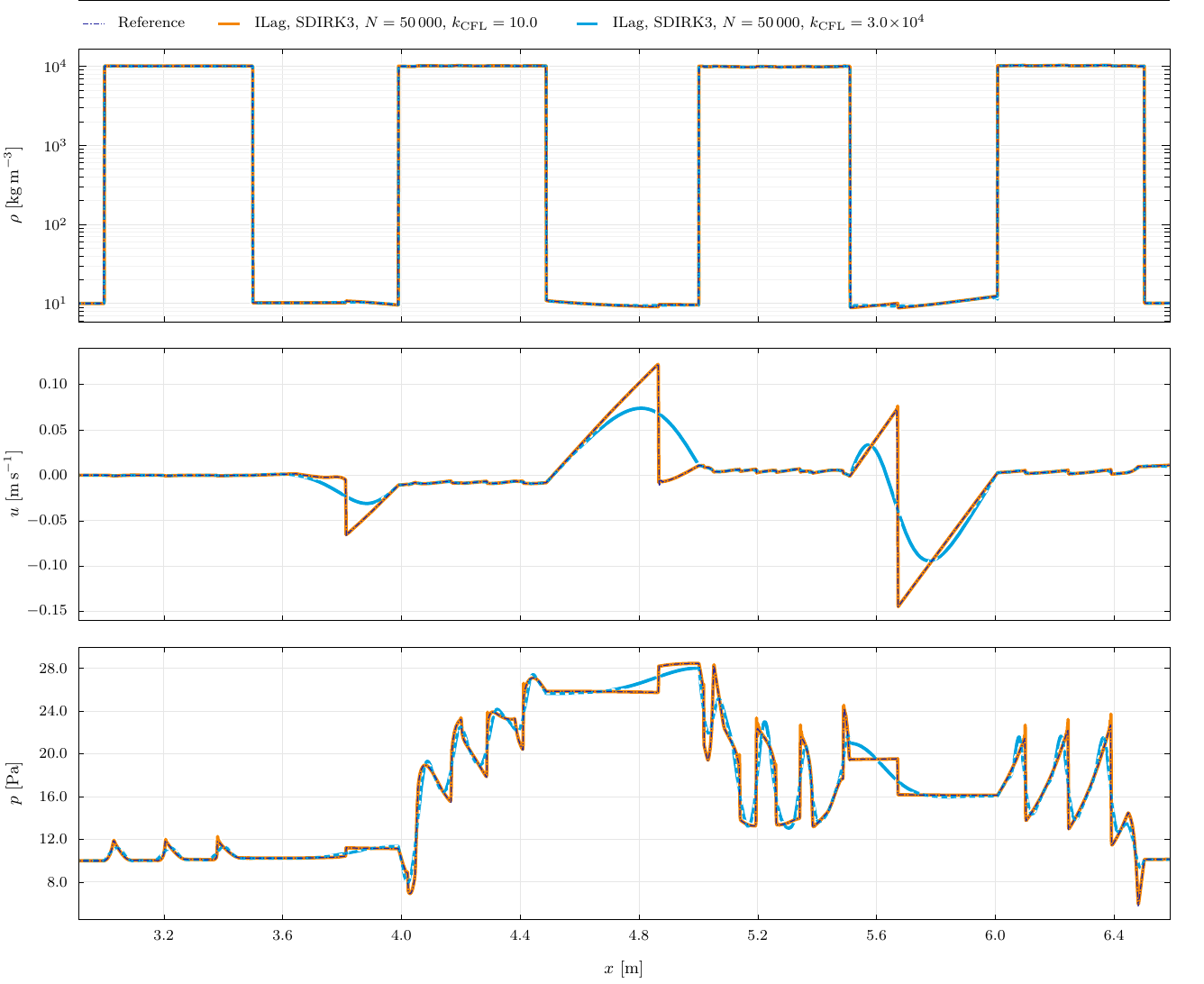}
\caption{Numerical results for the stratified medium problem SM2 (the density ratio is $\rho_1/\rho_2=1000$). 
For each flow variable, we
show the solutions given by the proposed implicit Lagrangian method on a grid of $N=50\,000$ cells
(2500 cells per layer) at $\cfl=10.0$ and $\cfl=30\,000$.
The first solution visually coincides with the reference, and with higher CFL the main flow structures
are captured without introducing instabilities.
The reference solution has been obtained by solving the Kapila two-phase flow model model 
with an explicit Eulerian method, on a fine grid of half a million cells.}
\label{fig:multilayer2}
\end{figure}

\begin{figure}[!b]
\centering
\includegraphics[width=0.99\textwidth]{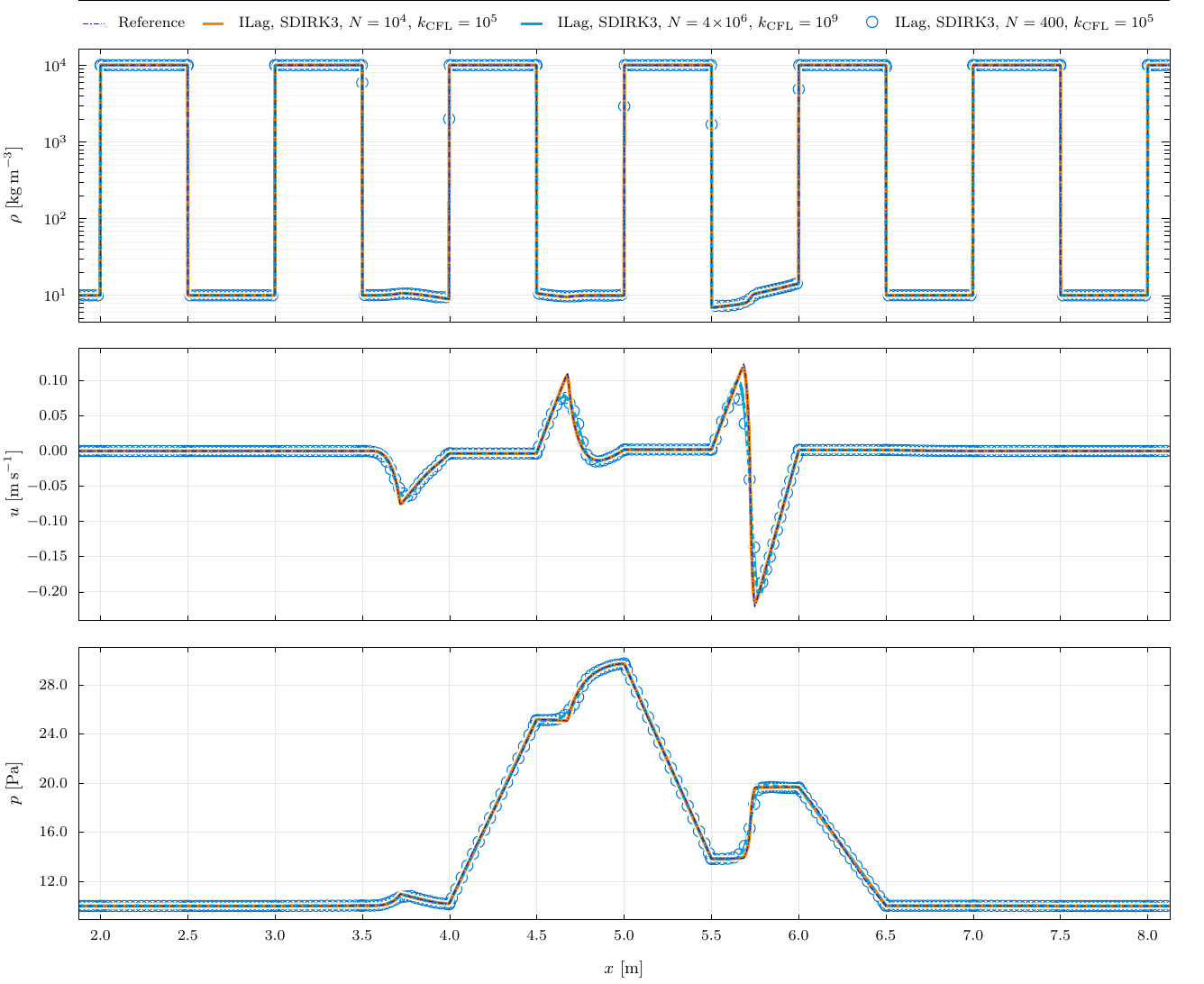}
\caption{Numerical results for the stratified medium problem SM3 (the density ratio is $\rho_1/\rho_2=1000$ and $\Pi_1=10^8$). 
For each flow variable, we show the solution given by the proposed implicit Lagrangian method on different grids and
at different CFL values. A well-resolved run employing $N=10\,000$ cells and $\cfl=100\,000$ is visually 
identical to the reference solution. A second run with four million cells and $\cfl=10^9$ shows
that if very fine meshes are to be used, then correspondingly high CFL values can be adopted. 
A third run with $N=400$ cells and $\cfl=10^5$ shows that the main flow features can be captured using modest grid sizes.}
\label{fig:multilayer3}
\end{figure}

\begin{figure}[!b]
\centering
\includegraphics[width=0.99\textwidth]{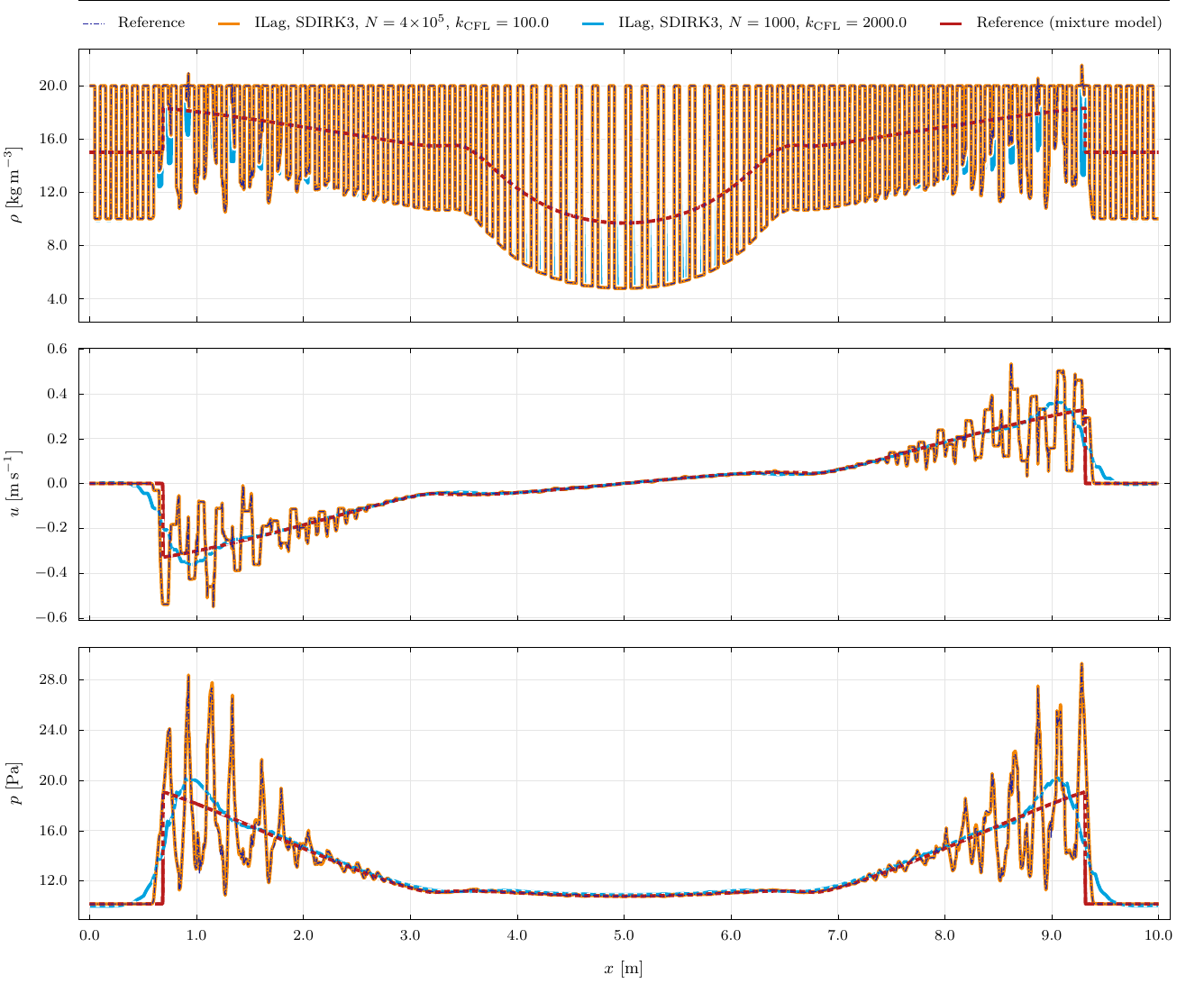}
\caption{Numerical results for the stratified medium problem SM4 (the layer count is $N_\up{l} = 200$). 
For each flow variable, we show the solution given by the proposed implicit Lagrangian method on different grids and
at different CFL values. A well-resolved run employing $N=400\,000$ cells (2000 per layer) and $\cfl=100.0$ matches 
the reference solution. A second run using $N=1000$ cells (5 per layer) and $\cfl=2000$ (that is, under-resolved 
in both space and time) yield a solution that tends towards that given by the multiphase model for the fluid mixture.}
\label{fig:multilayer4}
\end{figure}

The battery of tests is concluded by solving a series of applicative problems on multilayer systems, 
like in the paper by Phan et. al. \cite{phan2023numerical}, 
from which we derive the setup consisting of a set of alternating
layers of different density and material parameters, upon which an initially smooth pressure 
perturbation is allowed to propagate.

The alternating layers of different materials (characterized by 
different initial density and parameters of the equation of state), 
have constant length $L_\up{l} = (x_\up{R} - x_\up{L})/N_\up{l}$.
The leftmost layer will be assigned initial density $\rho_1$ and parameters of the equation of state $\gamma_1$ and $\Pi_1$, 
with the following layer switching to density $\rho_2$ and parameters $\gamma_2$ and $\Pi_2$, then back to $\rho_1$, $\gamma_1$ and $\Pi_1$, and so on.

Formally, we can set $c(x) = \up{mod}_2\left(1 + \lfloor x/L_\up{l} \rfloor \right)$, 
with $\up{mod}_2$ denoting the remainder operator for integer division by 2, 
and then define the initial data as 
\begin{equation} \label{eqn:multilayerrho}
        \rho(t=0,\ x) = c\,\rho_1 + \left(1 - c\right)\,\rho_2,\qquad
        \gamma(x) = c\,\gamma_1 + \left(1 - c\right)\,\gamma_2,\qquad
        \Pi(x) = c\,\Pi_1 + \left(1 - c\right)\,\Pi_2.
\end{equation}
Note that $c(x)$ can evaluate only to $c = 0$ or to $c = 1$, and thus no mixed states are produced in \eqref{eqn:multilayerrho}.
The fluid is initially at rest, that is, $u(t=0,\ x) = u_0 = 0.0\,\up{m\,s^{-1}}$, and
the pressure field is given by
\begin{equation}
p(t=0,\ x) = 
\left\{
    \begin{aligned}
        \left\{1 + 10\,\cos\left[2\,\pi\,\left(x - x_0\right)/\ell\right]\,\right\}\,p_0&\qquad\text{ if } \abs{x - x_0} < \ell/2,\\
        p_0 &\qquad\text{ if }\abs{x - x_0} \ge \ell/2,\\
    \end{aligned}
    \right.
\end{equation}
with $\ell = 2.5\,\up{m}$.
In this work we present several variants of this problem, differing from each other by the choice of initial densities $\rho_1$ and $\rho_2$, 
material parameters $\gamma_1$, $\gamma_2$,\ $\Pi_1$,\ $\Pi_2$ and number of layers $N_\up{l}$.
Table~\ref{tab:rpmultilayer} reports the specific parameters defining each one of the variants, labelled SM1, SM2, SM3, and SM4. 

\reviewera{While the numerical discretisation itself would allow for any distribution of layers (as in Section \ref{sec:bubble} for example), this setup with a constant layer width 
has been chosen as a matter of convenience when comparing the results with those of a homogenised model.}

The test labelled SM1 is a configuration like the one originally shown by Phan et al., with mild density and compressibility ratios, which 
nevertheless yield a sizeable mismatch of the wavespeeds in the two media, due to the Lagrangian framework.

Test SM2 introduces a new element of difficulty by increasing the density ratio from 2 to 1000. The lower density has 
been kept as $\rho_2=10.0$ while we increased $\rho_1$ to $\rho_1=10^4$, so to explore different regimes while changing only one parameter of the problem with 
respect to the original data given in \cite{phan2023numerical}, keeping in mind that the scaling of the parameters has to be considered with respect to each other.

A third variant, SM3 modifies one more parameter with respect to SM2, to obtain a density ratio of 1000 and at the same time increase the material pressure parameter $\Pi_1$ from $\Pi_1=100.0$ to $\Pi_1=10^8$, 
giving raise to a challenging situation in which the mismatch between the wavespeeds in the two media warrants the use of CFL numbers as high as $\cfl=10^9$.
Due to the density ratio and low compressibility of of the heavy phase, this test can be seen as a proxy for air-water configurations.

A final fourth alternative configuration, SM4 is a slightly stiffened version of SM1, where the pressure parameter $\Pi_1$ is increased to $\Pi_1=10^4$ and the number of layers is also increased to $N_1=200$. This variant is intended to show how the stratified media presented in this section behave if the number of layers is very high. 

The results in Figure~\ref{fig:multilayer1}, relative to the base configuration SM1, show that The solution given by the method using $N=10^4$ cells (500 per layer) and $\cfl=10$
matches the converged reference solution given by a semi-implicit Finite Volume solver for Kapila's reduced Baer--Nunziato model, and that increasing the timestep 
size by setting $\cfl=1000$ preserves the structure of the solution, while filtering out the higher frequencies, as ordinarily expected. 

In Figure~\ref{fig:multilayer2}, the same behaviour is observed for the problem with a higher density ratio SM2, the solution using $N=50000$ cells (2500 per layer) matches the converged reference from an independent solver/model and the increasing the timestep by setting $\cfl=3\times10^4$ maintains the structure of the solution, smoothing
out the sharper features.

In Figure~\ref{fig:multilayer3}, it is shown that the scheme can perform well, particularly in extremely stiff regimes, also if the number of cells is decreased. The solution 
using $N=10^4$ cells (500 per layer) matches the converged reference from the semi-implicit multiphase solver, and the same holds for a finer mesh of $N=10^6$ cells, using
$\cfl=10^9$. Additionally, a coarser simulation using $N=400$ total cells (20 per layer) and $\cfl=10^5$, can also capture the reference solution with 
relatively small deterioration in the accuracy of the results. This test represents a particularly challenging scenario for explicit solvers, due to heavy low-Mach effects and stringent timestep restrictions. The low-Mach or quasi-incompressible behaviour of the system can be clearly observed in Figure~\ref{fig:multilayer3}, which features linearly varying pressure in the heavy stiff medium, while slightly more complex profiles are observed in the lighter gas phase. The reference solution has been computed by means of a semi-implicit Eulerian, second order, path-conservative \cite{pares2006, castro2006} 
Finite Volume applied to Kaplia's multiphase flow model, specifically the one presented in \cite{chiocchettidumbser}.

In Figure~\ref{fig:multilayer4}, the results of the physical system SM4 and of the proposed numerical method show that systems with many layers, when dynamics are fully
resolved, exhibit an emergent oscillatory behaviour, which visually averages out to the physics captured by mixture models like the reduced Baer--Nunziato model by Kapila (red line), in which the fluid is modeled as a homogeneous one with an internal energy which is a volume-weighted average of the internal energy of the two adjacent layers, assuming pressure equilibrium of the two phases.
The same volume weighting defines the density of the mixture.
The propagation speed of the shock front is the same in both configurations, while the detailed solution of the stratified system also shows strong physical oscillations.
The physical nature of the oscillations and their accuracy has been verified by comparing the solution given by the proposed scheme using $N=4\times10^5$ cells (2000 per layer) and $\cfl=100$ with a reference obtained with a semi-implicit solver for the Kapila model (also resolving each layer separately) and $N=2\times10^6$ cells, showing that the solutions match, while using different model equations and different numerical methods. The reference scheme was also employed to obtain the reference solution 
for the system by treating the domain as if filled with a uniform mixture of the two fluids. As mentioned, this mixture solution can be observed to match the space-averaged behaviour of the stratified medium, and appears to be also captured by the proposed implicit Lagrangian method if the flow is under-resolved in space and time, here taking $N=1000$ cells (5 per layer) and $\cfl=2000$.

\section{Conclusions}
\label{sec:conclusions}
In this paper we introduced an original fully implicit finite volume scheme for the numerical solution 
of layered multifluid in Lagrangian mass coordinates. 

The method is based on a simple non-conservative (possibly oscillatory) predictor, and a 
quasi-conservative, non-oscillatory corrector, which is stable for arbitrarily large CFL numbers.
It is second order accurate in space, while it can be made arbitrarily high order accurate in time 
using SDIRK time integrators. We implement and test second and third order accurate schemes in time.

Given that the method is stable for large CFL numbers, the time step can be set based on just
accuracy requirements. This means that when solutions are expected to be sufficiently accurate, the
use of a third order method in time may be more efficient than a second order one, as it is shown in
the error-CPU plot (see Fig. \ref{fig:performance}).

Indeed, when using implicit time discretization, the time step restriction is no longer based on acoustic speeds. 
This means that some time step control would be useful to assess the optimal time step, which is therefore 
based on accuracy requirements,  rather than stability ones. Although several time-step controllers have 
been widely adopted for both explicit and implicit numerical solvers of systems of ordinary 
differential equations (see the 
books by Hairer and Wanner for example \cite{Hairer1993, Hairer1996}, as well as a recent application 
to RANS simulations with large timesteps \cite{massa2022}), it would be interesting to 
explore time step controllers suitable for the present implicit step, which makes use of the peculiar 
features of the novel scheme. This subject is currently under investigation.

The method is very flexible and robust, and is validated on a large battery of tests, which includes
all the problems in Toro's book on the Euler equations \cite{torobook}, the tests on multifluids by
Abgrall and Karni \cite{abgrall2001}, the tests on a finely stratified multilayer systems
\cite{phan2023numerical}, and other much harder variations on these, in which the two fluids have a
large ratio in the standard density and in the stiffness, such as, for example, a sequence of
several air-water pairs. 

With the proposed method it is possible to observe interesting physical phenomena. For example, in a
sequence of water-air pairs, with layers of approximately the same thickness, the internal energy
variation with respect to a rest configuration is larger in the air regions, as expected. 
Kinetic energy appears larger in the water layers, for smooth flow, while in presence of shocks the
kinetic energy may appear to be larger in the air layers. 

The efficiency and robustness of the method suggests several applications, which constitute current work in progress: 
\paragraph{Shock impact resilience} A screen made of composite material formed by several layers of two 
different materials may offer better protection from the impact of a projectile than a screen of a 
single material with the same thickness. The optimal distribution of layers that guarantees the best 
protection is still an open problem. The method developed in this paper may help designing optimal 
shields made out of suitable layers of two (or more) materials (see \cite{verreault2015optimization}). 
Note that to find the optimal (or a good) layering strategy is an optimization problem that may require 
several accurate solutions of the direct problem, and for this reason it is extremely important to 
rely on a fast and accurate method. 
\paragraph{Singularity formation in stratified material in slab geometry} In a remarkable paper by 
Zababakhin \cite{zababakhin1966shock} it is stated that it is possible to produce a singularity in 
the solution of compressible Euler equations, similar to the ones obtained with implosive shocks 
in a single gas in cylindrical or spherical symmetry, just in slab geometry, provided a shock 
propagates into suitable sequence of layers of two different fluids. The method developed in the 
current paper could be used to actually compute such kind of solution, and perhaps to formulate 
conjectures about possible layer distributions that lead to singularity formation in the multifluid system. 
\paragraph{Metamaterials composed by layers of several different constituents} In this paper we show 
that the numerical solutions obtained by our method for a stratified multifluid are in excellent 
agreement with a detailed numerical solution of the Kapila model in which the two phases are never 
really mixed (except numerically near an interface between the two fluids), so that in most of the 
computational domain the mass fractions of each fluid is either zero or one. Kapila model is more 
general than the multifluid solver developed in this paper for the treatment of fluid mixture, 
since it deals with intermediate states, however it is much less efficient for the specific case 
of immiscible multilayer. Furthermore, our method, like any Lagrangian method, can seamlessly deal 
with a multilayer system composed by an arbitrary number of different materials, while in its 
current form Kapila's model can handle only a mixture of two fluids. This opens enormous possibilities
in exploring nonlinear wave propagation in metamaterials composed of several different media, 
and in the numerical simulation of new layered multimaterials.

In conclusion, the new scheme opens new perspectives in the study of layered multimaterial, with promising applications in science and technology.

\section*{Acknowledgments}
G.~Russo and S.~Chiocchetti would like to thank the Italian Ministry of University and Research (MUR) to support this research with funds coming from 
PRIN Project 2022 (No. 2022KA3JBA entitled  ``Advanced numerical methods for time dependent parametric partial differential equations with applications''), and from the European Union’s NextGenerationUE – Project: Centro Nazionale HPC, Big Data e Quantum Computing, “Spoke 1” (CUP E63C22001000006).

S.~Chiocchetti acknowledges the support obtained by the Deutsche Forschungsgemeinschaft (DFG) via the project
DROPIT, grant no. GRK 2160/2, and from the European Union’s Horizon Europe Research and Innovation Programme under the
Marie Sk\l{}odowska-Curie Postdoctoral Fellowship MoMeNTUM (grant agreement No. 101109532).

S.~Chiocchetti would also like to thank M.~Dumbser for his encouragement and to acknowledge the fundamental insights he provided during
the early development of the scheme presented in this work.

\vspace{2ex}
\noindent\includegraphics[height=1.5cm]{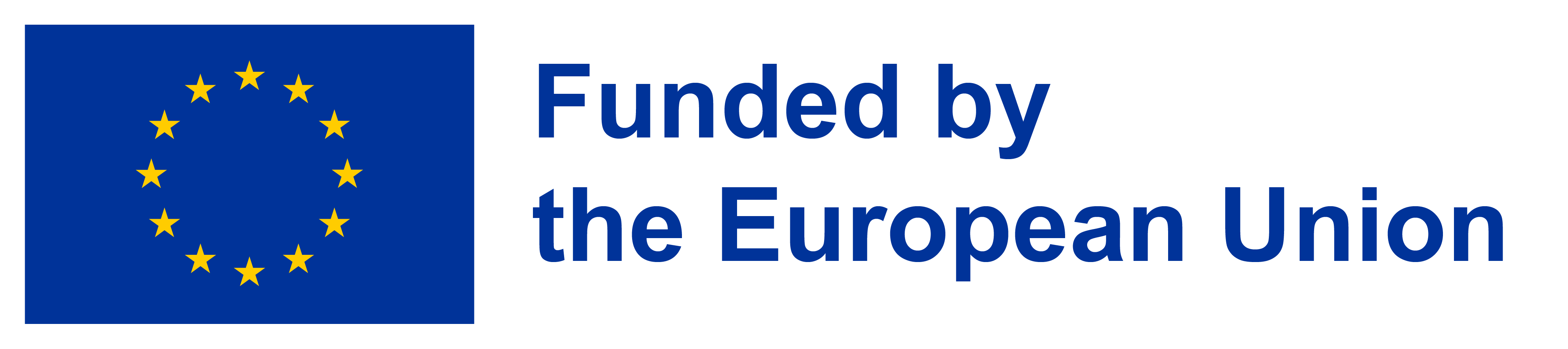}

\clearpage
\appendix
\section*{Appendices}
\addcontentsline{toc}{section}{Appendices}

\section{Additional Diffusion}
\label{appendix:extra_diffusion}
In this appendix, we report the details concerning an additional diffusion step 
applied selectively to address strong oscillations encountered in a two of the shock tube problems of Toro. 
Whenever such more aggressive diffusion is needed, for example in tests involving extremely strong shock waves (see Section~\ref{sec:rp}, problems \textit{Toro1} and \textit{Toro6} for the only two test cases in which this addition to the diffusion operator has been applied), then the target value $V_i^{n+1,\ast}$ for specific volume is further modified with a second filtering procedure.
To do so, we operate a second pass over the list of corrected specific volume values $V_i^{n+1,\ast}$, using 
a four cell window spanning from cell $i-1$ to cell $i+2$, and computing further corrected values
for $i$ and $i+1$, the two central cells of the window:
\begin{equation} \label{eqn:seconddiffusion}
\begin{aligned}
    &V_{i}^{n+1,\ast\ast} = (1 - w_i)\,V_i^{n+1,\ast} + w_i\,\left[V_{i-1}^{n+1,\ast} + \frac{m_i - m_{i-1}}{m_{i+2} - m_{i-1}}\,\left(V_{i+2}^{n+1,\ast} - V_{i-1}^{n+1,\ast}\right)\right],\\
    &V_{i+1}^{n+1,\ast\ast} = (1 - w_i)\,V_{i+1}^{n+1,\ast} + w_i\,\left[V_{i-1}^{n+1,\ast} + \frac{m_{i+1} - m_{i-1}}{m_{i+2} - m_{i-1}}\,\left(V_{i+2}^{n+1,\ast} - V_{i-1}^{n+1,\ast}\right)\right].
\end{aligned}
\end{equation}
Equation \eqref{eqn:seconddiffusion} states that these second-pass corrected values $V_{i}^{n+1,\ast\ast}$ and $V_{i+1}^{n+1,\ast\ast}$ are weighted averages of 1. the previous corrected
values $V_{i}^{n+1,\ast}$ and $V_{i+1}^{n+1,\ast}$ and 2. a simple linear interpolation between the two cell centers of the
outer cells of the four cell window.
The weight is computed as a function of the linear regression slopes $a_i$ (associated with cells $i-1$, $i$, and $i+1$) 
and $a_{i+1}$ (associated with cells $i$, $i+1$, $i+2$). Specifically we set
\begin{equation}
    w_i = \min\left(1,\ r_{i,i}\,r_{i,i+1} + \frac{2\,\abs{a_i - a_{i+1}}}{\abs{a_i} + \abs{a_{i+1}} + \epsilon}\right), \quad \epsilon = 10^{-14}
\end{equation}
with the linear regression coefficients being
\begin{equation}
    a_i = \frac{\sum_{k=-1}^1\left(m_{i+k} - m_{i,\up{m}}\right)\,\left(V_{i+k}^{n+1,\ast} - V_{i,\up{m}}^{n+1,\ast}\right)}{\sum_{k=-1}^1\left(m_{i+k} - m_{i,\up{m}}\right)^2},\qquad
    b_i = V_{i,\up{m}}^{n+1,\ast} - a_i\,m_{i,\up{m}},
\end{equation}
while having denoted
\begin{equation}
    m_{i,\up{m}} = \frac{1}{3}\sum_{k=-1}^1 m_{i+k},\qquad V_{i,\up{m}}^{n+1,\ast} = \frac{1}{3}\sum_{k=-1}^{1} V_{i+k}^{n+1,\ast},
\end{equation}
and having defined the normalised residuals of each one of the linear regression lines over three cells as
\begin{equation}
    r_{i,i} = \frac{\sum\limits_{k=-1}^1\abs{b_i + a_i\,\left(m_{i+k} - m_{i-1}\right) - V_{i+k}^{n+1,\ast}}}{\abs{V_{i+1}^{n+1,\ast} - V_{i-1}^{n+1,\ast}} + \epsilon}, \qquad
    r_{i,i+1} = \frac{\sum\limits_{k=0}^1\abs{b_{i+1} + a_{i+1}\,\left(m_{i+k} - m_i\right) - V_{i+k}^{n+1,\ast}}}{\abs{V_{i+2}^{n+1,\ast} - V_{i}^{n+1,\ast}} + \epsilon}.
\end{equation}

Each four-cell stencil spanning from cell $i-1$ to cell $i+2$ will generate a filtered value for both the
central cells $i$ and $i+1$. Since each cell is considered twice in this four-cell-stencil filtering procedure, 
once as the center-left cell and once as the center-right cell in the stencil, the specific volume $V_i^{n+1,\ast\ast}$ that
will be used to replace $V_i^{n+1,\ast}$ is taken to be the arithmetic average of the two.
We finally remark that this second step of artificial diffusion is in general not used for the results presented in 
this paper, unless as clearly indicated in two of the shock tube problems of Toro \cite{toro2009} in Section~\ref{sec:rp}.

\section{Automatic computation of the mesh parameter}
\label{app:mesh}

In this appendix we provide the details concerning the automatic procedure used for the selection of
the parameters adopted to construct the meshes used in this work.
These parameters can in principle be set manually, but removing such a need is beneficial
A manual selection of the mesh parameters $\beta$ and especially $z_0$ might be a bother for the user,
and the range of selectable values (constrained by positivity requirements on $\Delta m_i$) 
is not immediately apparent. For this reason we adopted a simple procedure that eliminates this issue and always yields valid
well-graded meshes.

The procedure starts by defining $\rho_\up{L}^\ast = \min(\rho_\up{L},\ \rho_\up{R})$ and 
$\rho_\up{R}^\ast = \max(\rho_\up{L},\ \rho_\up{R})$ so that the left layer is less dense than the right on
and aims at computing a value for the nondimensional parameter $\beta$ and a corresponding optimal value
for the shift $z_0$, 
such that, together with the mass conservation constraints \eqref{eq:meshmasscons}, the induced spacings 
$\Delta m_\up{L}$ and $\Delta m_\up{R}$ are both positive and their ratio $\Delta m_\up{L}/\Delta m_\up{R}$
is as similar as possible as the density ratio $\rho_\up{L}^\ast/\rho_\up{R}^\ast$.

Formally this corresponds to finding $\beta$ and $z_0$ that minimize a cost function of 
the type 
\begin{equation} \label{eq:meshquality}
    q = \abs{\frac{\Delta m_\up{L}\left(\beta,\ z_0\right)\,\rho_\up{R}^\ast}{\Delta m_\up{R}\left(\beta,\ z_0\right)\,\rho_\up{L}^\ast} - 1} + 
        \max\left\{0,\ H\,\up{sign}\left[\Delta m_\up{L}\left(\beta,\ z_0\right)\,\Delta m_\up{R}\left(\beta,\ z_0\right)\right]\right\},
\end{equation}
with $H = 10^{20}$ a penalty coefficient that excludes invalid meshes containing negative values for the spacing.
In practice, we find a value for $\beta$ in the interval between $\beta = 0$ 
(which corresponds to a uniform spacing mesh) and a larger initial value $\beta_0 = 1/4$ (corresponding to a rather gradual transition between the two meshes), 
starting from $\beta = \beta_0$ and bisecting the interval towards the largest admissible $\beta$.
At each iteration of the bisection procedure, for a fixed iteration value of $\beta$, the value 
of $z_0$ that minimizes \eqref{eq:meshquality} is found through a nested grid search.

At the end of the procedure, if $\rho_\up{L}^\ast$ and $\rho_\up{R}^\ast$ do not correspond to $\rho_\up{L}$ and $\rho_\up{R}$ respectively, 
that is, if they had to be switched, the sign of $z_0$ is switched, restoring the original orientation of 
the two layers. This allow to define a stopping criterion for the bisection process, saying that it
can terminate as soon as $q$ falls below a certain threshold $q_\up{f} = 0.25$, which states that the ratio between 
the mesh spacing $\Delta m_\up{L}/\Delta m_\up{R}$ and the density ratio $\rho_\up{L}^\ast/\rho_\up{R}^\ast$
are allowed to differed by $q_\up{f} = 0.25$ if the initial choice of $\beta$ was so large that
the bisection procedure had to be entered. For low to moderate density ratios, usually $\beta_0$ is
accepted immediately without requiring to be reduced by the bisection process.  

Once the parameters $z_0$ and $\beta$ are fixed, the formula for computing the final 
mesh spacing can be applied regardless of the orientation of the density jump, generating for any given number 
of cells $n$ in the half-layer pair, a well-graded variable spacing mesh that satisfies the mass constraints of the initial condition.
In Figure~\ref{fig:gradedmesh} we give an illustration of the procedure, showing a comparison
with the uniform Eulerian grid in the reference space $z$ and in the Eulerian coordinate system $x$ (in which 
the mesh appears completely different), together
with an example of how several half-layer pairs can be mirrored and concatenated to form the grid of 
a many-layer system. In a final note it must be specified that when concatenating many half-layer pairs, 
a single half-layers must be added at each one of the domain extremities: to do so, we simply replicate the mesh spacing 
of the last cell generated by the above procedure and then rescale each $\Delta m_i$, in the single half-layers at the 
extremities and in their directly neighbouring half-layers, by a constant factor to match the mass conservation constraints. The
effect of this rescaling is clear in the rightmost panel of Figure~\ref{fig:gradedmesh}.

\bibliographystyle{plain}
\bibliography{silag-bib}%
\end{document}